\documentclass[a4paper,10pt]{amsart}
\usepackage[foot]{amsaddr}
\usepackage[margin=2.7cm]{geometry}
\newtheorem{theorem}{Theorem}[section]

\theoremstyle{remark}
\newtheorem{remark}[theorem]{Remark}
\usepackage{systeme}
\usepackage{graphicx}
\usepackage[]{units}
\usepackage{color}
\usepackage{mathrsfs}
\usepackage{bm}
\usepackage[ansinew]{inputenc}
 \usepackage[breaklinks]{hyperref}
\usepackage{float}
\usepackage{amsmath}
\usepackage{isomath}
\usepackage{amsthm}
\usepackage{amssymb}

\usepackage{hyperref}

\newcommand{\blue}{\color{black}} %{\color{blue}}

\begin{document}
%\begin{frontmatter} 

\title[]{A Reduced Basis approach for PDEs on parametrized geometries based on the Shifted Boundary Finite Element Method {\blue and application to a Stokes flow}}

\author{Efthymios N. Karatzas\textsuperscript{1,*}}
\address{\textsuperscript{1}SISSA, International School for Advanced Studies, Mathematics Area, mathLab Trieste, Italy.}
\thanks{\textsuperscript{*}Corresponding Author.}
\email{efthymios.karatzas@sissa.it}

%    author two information
\author{Giovanni Stabile\textsuperscript{1}}%{1,**}}
%\thanks{\textsuperscript{**}Second Corresponding Author.}
\email{gstabile@sissa.it}

\author{Leo Nouveau \textsuperscript{2}}
\email{leo.nouveau@duke.edu}
\address{\textsuperscript{2}Civil and Environmental Engineering, Duke University, Durham, NC 27708, United States.}

\author{Guglielmo Scovazzi\textsuperscript{2}}
\email{guglielmo.scovazzi@duke.edu}

\author{Gianluigi Rozza\textsuperscript{1}}
\email{grozza@sissa.it}


\date{}

\dedicatory{}

\begin{abstract}
We propose a model order reduction technique integrating the Shifted Boundary Method (SBM) with a POD-Galerkin strategy. This approach allows {\blue to deal with} complex parametrized domains in an efficient and straightforward way.  
The impact of the proposed approach is threefold. 

First, problems involving parametrizations of complex geometrical shapes and/or large domain deformations can be efficiently solved at full-order by means of the SBM{. \blue This unfitted boundary method permits to avoid remeshing and the tedious handling of cut cells by introducing an approximate surrogate boundary.}

Second, the computational effort {\blue is reduced} by the development of a Reduced Order Model (ROM) technique based on a POD-Galerkin approach. 

Third, the SBM provides a smooth mapping from the true to the surrogate domain, and for this reason, the stability and performance of the reduced order basis are enhanced. This feature is the net result of the combination of the proposed ROM approach and the SBM.  Similarly, the combination of the SBM with a projection-based ROM gives the great advantage of an easy and fast to implement algorithm considering geometrical parametrization with large deformations. The transformation of each geometry to a reference geometry (morphing) is in fact not required. 

These combined advantages will allow the solution of PDE problems more efficiently. We illustrate the performance of this approach on a number of two-dimensional Stokes flow problems.
\end{abstract}
\maketitle
%\end{frontmatter}
\section{Introduction and Motivation}\label{sec:intro}
Dealing with the numerical simulation of problems characterized by deforming domains, it is mainly possible to rely on two different approaches: one based on an arbitrary Lagrangian-Eulerian (ALE) formulation \cite{HiAmCo74}, where the grid is deformed using a mesh motion algorithm, and one based on an embedded boundary approach, which relies on an undeformed background mesh into which the boundary is embedded or immersed (see \cite{MiIa05} and references within). In this work the attention is focused on the second type of approaches and in particular onto ROMs emanating from the recently proposed Shifted Boundary Method presented in \cite{MaSco17_2}.

In general, the family of immersed/embedded boundary methods can be very useful in many cases where conformal classic methods are inefficient. Using such methods, it is possible to avoid most of the issues related with mesh deformation and topological changes that normally occur when dealing with geometrical parametrization or fluid-structure interaction problems \cite{BaTaTe13,BaDaPeRo17,stabile_geo}. Starting with the pioneering work of Peskin \cite{PESKIN1972252}, great effort has been focused on embedded or immersed methods. 

However, for both immersed or conforming methods, there are still many cases where the solution of partial differential equations with standard discretization techniques (Finite Element Method, Finite Volume Method, Finite Difference Method) becomes unfeasible. Such situations occur, for example, whenever a large number of different system configurations need to be tested - as in uncertainty quantification, optimization, parametrization studies - or reduced computational times {\blue are required as} in real-time control problems. A possible way to overcome this limitation is the use of reduced order modeling techniques \cite{HeRoSta16,quarteroniRB2016,ChinestaEnc2017,BeOhPaRoUr17}. In particular we rely onto reduced basis methods. 

It is {\blue beyond} the scope of this paper to provide a comprehensive review on immersed/embedded methods and we report here only some of the most significant contributions on the topic. The reader interested {\blue in more detail} on the different existing approaches, namely Ghost-Cell finite difference methods, Cut-Cell finite volume approach, Immersed Interface, Ghost Fluid, Volume Penalty methods can refer to the review paper \cite{MiIa05} and references therein.

In \cite{MaSco17_2,MaSco17_3}, a new Embedded Boundary Method (EBM) - called Shifted Boundary Method (SBM) - was introduced and applied to heat transfer problems, advection-diffusion equations, Stokes flow and laminar and turbulent Navier-Stokes equations. The idea of the SBM is to shift the position of the boundaries from their true location to a surrogate one and to appropriately modify the value of boundary conditions on shifted boundaries, with the goal of simplicity, robustness and efficiency. The main characteristic of the SBM is that the surrogate boundary location is chosen so that cut cells are completely avoided while optimal convergence rates and low algebraic condition numbers are preserved. These features make the SBM a candidate for the efficient solution of fluid problems with parametrized geometries, fluid-structure interaction problems, evolutionary time systems, nonlinear problem cases, hyperbolic problems \cite{MaSco17_2,MaSco17_3,SoMaScoRi17}.  
 
The  overall objective of this paper is to investigate how the SBM may be employed to effectively and efficiently solve parameterized partial differential problems within the reduced basis (RB) context. To the best of our knowledge, the combination of the reduced basis method with full-order embedded boundary methods has not been investigated  with the aim of improving computational performance and allow offline-online computational strategies. We mention here the work of \cite{BaFa2014} for classical EBMs and model reduction, in which two regions are separated by an evolving interface and the snapshot compression problem is formulated as a binary weighted low-rank approximation.
In particular, in the present work, the attention is focused on reduced order methods on parametrized geometries. We remark that the use of an embedded method at full-order level allows us to avoid the need to map all the deformed configurations to reference domains (morphing), as it is traditionally done in systems with parametrized geometry, see e.g. \cite{HeRoSta16,RoVe07,Rozza2009,ballarin2015supremizer,RoHuMa13,Rozza2008229,BeOhPaRoUr17} and references therein. Mapping all the possible geometrical configurations to a reference domain, especially in cases with large deformations, may in fact produce highly distorted meshes and therefore lead to ill conditioned problems. 

The work is organized as follows: in \autoref{sec:HF} we introduce the formulation and the methods used for the full-order approximation of the equations. In \autoref{sec:ROM} the reduced order methodology is discussed in details, while in \autoref{sec:num_exp} {\blue the proposed ROM technique is tested on a numerical benchmark, dealing with the geometrical parametrized problem of the flow around a circular embedded domain.} Finally in \autoref{sec:conclusions} conclusions and perspectives for future improvements and developments are drawn.

\section{The physical problem and the full-order approximation} \label{sec:HF}
\subsection{Strong formulation of the steady Stokes problem} 
The Stokes equations describe the flow of a Newtonian, incompressible viscous fluid when the convective forces are negligible with respect to the viscous forces. Consider an open domain $\mathcal D$ in ${\mathbb R}^d$, with $d=2,3$ the number of space dimensions with Lipschitz boundary $\Gamma$ (split into two sub-boundaries $\Gamma_D$, $\Gamma_N$) and let a $k-$dimensional parameter space $\mathcal P$ with a parameter vector $\mu\in \mathcal P \subset \mathbb{R} ^k
$, the strong form of the stationary Stokes flow system of equations with Dirichlet and Neumman  boundary conditions on $\Gamma_D$ and $\Gamma_N$ respectively, geometrically parameterized by $\mu$, is given by:
\begin{eqnarray*}
-\nabla\cdot (2\nu{{\bm{\epsilon}} (\bm u(\mu))} - p(\mu){ \bm I}) &=& \bm g(\mu),\, \,\,\,\text{ in } {\mathcal{D}(\mu)}, 
\\
\nabla \cdot { \bm u(\mu)} &=& 0, \qquad \,\text{ in }{\mathcal{D}(\mu)}, 
\\ \bm u(\mu) &=& \bm g_D(\mu), \text{ on } \Gamma_D(\mu),
\\
(2\nu { \bm \epsilon ( \bm u(\mu))} - p(\mu){ \bm I}) \cdot { \bm n} &=& { \bm g_N(\mu)} , \text{ on } \Gamma _N(\mu).
\end{eqnarray*}
We denote by $\nu$ the viscosity, ${ \bm  \epsilon ( \bm u)} = 1/2( \nabla   {  \bm u}+ \nabla {  \bm u}^T )$ the velocity strain tensor (i.e., the symmetric gradient of the velocity $\bm u$), $p$ is the pressure, $ \bm g$ a body force, $ \bm  g_{D}$ the values of the velocity on the Dirichlet boundary and $  \bm g_{N}$ is the normal stress on the Neumman boundary. 
The first equation represents the conservation of the linear momentum of the fluid, while the second equation is the incompressibility condition and describes the mass conservation.
\subsection{Full-order parametrized Shifted Boundary Method formulation}
In the following, we recall the full-order Shifted Boundary Method (SBM) for the steady Stokes equations following \cite{MaSco17_2}. 
{\blue{The Shifted Boundary Method is an embedded (immersed, non-conformal) finite element method, which relies on a Nitsche-type approach for the weak imposition of  Dirichlet boundary conditions~\cite{MaSco17_2,MaSco17_3,SoMaScoRi17}. Weak boundary conditions require less complicated data structures with respect to strongly imposed boundary conditions, and, for this reason, the SBM relies on the Nitsche approach.}}
For the sake of simplicity, in this subsection, we will omit the parameter dependency with respect to $\mu$. We also denote $\upsilon :=\upsilon _h$. 

\begin{figure}
\centering
\begin{minipage}{\textwidth}
\centering
\includegraphics[width=0.2\textwidth]{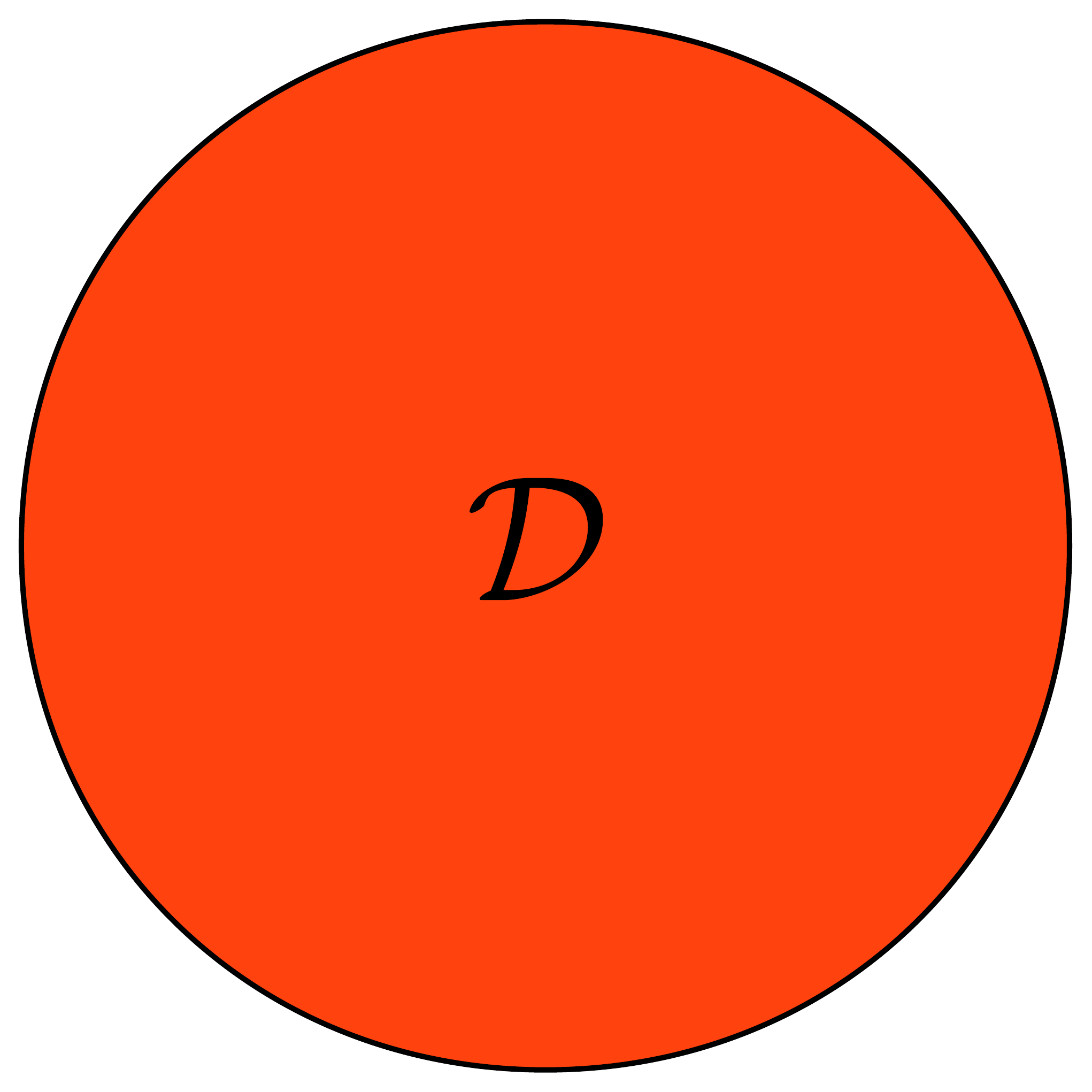}(a)
~~~~~~
\includegraphics[width=0.2\textwidth]{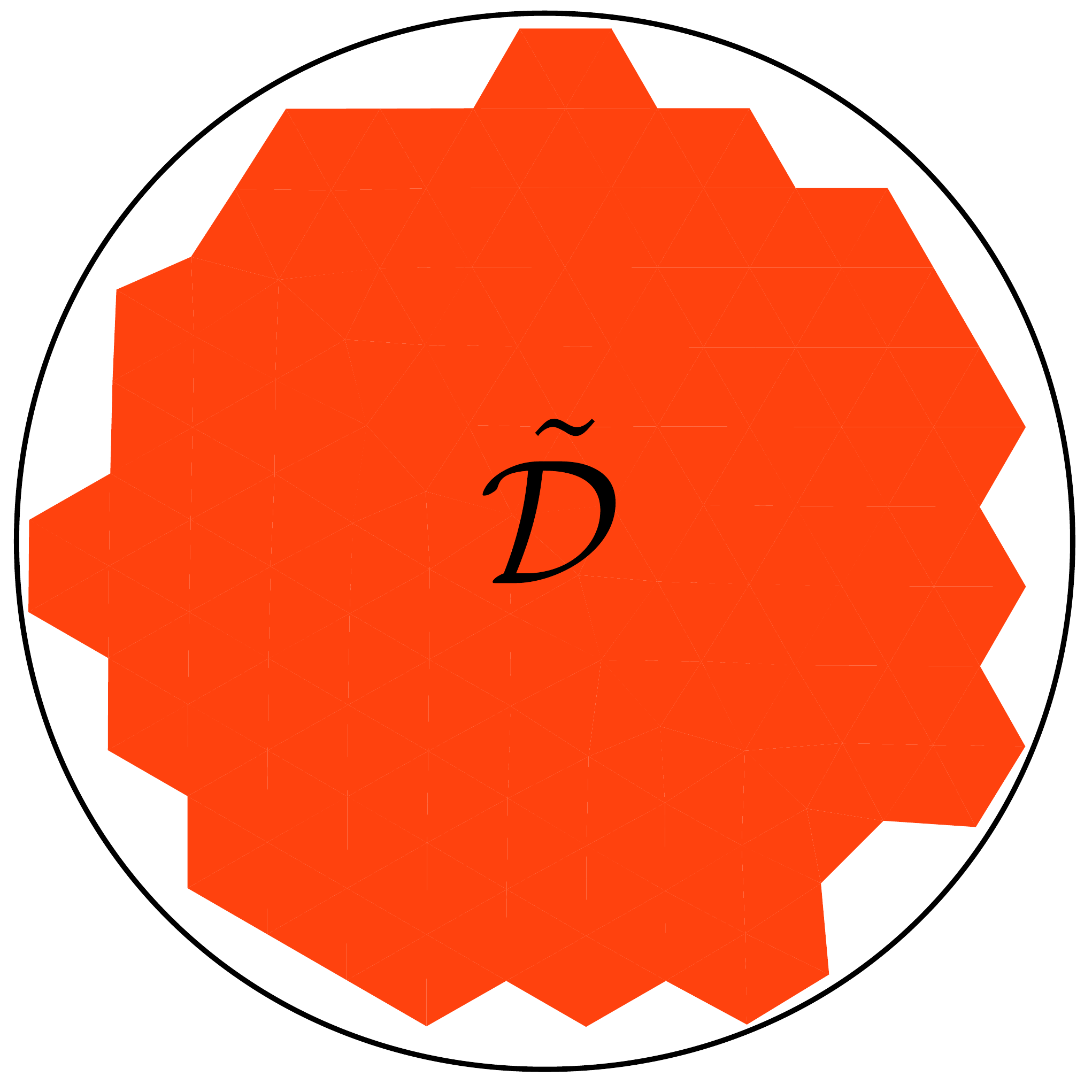}(b)
~~~~~~
\includegraphics[width=0.2\textwidth]{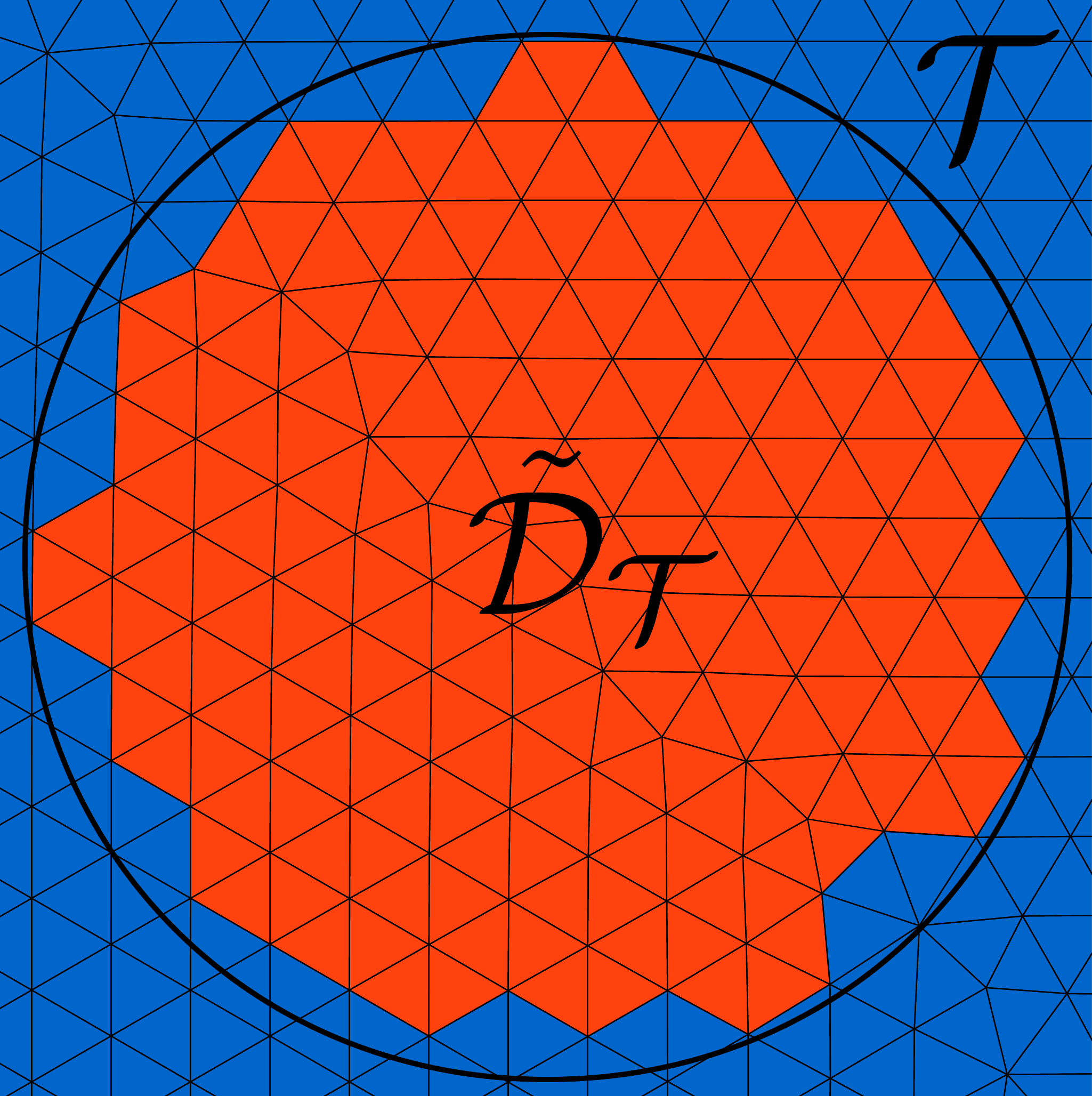}(c)\\[.5cm]
\includegraphics[width=0.27\textwidth]{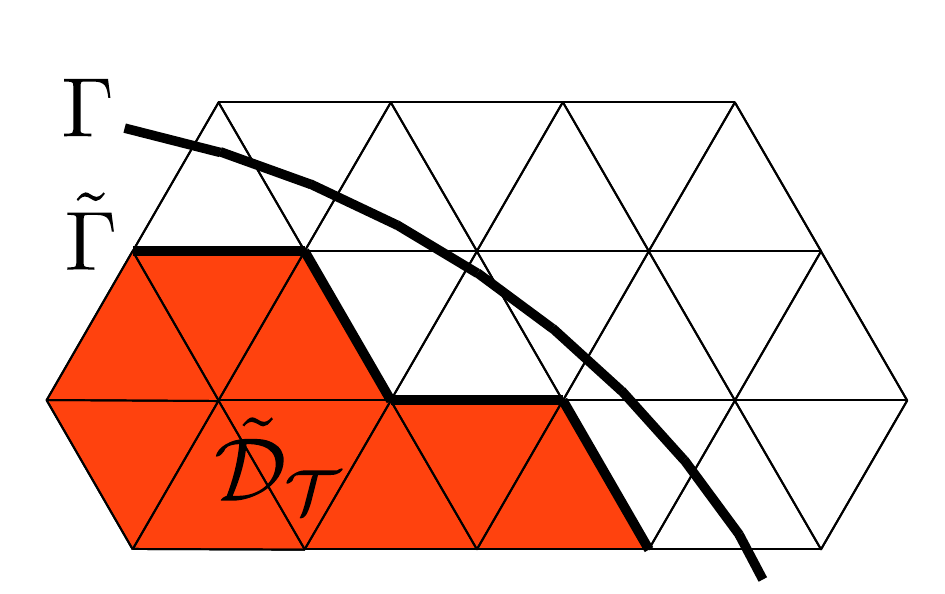}(d)
~~~~~~
\includegraphics[width=0.24\textwidth]{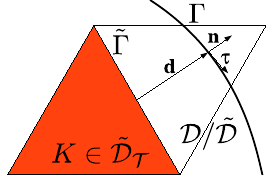}(e)
\end{minipage}
\caption{(a) The geometry of a disk, (b) the SBM surrogate geometry, (c) a zoom on the background mesh together with the surrogate SBM discretized geometry, (d)  SBM mesh and surrogate boundary  and (e) the normal and distance vector considering one element.}
\label{SurrogateMesh}
\end{figure}

We define now an approximate computational domain and an approximate boundary starting from the true boundary and the true computational domain, with the purpose of avoiding cut cells.
As seen in Figure~\ref{SurrogateMesh}, we denote by $\tilde  \Gamma$ the surrogate boundary composed of the edges/faces of the mesh that are the closest to the true boundary $\Gamma$. The closest faces/edges of $\tilde \Gamma$ to $\Gamma$ are detected using the closest-point projection algorithm.

The surrogate boundary $\tilde \Gamma$ encloses the surrogate domain $\tilde {\mathcal{D}} \subset \mathcal{D}$, where $\mathcal{D}$ is the original (true) computational domain. Furthermore, $\bm{\tilde{n}}$ indicates the unit outward-pointing normal to the surrogate boundary $\tilde \Gamma$, and it differs from the outward-pointing normal  $\bm  n$  of ${\Gamma}$. 
Notice also that the closest-point projection, in spite of the segmented/faceted nature of the surrogate boundary $\tilde \Gamma$, is actually a {\blue{continuous, piecewise-smooth}} mapping $\bm M$ from points in $\tilde \Gamma$ to points in $ \Gamma$:
\begin{eqnarray*}
{\bm{M}}: {\tilde{\bm x}} \in{\tilde\Gamma}\to {\bm x} \in {\Gamma}   .
\end{eqnarray*}
The mapping ${\bm{M}}$ defines a distance vector function
$$\bm d \equiv d_M ( {\bm x}) = {\bm x} - \tilde {\bm x} = [M - I]( \tilde {\bm x}), $$
where the distance vector $\bm d = ||\bm d||\bm n$ is aligned along $\bm n$, due to the closest point projection properties, see \autoref{SurrogateMesh}. 
{\blue{Since the distance vector $\boldsymbol{d}$ is aligned with the normal to the true surface $\Gamma$ and the true surface is smooth between edges and corners, it is legitimate to assume that ${\bm{M}}$ is continuous, and piecewise smooth.}}
The unit normal vector $\bm n$ and the unit tangential vectors $\bm \tau_i$   $(1 < i < d - 1)$ to the boundary $\Gamma$ can be easily extended to the boundary $\tilde \Gamma$ since
$\bm \bar {\bm n}(\tilde {\bm x}) \equiv {\bm n}(\bm M(\tilde {\bm x}))$,
 ${\bar{\bm\tau} _i({\bm x})} \equiv{\bm \tau}_i (\bm M(\tilde {\bm x}))$.
 In the following, whenever we write $\bm n(\tilde {\bm x})$ we actually mean $\bm \bar {\bm n}(\tilde {\bm x})$ at a point ${ \tilde {\bm x}}\in \tilde\Gamma$,
 and similarly for $ {\bm\tau} _i ( {\tilde {\bm x}})$ and $ \bar{\bm\tau}_i ( {\tilde {\bm x}})$. We can also introduce the derivative of a function $\bm g$ along the direction $\bar{\bm\tau _i}  $ at a point ${\tilde {\bm x}} \in \tilde\Gamma$ as $ \nabla _{\bar{\bm \tau} _i}\bar {\bm g} =\nabla \bar 
{\bm g}({\bm\tilde x}) \cdot \bar{{\bm\tau} _i}  ({\tilde {\bm x}})$.
The above constructions are the key ingredients in the extension of the boundary conditions on $\Gamma$ to the boundary $\tilde  \Gamma$ of the surrogate domain. 

Under the important assumption $\bm n\cdot \tilde {\bm n} \ge 0$, which relates to minimal grid resolution, we can now introduce the SBM variational formulation. 
Consider a discretization of the continuous boundary value problem with a mesh ${\mathcal{{\tilde{D}}_T}}$ consisting of simplices $K$ belonging to a tessellation $\mathcal T$. Moreover, we introduce  the discrete spaces ${\bm V}_h$ and $Q_h$, for velocity and pressure and we assume that a stable and convergent base formulation for the Stokes flow exists for these spaces in the case of conformal grids.
For the piecewise linear spaces
$$
\bm V_h = \left\{\bm v \in (C^0 ({\tilde{\mathcal{D}}}(\mu)))^{2} : \bm v|_K \in ({P}^1 (K))^{2}, \forall K \in {\mathcal{{\tilde{D}}_T}} (\mu)\right\},
 $$
$$
Q_h = \left\{ v \in C^0 ({\tilde{\mathcal{D}}(\mu)}) : v|_K \in  {P}^1 (K), \forall K \in {\mathcal{{\tilde{D}}_T}}(\mu)\right \},
$$
the stabilized formulation of \cite{Hu87} satisfies these assumptions, and will be used in what follows, although alternative choices are possible.  
Considering the standard notation $(\cdot, \cdot)_{\tilde{\mathcal{D}}}$, $\langle \cdot, \cdot \rangle _{\tilde\Gamma_D}$ , $\langle \cdot, \cdot \rangle _{\tilde\Gamma_N}$ for the $L^2(\tilde{\mathcal{D}})$, $L^2(\tilde\Gamma_D)$ and $L^2(\tilde\Gamma_N)$ inner products onto the surrogate geometry $\tilde{\mathcal{D}}$, $\tilde\Gamma_D$ and $\tilde\Gamma_N$, respectively, 
and that for every element $K\in {\mathcal T}$, we associate a parameter $h_{K}$, denoting the diameter of the set $K$, and the size of the mesh is denoted by $h= \max _{K\in{\mathcal T}} h_K$,
the SBM weak form reads: \\[.3cm]
{\it Find ${\bm u} \in  {\bm V}_h $ and $p \in Q_h $ such that, $\forall {\bm w} \in {\bm V}_h $ and $\forall q \in Q_h,$}
\begin{eqnarray}\label{eq:SBM}
\begin{cases}
a({\bm w},{\bm u};\mu)+b({\bm w},p;\mu)=\ell _g(\bm w;\mu),
\\
b({\bm u},q;\mu)+\hat b({\bm u},q;\mu) + (\mbox{stab. term.}) =\ell _q(q;\mu),
\end{cases}
\end{eqnarray}
where
\begin{eqnarray}
a({\bm w},{\bm u};\mu)
&=&
( {\bm \epsilon ( \bm w)}, 2\nu {\bm \epsilon (\bm u)})_{\tilde {\mathcal{D}}(\mu)} 
- \left\langle {\bm w} \otimes {\bm {\tilde n}}, 2\nu {\bm \epsilon (\bm u)}\right \rangle _{\tilde\Gamma _D(\mu)}
- \left\langle 2\nu {\bm \epsilon (\bm w)} , ({\bm u} + (\nabla {\bm u}){\bm d} ) \otimes {\bm{ \tilde n}}\right\rangle  _{\tilde\Gamma _D(\mu)} 
\nonumber\\
&&
+ \alpha\left \langle 2\nu/h({\bm w} + ( \nabla{\bm w}){\bm d}), {\bm u} + (\nabla  \bm u){\bm d} \right \rangle  _{\tilde\Gamma _D(\mu)}
+ \beta \left\langle 2\nu h \nabla _{\bar {\bm\tau} _i} {\bm w}, \nabla _{\bar {\bm\tau} _i} {\bm u}  \right\rangle_{\tilde\Gamma _D(\mu)} 
,
\\
b({\bm w},p;\mu)
&=&
 - (\nabla \cdot {\bm w}, p)_{\tilde {\mathcal{D}}(\mu)}+\left\langle  {\bm w} \cdot  {\bm{\tilde n}}, p\right \rangle _{\tilde\Gamma _D(\mu)}
 ,
 \\
 \ell _g(\bm w;\mu)
 &=&
 -({\bm w}, {\bm g})_{\tilde {\mathcal{D}}(\mu)}
+ \left\langle {\bm w}, {\bm g_N }  \right\rangle_{\tilde\Gamma _N(\mu)}
- \left\langle 2\nu {\bm \epsilon (\bm w)},   {\bm {\bar g_D}} \otimes {\bm {\tilde  n}}\right\rangle  _{\tilde\Gamma _D(\mu)} 
\nonumber \\
&&
+ \alpha\left \langle 2\nu/h({\bm w} + ( \nabla{\bm w}){\bm d}),  {\bm {\bar g_D}}\right \rangle  _{\tilde\Gamma _D(\mu)}
,
\\
 \hat b({\bm u},q;\mu)
 &=&
 \left\langle q {\bm d} \otimes {\bm{\tilde n}}, \nabla {\bm u}\right \rangle _{\tilde\Gamma _D(\mu)}
 ,
 \\
\ell _q(q;\mu)
&=&
\left\langle {\bm{\bar  g_D}}\cdot  \bm {\tilde{n}} ,q
\right \rangle _{\tilde\Gamma _D(\mu)} .
\end{eqnarray}
These equations yield the following algebraic system of equations:
\begin{equation}\label{eq:system_linear}
\begin{bmatrix} \bm{A}(\mu) & \bm{B}^T(\mu) \\ \bm{B}(\mu)+ {\hat{\bm{B}}(\mu)} & {\bm{C}}(\mu) \end{bmatrix} \begin{bmatrix} \bm{u}(\mu) \\ p(\mu) \end{bmatrix} = \begin{bmatrix} \bm{F}_g(\mu) \\ {\bm F}_q(\mu) \end{bmatrix}\mbox{.}
\end{equation}
In the above system of equations it is important to highlight several factors that will play a crucial role in \autoref{sec:ROM}. First of all, the discretized differential operators $\bm{A}$, $\bm{B}$ and $\bm{\hat{B}}$ are parameter dependent. Secondly, in the typical saddle point structure of the Stokes problem, the incompressibility equation is partially relaxed adding a stabilization term $\bm{C}$, {\blue derived according to~\cite{hughes1986new}}. This stabilization term, at full-order level, permits to circumvent the fulfillment of the ``inf-sup'' condition and the use of otherwise unstable pair of finite elements, such as $\mathbb{P}_1-\mathbb{P}_1$, and, as reported in \autoref{sec:ROM}, it helps to partially preserve the stability of the reduced order model. {\blue{In principle, the SBM method can also be combined with``inf-sup''-stable elements, such as the Taylor-Hood element: the $\mathbb{P}_1-\mathbb{P}_1$ was chosen for its simplicity.}}

The presented formulation is used to solve the full-order problem during an offline stage and to produce the snapshots necessary for the construction of the ROM of \autoref{sec:ROM}.

\section{Reduced order model with a POD-Galerkin method}\label{sec:ROM}
A Reduced Order Model (ROM) is a simplification of a Full Order Model (FOM) model that preserves its essential behavior and dominant effects with the purpose of reducing solution time or storage capacity. In this section, the projection-based ROM, that in this case relies on a POD-Galerkin approach \cite{HeRoSta16,Rozza2004ReducedBM}, is described and all the relevant issues are highlighted.

The interest is focused on partial differential equations with parametrized geometries that govern fluid dynamics problems and, in particular, on the advantages related with the use of the SBM. The attention is devoted on Reduced Basis (RB) ROM generated starting from high dimensional SBM approximations. Looking at pioneering works dealing with the RB method, starting from FEM full-order approximations, linear elliptic equations have been treated in \cite{Rozza2008229}, linear parabolic equations in \cite{grepl2005}, and non-linear problems in \cite{Veroy2003,Grepl2007}. In spite of the great number of works on reduced order modeling for fluid dynamics problems, to the best of the authors' knowledge, only very few research works can be found dealing with embedded boundary methods~\cite{BaFa2014}.

Before going into details, we remind the basics of the reduced basis method. The first step consists of the construction of a set of FOM solutions of the parametrized problem under the variation of the input parameters. The final goal of RB methods is to approximate any member of this solution set, {\blue and in general} of the solution manifold, with a reduced number of basis functions. The method consists of a two-stage procedure: the offline and the online one. During the  costly {\it{offline stage}}, one assembles the solution set and examines its components in order to construct a reduced basis that approximates any member of the solution set to a prescribed accuracy. This phase involves the solution of a possibly large number of FOM problems and the cost is usually high. During a second stage, namely the {\it{online stage}}, after the Galerkin/Petrov-Galerkin projection of the full-order differential operators describing the governing equations onto the reduced basis spaces, it is possible to solve a reduced problem for any new value of the input parameters. This offline-online procedure is  effective in many scenarios. We mention here the case when a large number of parameter values are in need of being tested with the consequent repeated evaluation of the system {\blue response, and} the case when a reduced computational time is required or only limited computational power and memory is available. 

We recall here that POD-Galerkin ROMs for the incompressible Stokes and Navier-Stokes equations suffer from stability issues~\cite{Caiazzo2014598,ballarin2015supremizer,Gerner2012,RoVe07,RoHuMa13,stabile_stabilized, StaHiMoLoRo17}, due to pressure instabilities, while for dynamic instabilities on transient problems see for instance \cite{Iollo2000,Akhtar2009,Bergmann2009516,Sirisup2005218,taddei2017}. More details will be given in \autoref{sec:stab_issues}.

From a reduced order modeling point of view, the first aim of this work is to investigate how ROMs are {\blue applied with an SBM approach}, and which are the perspectives. Since our main interest is to generate ROMs for flow problems on parametrized geometries {\blue we consider} the case of a Stokes flow problem. {\blue However, for the case of a Poisson problem on a parametrized geometry, which, originally, served as a simple feasibility test we refer to \cite{KaStaNoRoSco18_Heat}}. The SBM unfitted/surrogate mesh Nitsche finite element method is used to apply parametrization and reduced order techniques considering Dirichlet or Dirichlet combined with Neumann boundary conditions. 

An objective of this work is also to test the efficiency of a geometrically parametrized ROM method without the use of the transformation to reference domains and to highlight the advantages of having a fixed background mesh. Moreover, we will highlight and investigate how common strategies for pressure stabilization can be transferred to the present framework,  \cite{RoVe07,Quarteroni2007,ballarin2015supremizer}.

\subsection{The Proper Orthogonal Decomposition (POD)}\label{subsec_POD_theory}
In order to generate the reduced basis spaces, for the projection of the governing equations, one can find in the literature several techniques such as the Proper Orthogonal Decomposition (POD), the Proper Generalized Decomposition (PGD) and the Reduced Basis (RB) with a greedy sampling strategy. For more details about the different strategies the reader may see \cite{Rozza2008229,ChinestaEnc2017,Kalashnikova_ROMcomprohtua,quarteroniRB2016,Chinesta2011,Dumon20111387}. In this work a POD strategy is exploited and is chosen to apply the POD onto the full snapshots matrices that include parameter dependency. The full-order model is solved for each $\mu^k \in \mathcal{K}=\{ \mu^1, \dots, \mu^{N_k}\} \subset \mathcal{P}$ where $\mathcal{K}$ is a finite dimensional training set of samples chosen inside the parameter space $\mathcal{P}$. 
The  number of snapshots is denoted by $N_s $ and the number of degrees of freedom for the discrete full-order solution by $N_u^h$,  $N_p^h$ for the velocity and pressure, respectively.  The snapshots matrices $\bm{\mathcal{S}_u}$ and $\bm{\mathcal{S}_p}$, for velocity and pressure, are then given by $N_s$ full-order snapshots:
\begin{gather}
\bm{\mathcal{S}_u} = [\bm{u}(\mu^1),\dots,\bm{u}(\mu^{N_s})] \in \mathbb{R}^{N_u^h\times N_s},\quad\bm{\mathcal{S}_p} = [p(\mu^1),\dots,p(\mu^{N_s})] \in \mathbb{R}^{N_p^h\times N_s}.
\end{gather}
Given a general scalar or vectorial function $\bm{u}:{\mathcal D} \to \mathbb{R}^d$, with a certain number of realizations $\bm{u}_1,\dots, \bm{u}_{N_s}$, the POD problem consists in finding, for each value of the dimension of POD space $N_{POD} = 1,\dots,N_s$, the scalar coefficients $a_1^1,\dots,a_1^{N_s},\dots,a_{N_s}^1,\dots,a_{N_s}^{N_s}$ and functions $\bm{\varphi}_1,\dots,\bm{\varphi}_{N_s}$ that minimize the quantity:
\begin{eqnarray}\label{eq:pod_energy}
E_{N_{POD}} = \sum_{i=1}^{N_s}||{\bm{u}}_i-\sum_{k=1}^{N_{POD}}a_i^k {\bm{\varphi}}_k||^2_{L^2({\mathcal D})}, \,&&\forall
N_{POD} = 1,\dots,N,\\\nonumber
&& \mbox{ with } ({\bm{\varphi}_i,\bm{\varphi}_j}) _{{\mathcal D }} = \delta_{ij}, \mbox{\hspace{0.5cm}}  \mbox{ }\forall
 i,j = 1,\dots,N_s .
\end{eqnarray}
In this case the velocity field $\bm{u}$ is used as example. It can be shown \cite{Kunisch2002492} that the minimization problem of equation~\eqref{eq:pod_energy} is equivalent of solving the following eigenvalue problem:
\begin{gather}
{\bm {\mathcal{C}}^u}\bm{Q}^u = \bm{Q^u}\bm{\lambda^u} ,\quad \mbox{\hspace{0.5cm} for }\mathcal{C}^u_{ij} = ({\bm{u}_i,\bm{u}_j}) _{{\mathcal D }} \mbox{,\, } i,j = 1,\dots,N_s ,\nonumber
\end{gather}
where $\bm{\mathcal{C}^u}$ is the correlation matrix obtained starting from the snapshots $\bm{\mathcal{S}_u}$, $\bm{Q^u}$ is a square matrix of eigenvectors and $\bm{\lambda^u}$ is a diagonal matrix of eigenvalues. 

The basis functions can then be obtained with: 
\begin{equation}
\bm{\varphi_i} = \frac{1}{N_s{\lambda_{ii}^u}^{1/2}}\sum_{j=1}^{N_s} \bm{u}_j Q^u_{ij}.
\end{equation}
The same procedure can be repeated also for the pressure field considering the snapshots matrix consisting of the snapshots $
p_1,p_2,\dots,p_{N_s}$. One can compute the correlation matrix of the pressure field snapshots $\bm{C^p}$ and solve a similar eigenvalue problem ${\bm{C}^p}{\bm{Q}^p}={\bm{Q}^p}{\bm{\lambda}^p}$. The POD modes ${\chi_i}$ for the pressure field can be computed with:
\begin{equation}
 {\chi_i}=\frac{1}{N_s{{\lambda^p _{ii}}^{1/2}}}\sum_{j=1}^{N_s} p_j {Q^p_{ij}}.
\end{equation}
The POD spaces are constructed for both velocity and pressure using the aforementioned methodology resulting in the spaces:
\begin{equation}
\begin{split}
\bm{L}_u = [{\bm{\varphi}}_1, \dots , {\bm{\varphi}}_{N_u^r}] \in \mathbb{R}^{N_{u}^h \times N_u^r},\quad%\\
\bm{L}_p = [{\chi_1}, \dots , {\chi_{N_p^r}}] \in \mathbb{R}^{N_{p}^h \times N_p^r}.
\end{split}
\end{equation}
where $N_u^r$, $N_p^r < N_s$ are chosen according to the eigenvalue decay of $\bm{\lambda}^u$ and $\bm{\lambda}^p$, \cite{Rozza2008229,BeOhPaRoUr17}.

\begin{figure} 
\centering
\includegraphics[width=0.4\textwidth]{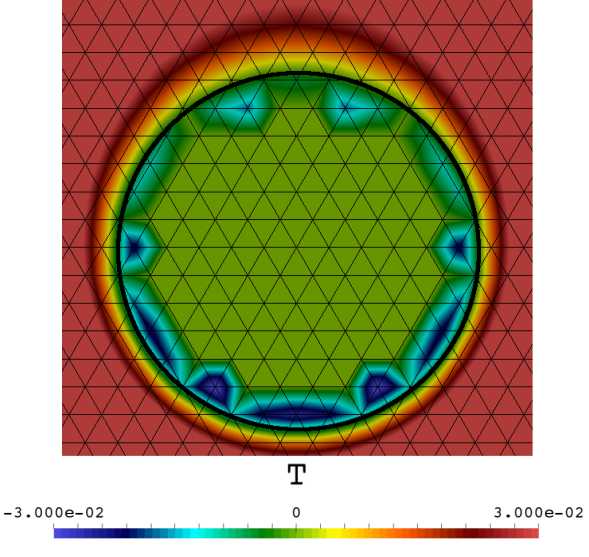}
\caption{A zoom onto the cylinder % of the Poisson numerical example 
in order to show the smoothing procedure employed by the SBM method inside the ghost area.}
\label{fig:poisson_zoom}
\end{figure}

\begin{remark}
The construction of the reduced order basis is based on the whole background domain. For this reason, the manipulation  of the out of interest - outside - the true geometry area, namely ``ghost area", needs particular care. In \cite{BaFa2014} the ghost area solution is set to zero. In the present work we use the solution values as they are computed {\blue applying the shifted boundary method and} the smooth map $\bm M$ from the true to the surrogate domain. This allows a smooth extension of the solution to the neighboring ghost elements with values which are decreasing smoothly to zero, see for instance the zoomed image in \autoref{fig:poisson_zoom}. This approach guarantees a regular ``solution'' in the background domain and { \blue permits the} construction of a reduced basis with better approximation properties. For more details and a full investigation of the possible choices and of the handling of the ghost area we refer to \cite{KaBaRO18}.
\end{remark}

\subsection{The geometrical parametrization: main differences with respect to a reference domain approach}
The standard procedure to deal with geometrical parametrization in a reduced order modeling framework is to map all the deformed configurations $\mathcal{D}(\mu)$ to a fixed reference domain $\mathcal{D}$ using a map $\mathcal{M}(\mu):\mathcal{D}\to\mathcal{D}(\mu)$ see e.g. \cite{HeRoSta16,RoVe07,Rozza2009,ballarin2015supremizer,RoHuMa13,Rozza2008229,BeOhPaRoUr17}. For the case at hand, considering a standard Galerkin formulation, without taking into consideration the stabilization terms, the weak formulation written on a reference domain reads:
\begin{eqnarray}
\begin{cases}
\tilde a({\bm w},{\bm u};\mu)+\tilde b({\bm w},p;\mu)=\tilde \ell _g(\bm w;\mu)
\\
\tilde b({\bm u},q;\mu) =0
\end{cases},
\end{eqnarray}
where the linear and bilinear forms are now written into a common reference domain:
\begin{align}
&\tilde a({\bm w},{\bm u};\mu) = \int_{\mathcal{D}} \bm{\epsilon}(\bm{w}) 2 \nu (J_T(\mu) )^{-1}(J_T(\mu) )^{-T}|J_T(\mu)| \bm{\epsilon}(\bm{u}) \mbox{d} \bm{x} , 
\\
&\tilde b({\bm{w},p;\mu}) = - \int_D \nabla \cdot \bm{w} (J_T(\mu) )^{-1}|J_T(\mu)|  p \mbox{d} \bm{x}, 
\\&
\tilde\ell _g(\bm w;\mu) = \int_\mathcal{D} |J_T(\mu)| \bm{g} \cdot \bm{w} \mbox{d} \bm{x},
\end{align}
and $J_T$ is the Jacobian of the map $\mathcal{M}(\mu)$ and $|J_T(\mu)|$ is its determinant. Moreover, in cases where it is possible to obtain a disjoint decomposition of the domain $\{\mathcal{D}_r\}_{r=1}^R$ such that for each $r \in \{1, \dots, R\}$ the map $\mathcal{M}_r(\mu):\mathcal{D}_r\to\mathcal{D}_r(\mu)$ is affine, it is possible to rewrite the above integrals as a sum of integrals over each subdomain and to take out from the integral the parametric dependent part. In this way it is possible to have an efficient offline-online splitting. However, especially for complex geometrical deformation, the operation of finding such decomposition is not trivial and may lead to a large number of subdomains. Moreover, in cases with large geometrical deformations, writing the equations on a reference domain may lead to distorted elements and therefore to numerical instabilities.

In the proposed approach, because we are always working onto the same physical domain, there is no need to introduce the change of variables in order to map the integrals to a common reference domain. The linear and bilinear forms of \autoref{eq:SBM} are ``naturally'' parametrized through the surrogate boundary $\tilde{\Gamma}(\mu)$ and the surrogate domain $\tilde{\mathcal{D}}(\mu)$. 

{\blue Our method, similarly to the reference domain approach, does not permit to avoid the non-affinity issue of the discretized differential operators. However, in order to obtain an efficient offline-online spliting, it is possible to rely on hyper reduction techniques such as the empirical interpolation method \cite{BARRAULT2004667}, the GNAT \cite{Carlberg2013623} or the Gappy-POD \cite{Everson1995}.}

\subsection{The projection stage and the generation of the ROM}
Once the POD functional spaces are set, the reduced velocity and pressure fields can be approximated with: 
\begin{equation}\label{eq:aprox_fields}
\bm{u^r} \approx \sum_{i=1}^{N_u^r} a_i(\mu) \bm{\varphi}_i(\bm{x}) = \bm{L}_u \bm{a}(\mu), \mbox{\hspace{0.5 cm}}
p^r\approx \sum_{i=1}^{N_p^r} b_i(\mu)\chi_i(\bm{x}) = \bm{L}_p \bm{b}(\mu).
\end{equation} 
The reduced solution vectors $\bm{a} \in \mathbb{R}^{N^r_u\times1}$ and $\bm{b}\in \mathbb{R}^{N^r_p\times1}$ depend only on the parameter values and the basis functions $\bm{\varphi}_i$ and ${\chi}_i$ depend only on the physical space. The unknown vectors of coefficients $\bm{a}$ and $\bm{b}$ can be obtained through a Galerkin projection of the full-order system of equations onto the POD reduced basis spaces and with the resolution of a consequent reduced algebraic {\blue system:}
\begin{equation}\label{eq:system_linear_proj}
\begin{bmatrix} \bm{L}_u & \bm{0} \\ \bm{0} & \bm{L}_p \end{bmatrix} \begin{bmatrix} \bm{A}(\mu) & \bm{B}^T(\mu) \\ \bm{B}(\mu) + \bm{\hat{B}}(\mu) & {\bm{C}}(\mu) \end{bmatrix} \begin{bmatrix} \bm{L}^T_u & \bm{0} \\ \bm{0} & \bm{L}^T_p \end{bmatrix} \begin{bmatrix} \bm{u}(\mu) \\ p(\mu) \end{bmatrix} = \begin{bmatrix} \bm{L}^T_u & \bm{0} \\ \bm{0} & \bm{L}^T_p \end{bmatrix} \begin{bmatrix} \bm{F}_g(\mu) \\ \bm{F}_q(\mu) \end{bmatrix}\mbox{,}
\end{equation}
which leads to the following algebraic reduced system:
\begin{equation}\label{eq:system_linear_reduced}
\begin{bmatrix} \bm{A}^r(\mu) & {\bm{B}^r}^T(\mu) \\ \bm{B}^r(\mu) + \bm{\hat{B}}^r(\mu) & {\bm{C}}^r (\mu) \end{bmatrix} \begin{bmatrix} \bm{a} \\ \bm{b} \end{bmatrix} = \begin{bmatrix} \bm{F}_g^r(\mu) \\ {\bm{F}^r_q(\mu)} \end{bmatrix}\mbox{,}
\end{equation}
where $\bm{A}^r(\mu) \in N_u^r \times N_u^r$, $\bm{B}^r(\mu), \bm{\hat{B}}^r(\mu) \in N_p^r \times N_u^r$, ${\bm{C}}^r(\mu) \in N_p^r \times N_p^r$ and $\bm{F}_g^r(\mu) \in N_u^r \times 1$, $\bm{F}^r_q(\mu) \in N_p^r \times 1$ are the reduced discretized operators and reduced forcing vectors respectively. The dimension of the reduced operators, as seen also in the numerical examples, is usually much smaller with respect to the dimension of the full-order system of equations and therefore less expensive to solve. We remark here that the full-order discretized differential operators that appear in \autoref{eq:system_linear} are parameter dependent and therefore, also at the reduced order level, we need to assemble the FOM problem in order to compute the reduced differential operator. Possible ways to avoid such potentially expensive operation, relying on an approximate affine approximation of the full-order differential operator, could be to use hyper reduction techniques,  \cite{Xiao20141,BARRAULT2004667,Carlberg2013623}. In this work, since the attention is mainly devoted to the methodological development of a reduced order method in an embedded boundary setting rather than in its efficiency, we do not rely on such hyper reduction techniques and we assemble the full-order differential operators also during the online stage. Considering that the most demanding computational effort is spent during the resolution of the full-order problem rather than in the assembly of the differential operators, as reported in \autoref{sec:num_exp}, it is anyway possible to achieve a good computational speedup. 
We remark here that, during the online stage, also the stabilization term ${\bm{C}}$ is projected onto the reduced basis space, the implications of such a choice are reported in the following section. 
\subsubsection{Stability Issues}\label{sec:stab_issues}
The reduced problem, as formulated in \autoref{sec:ROM}, may present stability issues. It is well known in fact that, using a mixed formulation for the approximation of the incompressible Stokes equations, the approximation spaces need to satisfy the ``inf-sup'' (Ladyzhenskaya-Brezzi-Babuska) condition, see \cite{BREZZI199027,boffi_mixed}. It is required that there should exist a constant $\beta > 0 $, independent to the discretization parameter $h$, such that:
\begin{equation}
\inf_{0\neq q \in {Q_h}} \sup_{0\neq \bm{v} \in {\bf V}_h} \frac{\langle \nabla \cdot \bm{v}, q \rangle}{||{\nabla \bm{v}||} \, ||{q}|| }\ge \beta > 0.
\end{equation}
Dealing with classic finite element methods, for what concerns the full-order level, this requirement can be met choosing appropriate finite element spaces such as the standard Taylor-Hood ($\mathbb{P}_2-\mathbb{P}_1$). In our case, at full-order level, as explained in \autoref{sec:HF}, in order to tackle this issue, we rely instead on a stabilized finite element technique. Since also the reduced order problem has a saddle point structure we need to ensure that a \emph{reduced} version of the LBB condition is fulfilled; 
Regardless the full-order discretization technique, even though the snapshots have been obtained by stable numerical methods, there is no guarantee that the original properties of the full-order system are preserved after the Galerkin projection onto the RB spaces \cite{RoVe07,Gerner2012,ballarin2015supremizer}. 
To overcome this issue, most of the contributions available in literature do not attempt to recover the pressure field and, at reduced order level, resolve only the momentum equation neglecting the contribution of the gradient of pressure. Although, as highlighted in \cite{noack2005}, in many applications the pressure term is needed  and cannot be neglected. In the proposed approach the velocity space is enriched in order to satisfy a reduced version of the ``inf-sup'' condition, see \cite{ballarin2015supremizer,RoVe07,quarteroni2007numerical}. 
\begin{remark}
It is important to remark that, during the projection stage, we perform also the projection of the stabilization terms. As we will see in the numerical experiments this helps to employ a relatively smaller number of additional supremizer modes. 
\end{remark}

\subsubsection{\it Supremizer enrichment}\label{subsec:sup_enrich}
As mentioned in \autoref{sec:stab_issues}, the problem, formulated using a mixed formulation, needs either to meet the ``inf-sup'' condition or to be properly stabilized. To achieve this, it is proposed the fulfillment of a reduced (and also parametric) version of the ``inf-sup'' condition. 
Within this approach, the velocity supremizer basis functions $\bm{L}_{\text{sup}}$ are computed and added to the reduced velocity space, which is transformed into $\tilde{\bm{L}}_u$, namely:
\begin{equation}
\begin{split}
\bm{L}_{\text{sup}} =[\bm{\eta}_1, \dots, \bm{\eta}_{N_{\text{sup}}^r}] \in \mathbb{R}^{N_u^h \times N_{\text{sup}}^r}, \mbox{ } \quad
\tilde{\bm{L}}_u = [\bm{\varphi}_1,\dots,\bm{\varphi}_{N_u^r},\bm{\eta}_1,\dots,\bm{\eta}_{N_{\text{sup}}^r}] \in \mathbb{R}^{N_u^h \times (N_u^r+N_{\text{sup}}^r)}.
\end{split}
\end{equation}
These basis functions are chosen solving a supremizer problem that ensures the fulfillment of a reduced version of the ``inf-sup'' condition. The supremizer solution ${\bm{s}_i}$, corresponding to the parameter value $\mu^i$, given a certain pressure basis function $\chi_i$, is the ``ingredient'' that permits the realization of the ``inf-sup'' condition. For each pressure basis function the corresponding supremizer element can be found solving the following problem:
\begin{equation}\label{eq:sup_problem}
\Delta {\bm{s}_i} = - {\nabla} \chi_i\mbox{ in } \mathcal D(\mu ^i), \quad{\bm{s}_i}=\bm{0}\mbox{ on } \Gamma(\mu ^i).
\end{equation}
In this case, since we want to rely on the same FOM solver, also the supremizer problem is solved with an SBM approach using the Poisson solver presented in \cite{MaSco17_2}. For more details regarding the derivation of the supremizer stabilization method one may see \cite{RoVe07,Gerner2012,ballarin2015supremizer}. According to the latter references, in order to have a parameter independent velocity enrichment, an \emph{approximated supremizer enrichment} procedure will be followed. This means that the supremizer problem is solved for each pressure snapshot $p(\mu)$ and a snapshots matrix of supremizer solutions is assembled:
\begin{equation}\label{eq:sup_snap}
\bm{S_{\text{sup}}} = [\bm{s}(\mu^1),\dots,\bm{s}(\mu^{N_s})] \in \mathbb{R}^{N_u^h\times N_s}.
\end{equation}
A proper orthogonal decomposition procedure is then applied to the resulting snapshots matrix in order to obtain a supremizer POD basis functions $\bm{\eta}_i$. 
Furthermore, the supremizer basis functions do not depend on the particular pressure basis functions but are computed during the offline phase, employing  the pressure snapshots.
\section{Numerical experiments}\label{sec:num_exp} 
\begin{figure} \centering
  \includegraphics[width=0.8\textwidth]{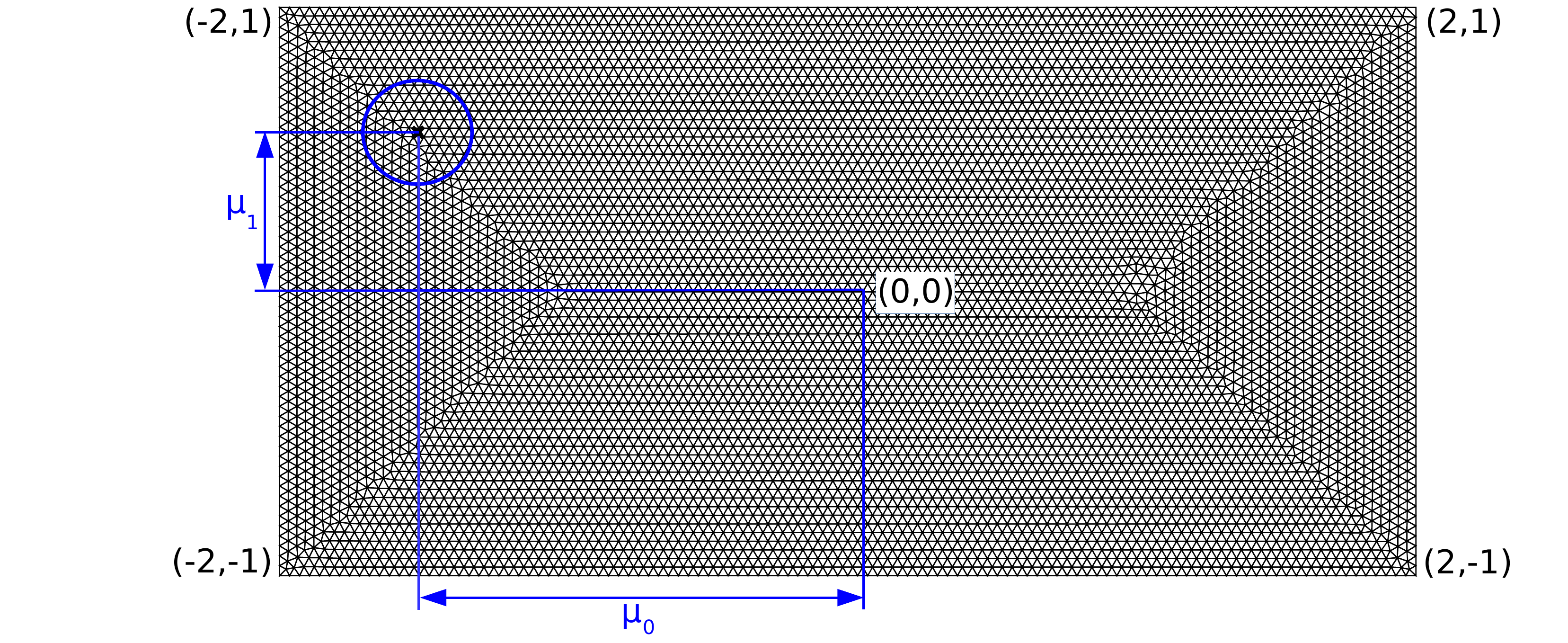}
  \caption{The background mesh together with a sketch of the embedded domain and the parameters considered in the numerical examples.}
  \label{background_mesh}
\end{figure}

\begin{figure} \centering
  \includegraphics[width=\textwidth]{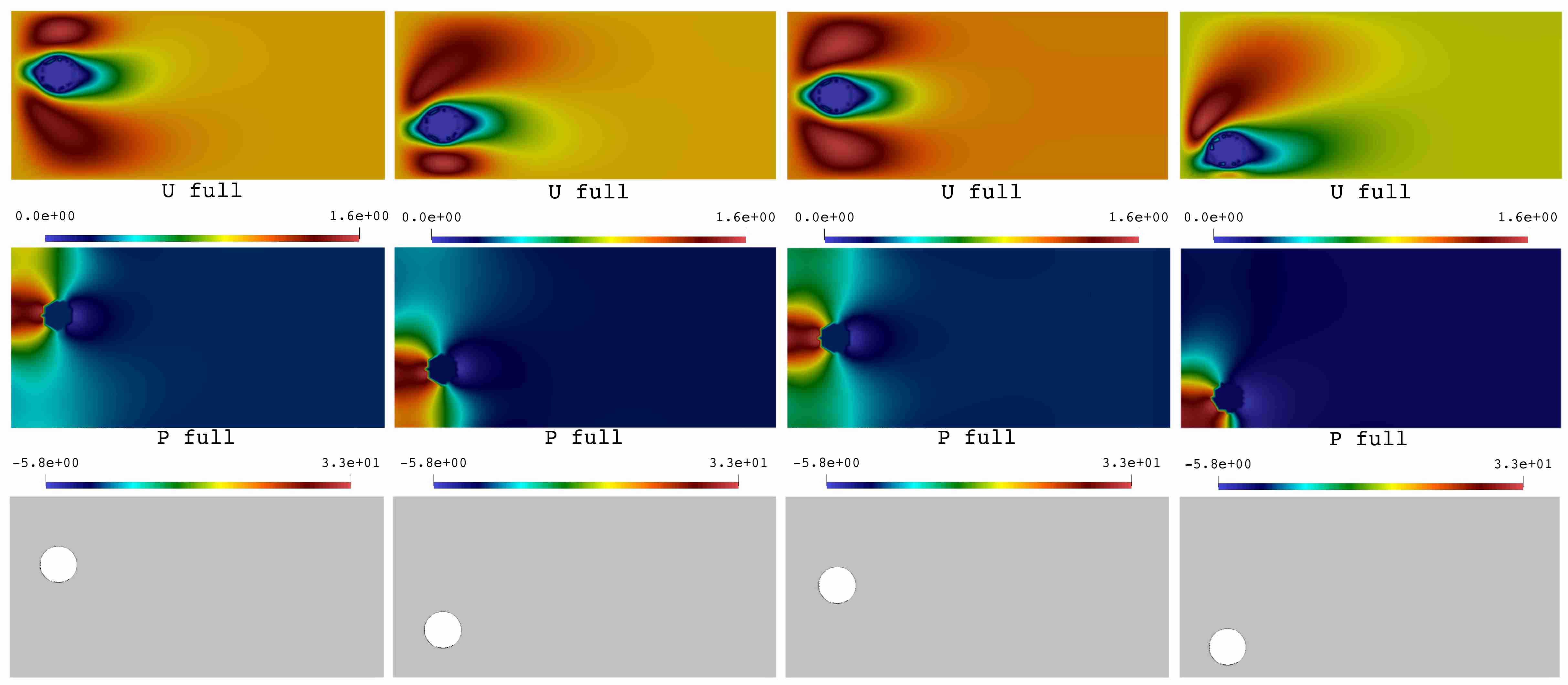}
  \caption{Shifted Boundary method: the full-order solution snapshots for $4$ different values of the input parameter in the case of 1D geometrical parametrization. In the figure we report the velocity $\bm u $ (first row), the pressure $p$ (second row) and the true geometry (third row). Each column corresponds to a different parameter value $\mu_1=[0.2439, -0.4764, 0.0233, 0.6468]$, while $\mu_0$ is fixed to the value $-1.5$.}
  \label{StokesParameterized} 
\end{figure}

In the present section we will test the presented methodology on numerical tests investigating a steady Stokes flow around an embedded circular cylinder {\blue and a Stokes flow around a more complex geometry morphed using the Free Form Deformation.} In the first two numerical example the embedded domain {\blue is parametrized} through $\mu=(\mu_0,\mu_1)$ according to the following expression:
$$
(x-\mu_0)^2 +(y-\mu_1)^2 \le R^2,
$$
where the two parameters $\mu_0$ and $\mu_1$ describe the $x$ and $y$ coordinates of the center of the embedded circular domain, see e.g. \autoref{background_mesh}. We consider two different test cases, the first one consists of an 1D geometrical parametrization where the first parameter $\mu_0$ is fixed. The second one consists of a 2D geometrical parametrization where both parameters $\mu_0$ and $\mu_1$ are left free. The problem domain is the rectangle ${\mathcal{D}}= [-2, 2] \times [-1, 1]$, in which a cylinder with constant radius $0.2$ is embedded. The viscosity $\nu$ is set to $1$. A constant velocity in the $x$ direction, $u_{\text{in}} = 1$ is applied at the left side of the domain, and an open boundary condition with $\nabla \bm u \cdot \bm n = p_{\text{out}} = 0$ is applied on the right. A slip (no penetration) boundary condition is applied on the top and bottom edges. On the boundary of the embedded cylinder a no slip boundary condition is applied.
The results for the test problems have been obtained with a mesh size of $h = 0.0350
$ for the background mesh, see e.g. \autoref{background_mesh}, using 15022 triangles for the discretization and $\mathbb{P}1/\mathbb{P}1$ finite elements in space with stabilization as described in \cite{Hu87}. All the numerical experiments have been tested with and without supremizer basis enrichment.
Some representative snapshots are visualized in \autoref{StokesParameterized}. 

Both the full-order and the ROM simulation were run in serial on a personal computer with an Intel\textsuperscript{\textregistered} Core\textsuperscript{TM}  i7-4770HQ 3.70GHz CPU. In order to test the accuracy of the proposed method we compared the reduced order solutions against the full-order ones using $10$ different samples of the parameters that were not previously tested during the training stage. In both cases the ROMs are constructed using the methodologies described in \autoref{sec:ROM}. 

For the generation of the POD spaces, we considered $1024$ full-order snapshots for the velocity and pressure fields. The snapshots are collected and $N^r_u=j$ velocity modes, $N^r_p=6j$ pressure modes, and whenever we use supremizers $N^r_{\text{sup}}=4j$ modes have been selected with $j=8, 12,16,20,25,30,35,40,45,50$. We mention here that we chose to use a larger number of pressure and supremizer modes compared to the number of velocity modes after testing different configurations that we do not report here for sake of brevity. This decision was, in fact, taken  after testing other cases e.g. $N^r_p=j, 2j, 3j, 4j, 6j$, $N^r_{\text{sup}}=j,2j,3j,4j$ and its various combinations. We highlight here that a number of supremizer modes which is different with respect to the number of pressure modes is not a classical choice, since in the literature $N^r_p=N^r_{\text{sup}}$ is usually used. As mentioned in \autoref{sec:ROM}, it has been in fact heuristically verified \cite{ballarin2015supremizer,stabile_stabilized} that, using an approximated supremizer approach, the minimum number of supremizer modes is equal to the number of pressure basis functions. In the present case, since at the reduced order level we are projecting also the stabilization term, the loss of stability given by velocity divergence free modes is partially circumvented by a ``reduced version'' of the stabilization terms. \footnote{\blue We remark that this topic deserves further investigation. In fact, we notice that projection-based ROMs are quite sensitive with respect to the selection of the FOM stabilization technique. In \cite{Shafqat}, with a similar residual based SUPG stabilization, especially for the case of geometrical parametrization using a reference domain approach, we notice that the supremizer stabilization leads to a better approximation of the pressure field. On the other hand, when a variational multiscale stabilization method is used, we notice that the supremizer stabilization is not needed and the ROM coefficients multiplying the supremizer basis functions turned out to be zero \cite{stabile_vms}.}

\subsection{1D geometrical parametrization}
\begin{figure} \centering
\begin{minipage}{\textwidth}
\centering
\begin{minipage}{0.24\textwidth}
  \includegraphics[width=\textwidth]{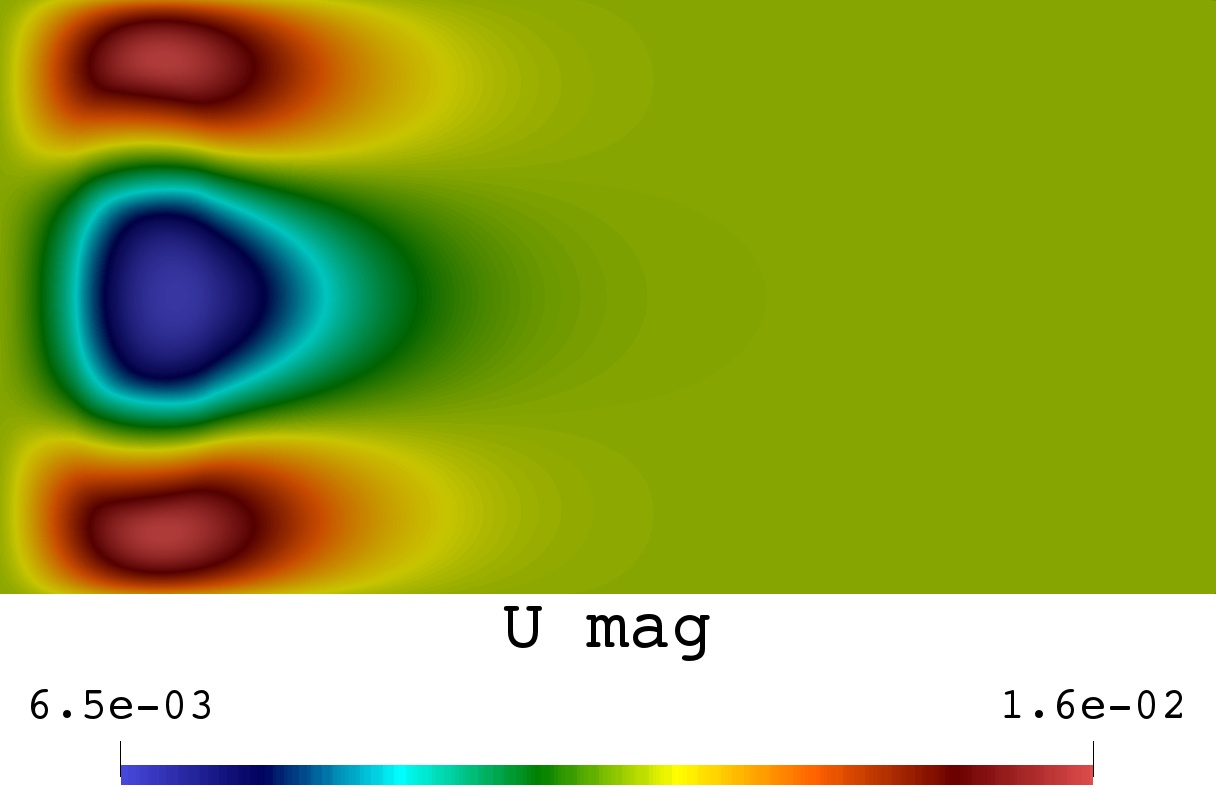}
\end{minipage}
\begin{minipage}{0.24\textwidth}
  \includegraphics[width=\textwidth]{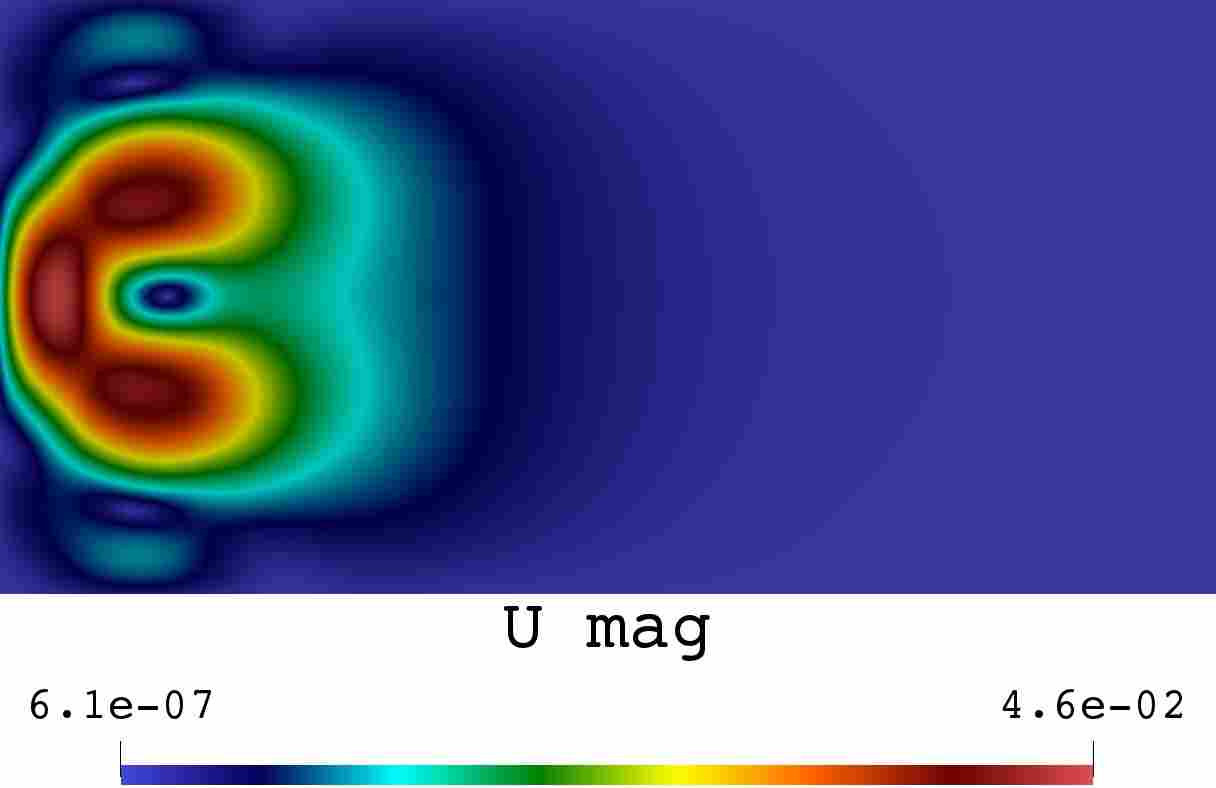}
\end{minipage}
\begin{minipage}{0.24\textwidth}
  \includegraphics[width=\textwidth]{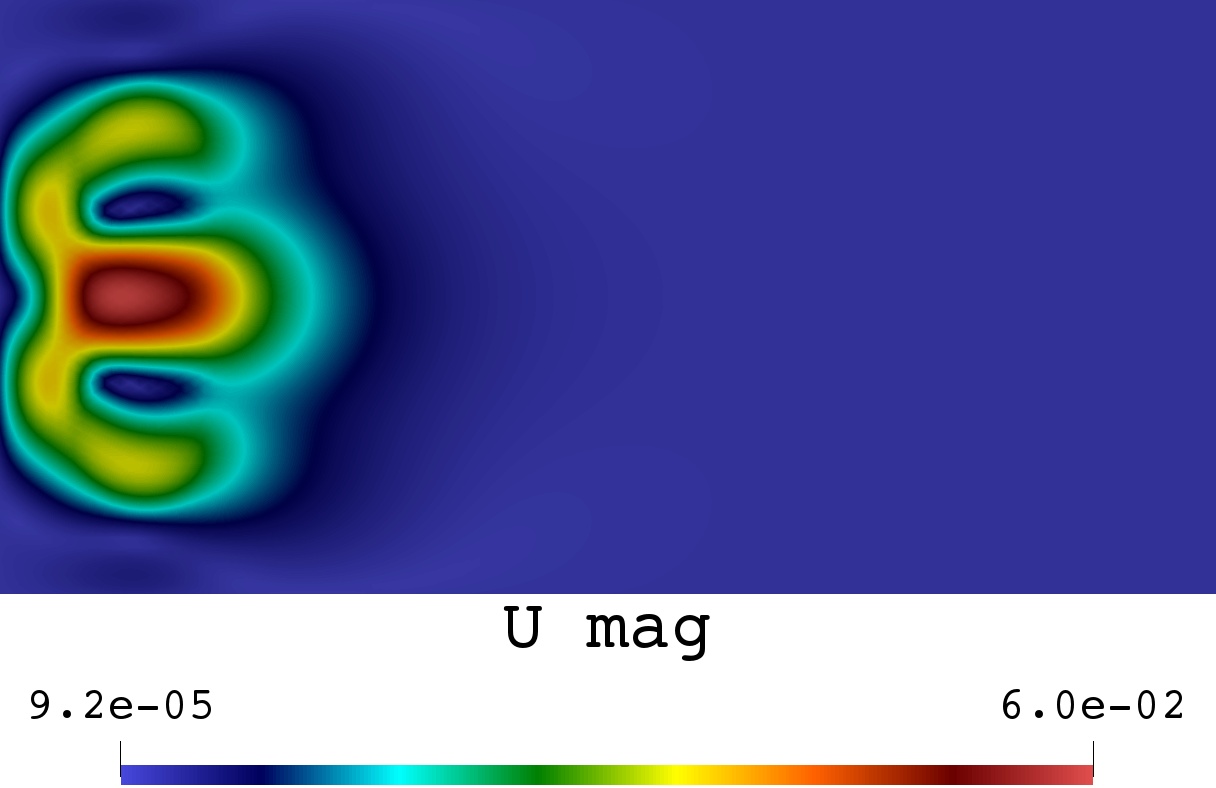}
\end{minipage}
\begin{minipage}{0.24\textwidth}
  \includegraphics[width=\textwidth]{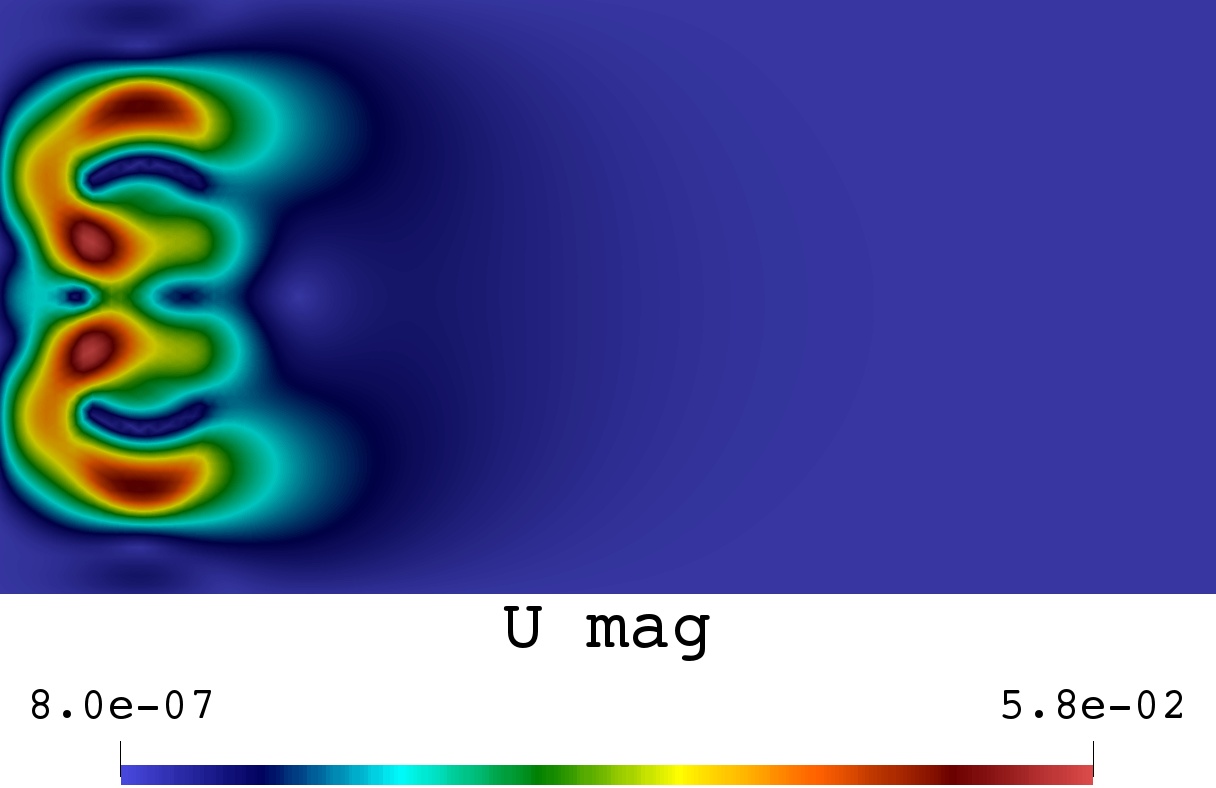}
\end{minipage}
\begin{minipage}{0.24\textwidth}
  \includegraphics[width=\textwidth]{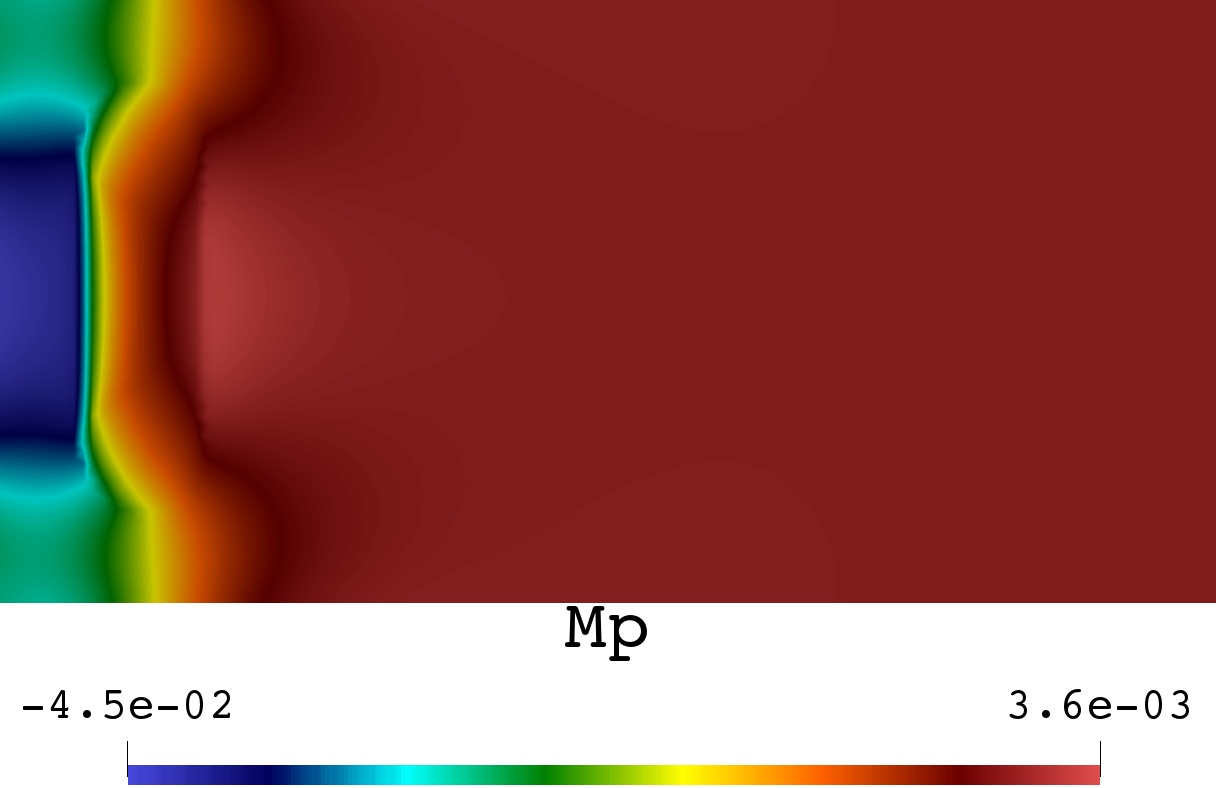}
\end{minipage}
\begin{minipage}{0.24\textwidth}
  \includegraphics[width=\textwidth]{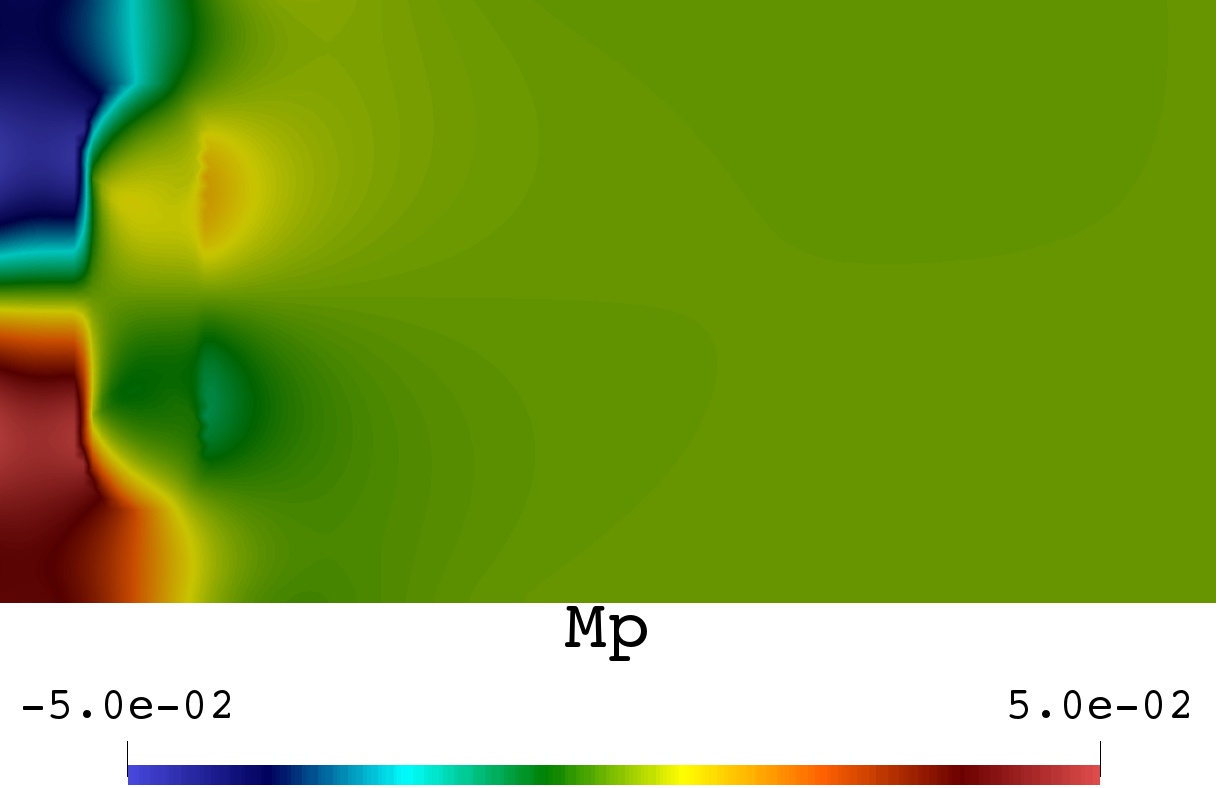}
\end{minipage}
\begin{minipage}{0.24\textwidth}
  \includegraphics[width=\textwidth]{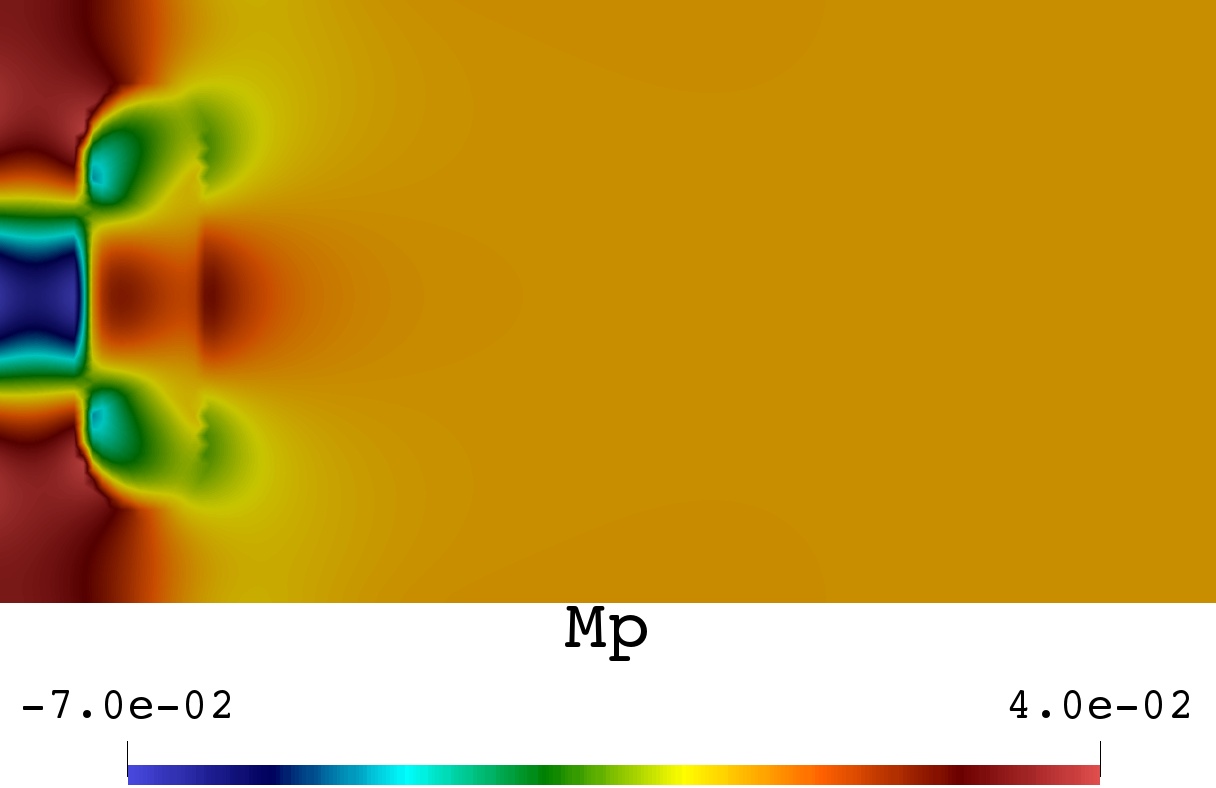}
\end{minipage}
\begin{minipage}{0.24\textwidth}
  \includegraphics[width=\textwidth]{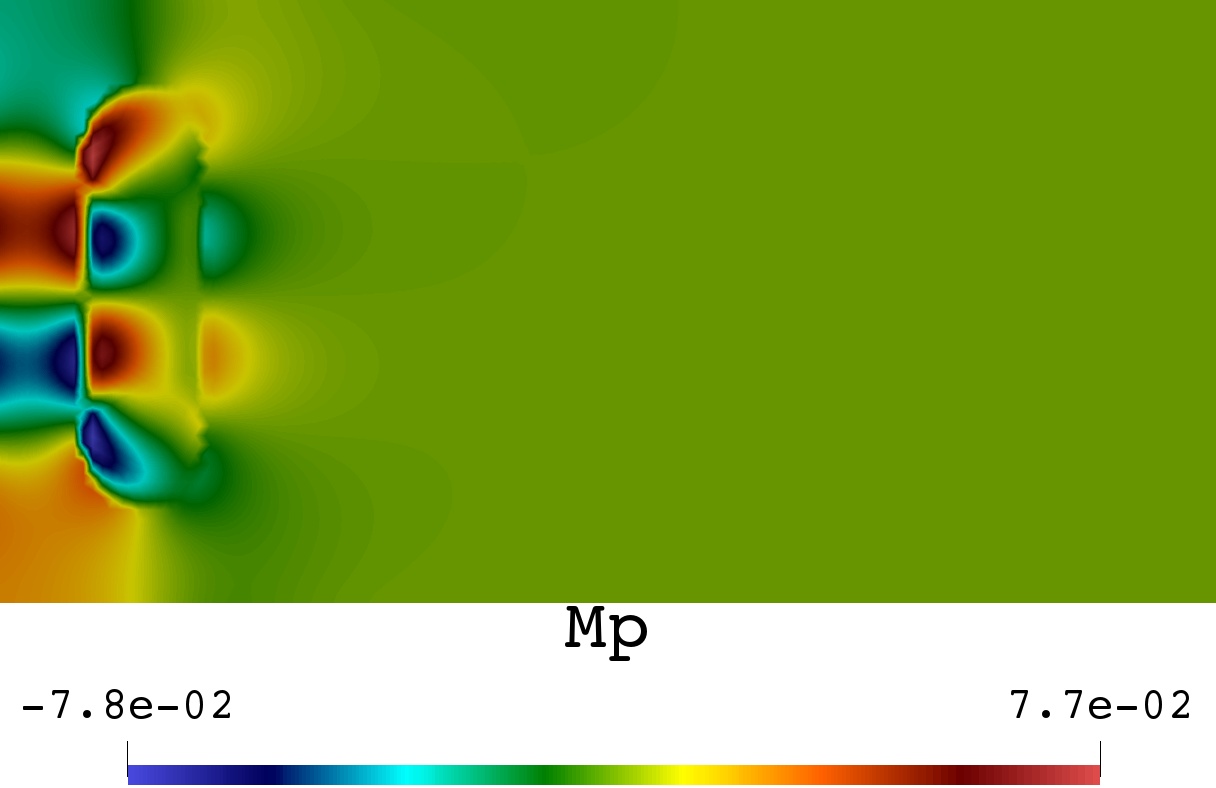}
\end{minipage}
\end{minipage}
  \caption{The first $4$ basis components for velocity (first row) and pressure (second row) in the 1D geometrical parametrization case with $\mu_1 \in [-0.65, 0.65]$.}
  \label{Stokes_Pressure_and_Velocity_Components_Modes}
\end{figure}

In the first test case the first parameter is fixed $\mu_0=-1.5$ while $\mu_1$ is left free and parametrized. The training set is chosen with an equally spaced distribution $\mu_1\in [-0.65, 0.65]$. 
The ROM results are compared against full-order results, see e.g.  \autoref{table:errors_no_supremizers} where the comparison is shown directly on the full-order and reduced-order velocity and pressure fields. In \autoref{Stokes_Pressure_and_Velocity_Components_Modes} we report the first four modes for the velocity, and the pressure fields.

\subsubsection{The role of the supremizer stabilization}
In this subsection we investigate the role of the ``inf-sup'' stabilization technique applied on an SBM stabilized full-order solver which, from the full-order point of view, allows the use of $\mathbb{P}1/\mathbb{P}1$ finite elements. The aim is to understand if the full-order stabilization is propagated also at the reduced order level. In general, considering incompressible flow applications, the quantity of interest is in both flow velocity and pressure. However, in many research contributions available in literature, due to the numerical instabilities given by spurious pressure modes, the pressure term is usually neglected \cite{noack2005}. In this work the attention is given to both velocity and pressure and in these first ROM-SBM experiments we want to test if the reduced version of the full-order Shifted Boundary Method solvers, without any additional stabilization technique, produce acceptable results. 

The results in \autoref{table:errors_no_supremizers} and \autoref{fig:Supremizer_solution_1000_100_1000_100_1024_j_3J_4j} (left plot) show  the mean relative error $||u-u_r||_{L^2(\mathcal D)}/ ||u||_{L^2(\mathcal D)}$ and  $||p-p_r||_{L^2(\mathcal D)}/ ||p||_{L^2(\mathcal D)}$ for a 10 samples testing set and for various different dimensions of the reduced basis spaces.
The ROM has been generated starting from $1024$ full-order snapshots with $j$ velocity modes, $6j$ pressure modes with and without $4j$ supremizer modes.
{\blue The comparison shows} that even without an additional online stabilization approach, as mentioned in \autoref{subsec_POD_theory}, the pressure relative error reaches acceptable values and has a better behavior respect to cases according to the classic FEM-ROM literature without pressure stabilization methods \cite{ballarin2015supremizer}. For the interested reader we refer for instance to \cite{RoVe07}. 
Also in the case without supremizer stabilization, even with a relatively low number of modes, it is possible to observe good approximation properties for both the velocity and the pressure field. 

The supremizer enrichment strategy permits to further reduce the approximation error of the pressure field, without increasing the approximation error of the velocity field. To summarize, in all the numerical experiments, the reduced velocity and pressure solutions reach a reasonable level of accuracy and in case with a supremizer enrichment stabilization it is possible to further decrease the approximation error of the pressure field without deteriorating the accuracy of the reduced velocity solution. 

\begin{figure} 
\centering
\begin{minipage}{\textwidth}
\centering
\begin{minipage}{0.24\textwidth}
  \includegraphics[width=\textwidth]{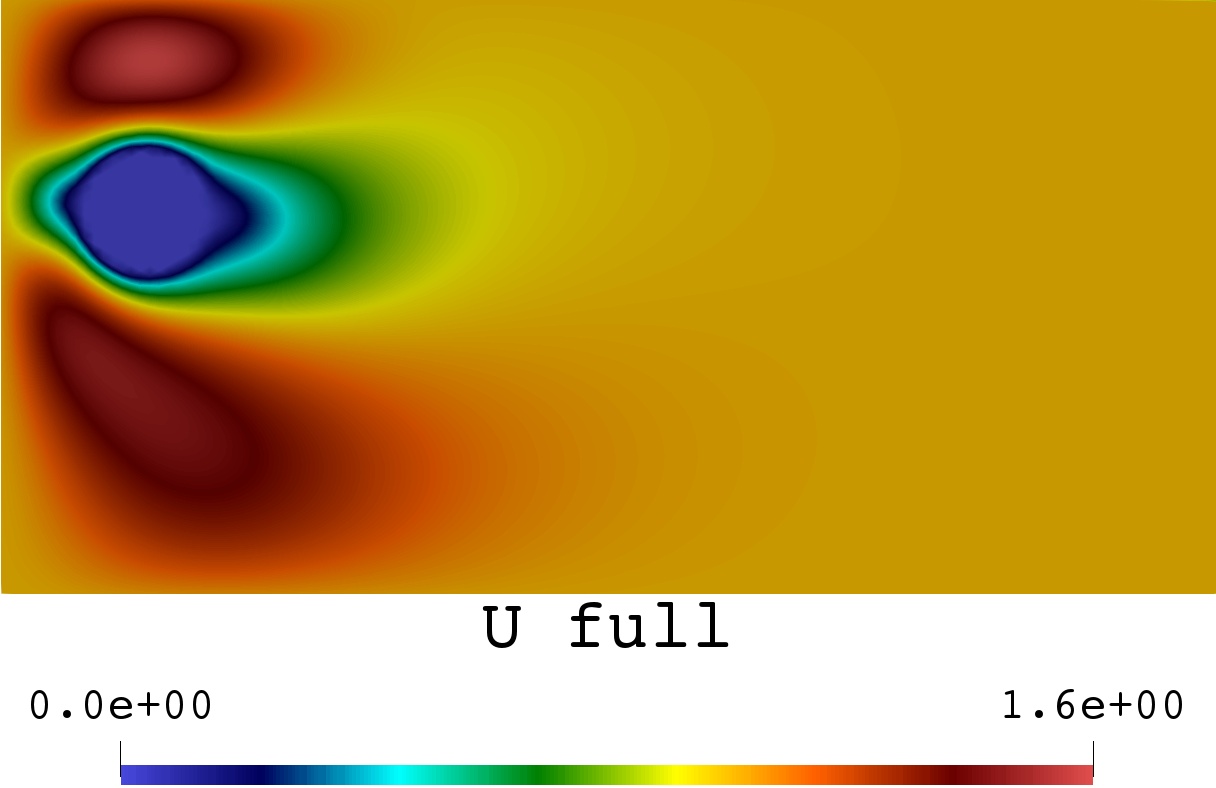}
\end{minipage}
\begin{minipage}{0.24\textwidth}
  \includegraphics[width=\textwidth]{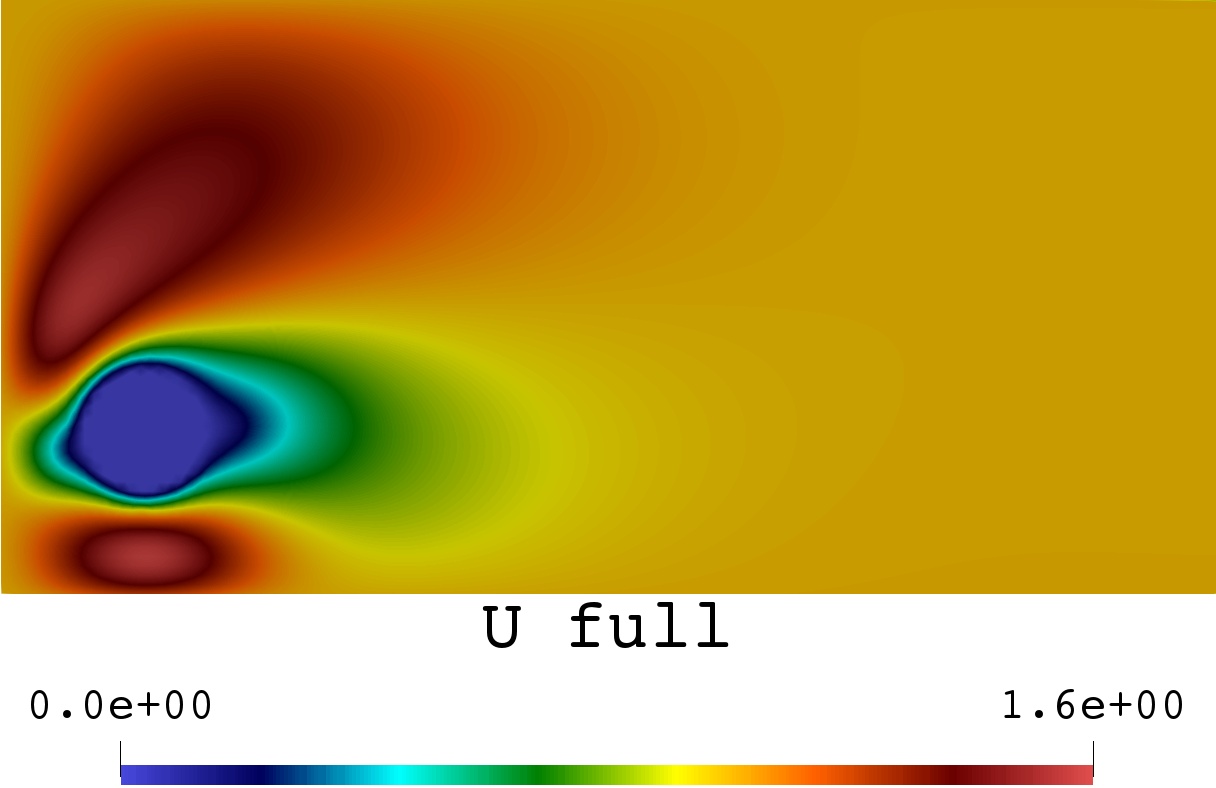}
\end{minipage}
\begin{minipage}{0.24\textwidth}
  \includegraphics[width=\textwidth]{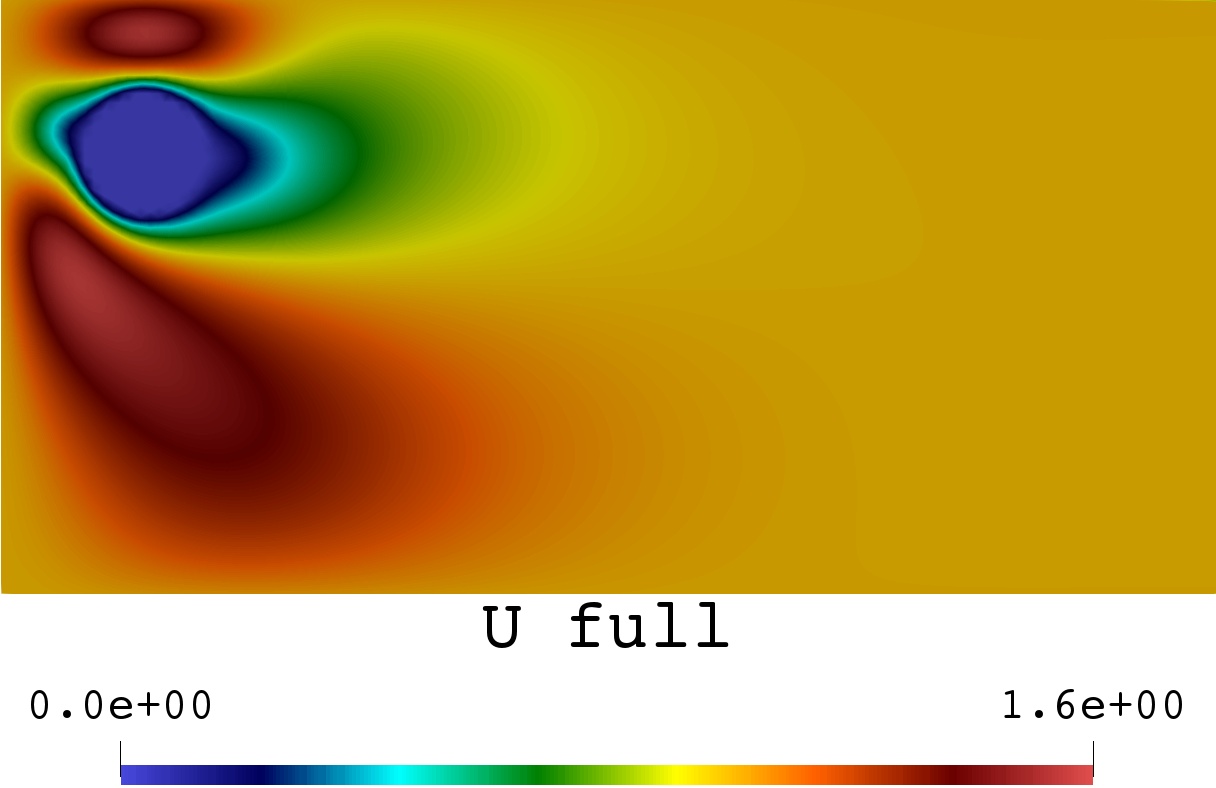}
\end{minipage}
\begin{minipage}{0.24\textwidth}
  \includegraphics[width=\textwidth]{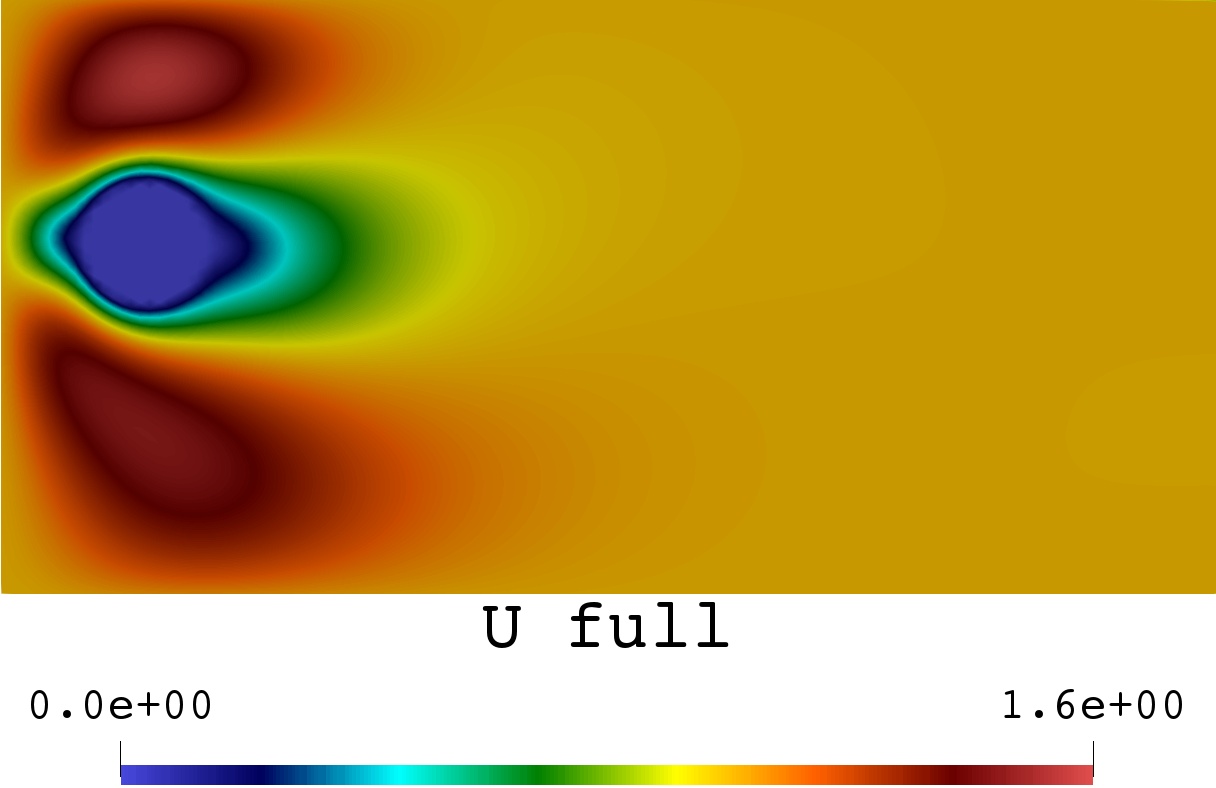}
\end{minipage}
\begin{minipage}{0.24\textwidth}
  \includegraphics[width=\textwidth]{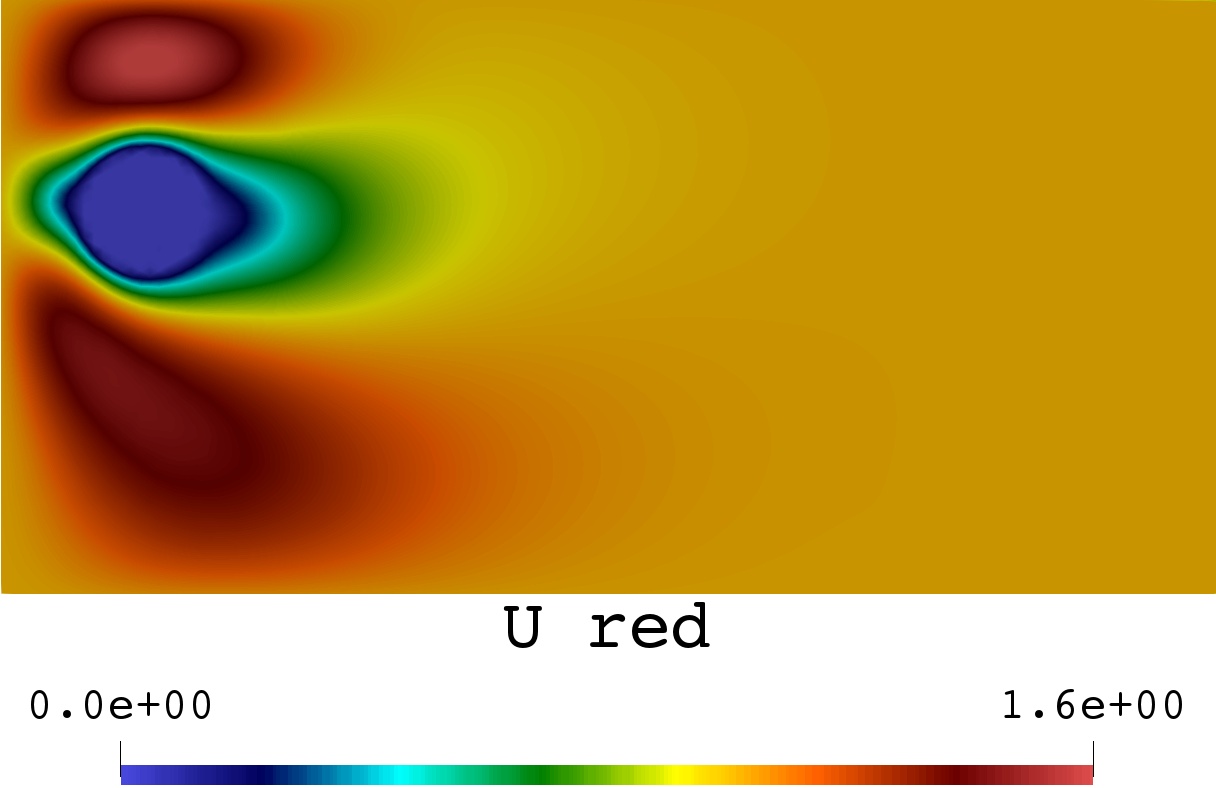}
\end{minipage}
\begin{minipage}{0.24\textwidth}
  \includegraphics[width=\textwidth]{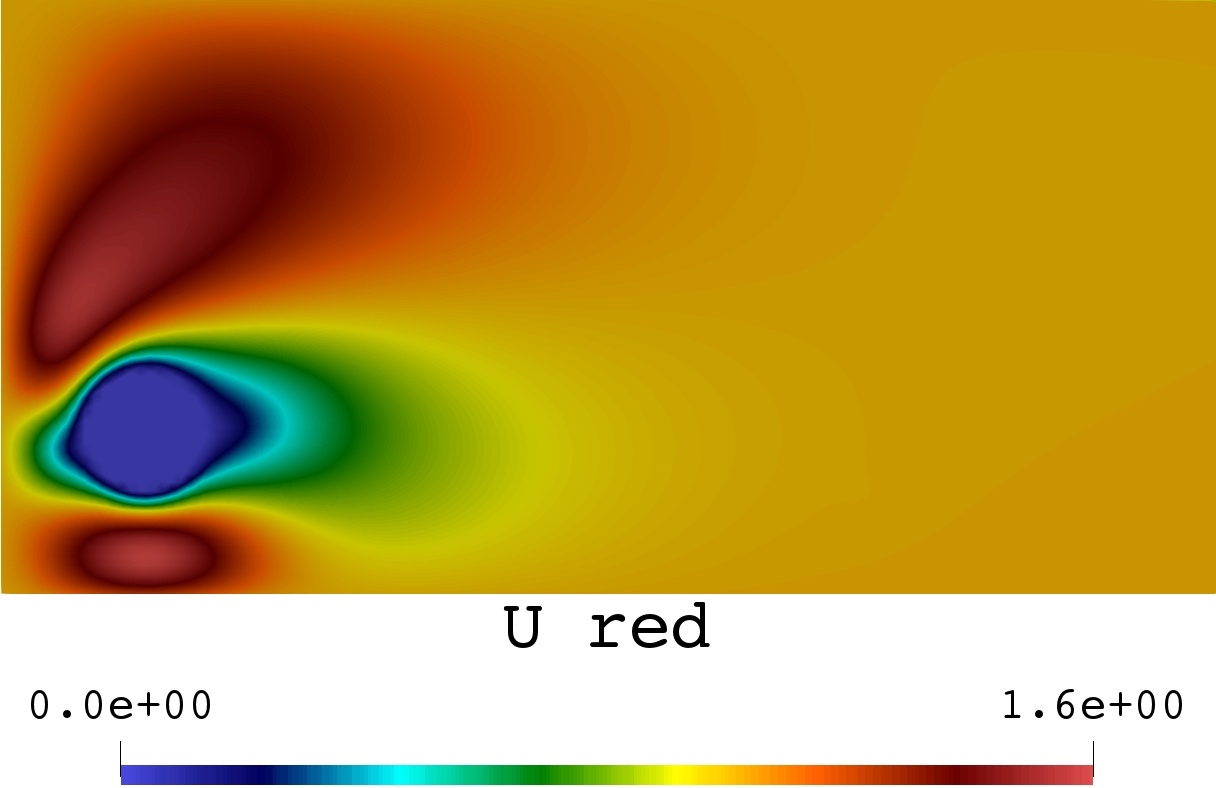}
\end{minipage}
\begin{minipage}{0.24\textwidth}
  \includegraphics[width=\textwidth]{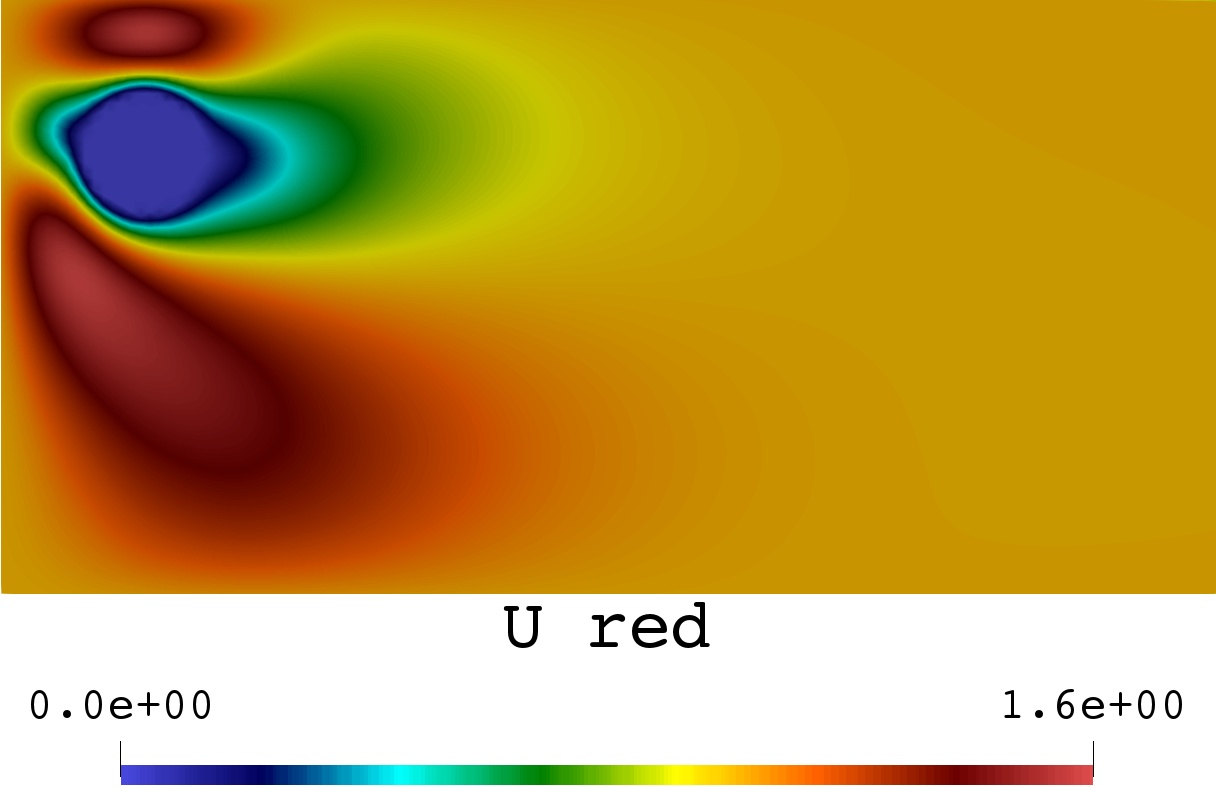}
\end{minipage}
\begin{minipage}{0.24\textwidth}
  \includegraphics[width=\textwidth]{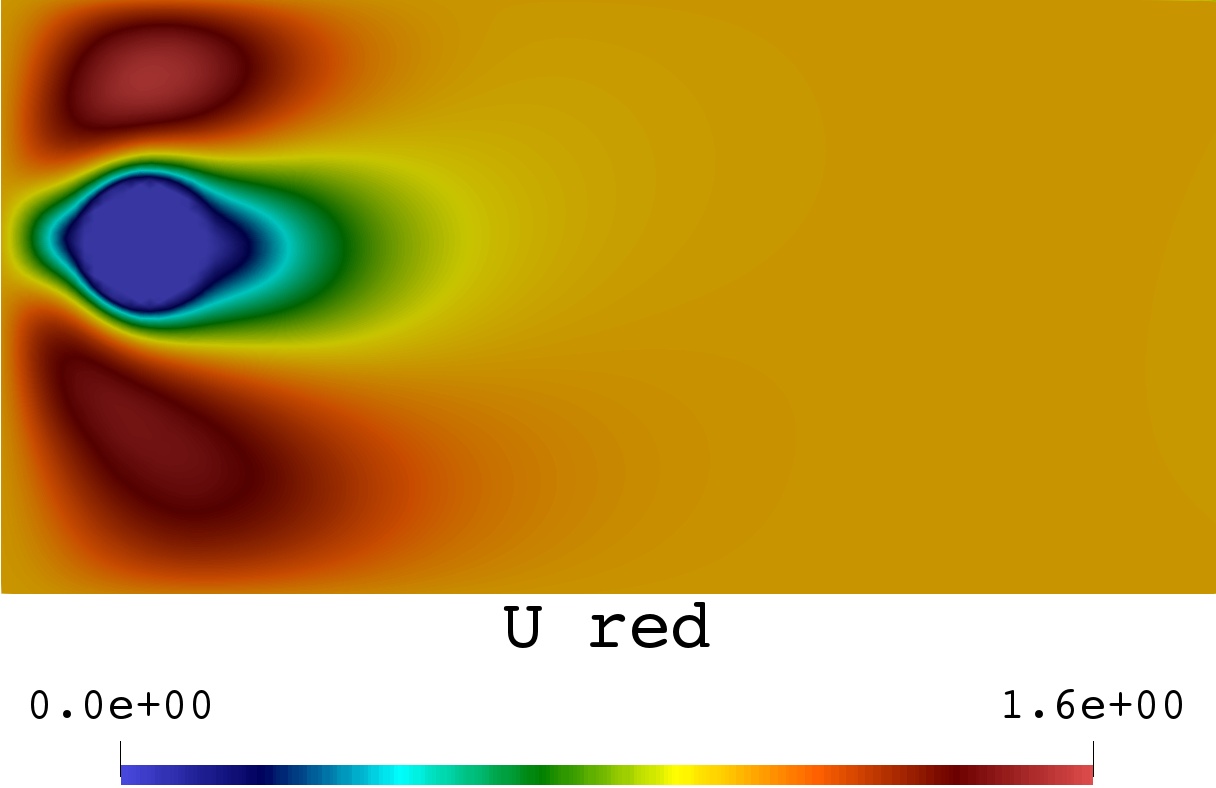}
\end{minipage}
\begin{minipage}{0.24\textwidth}
  \includegraphics[width=\textwidth]{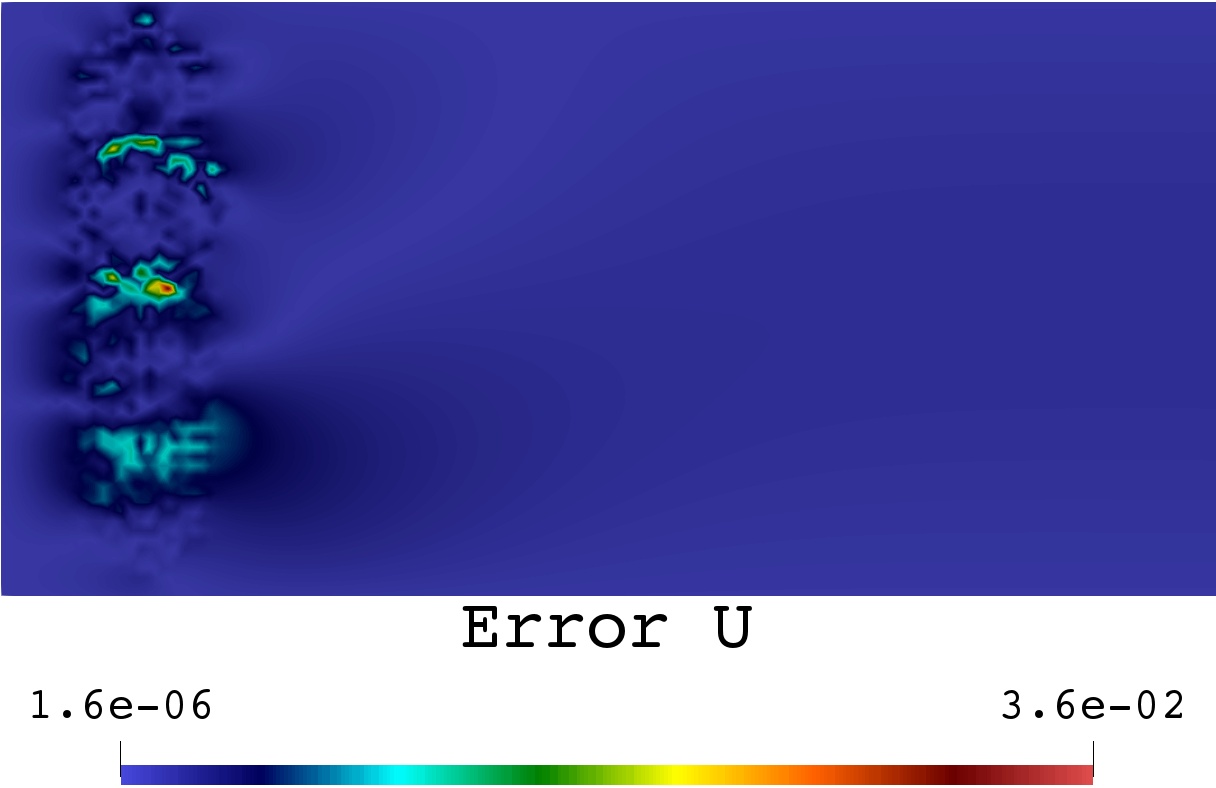}
\end{minipage}
\begin{minipage}{0.24\textwidth}
  \includegraphics[width=\textwidth]{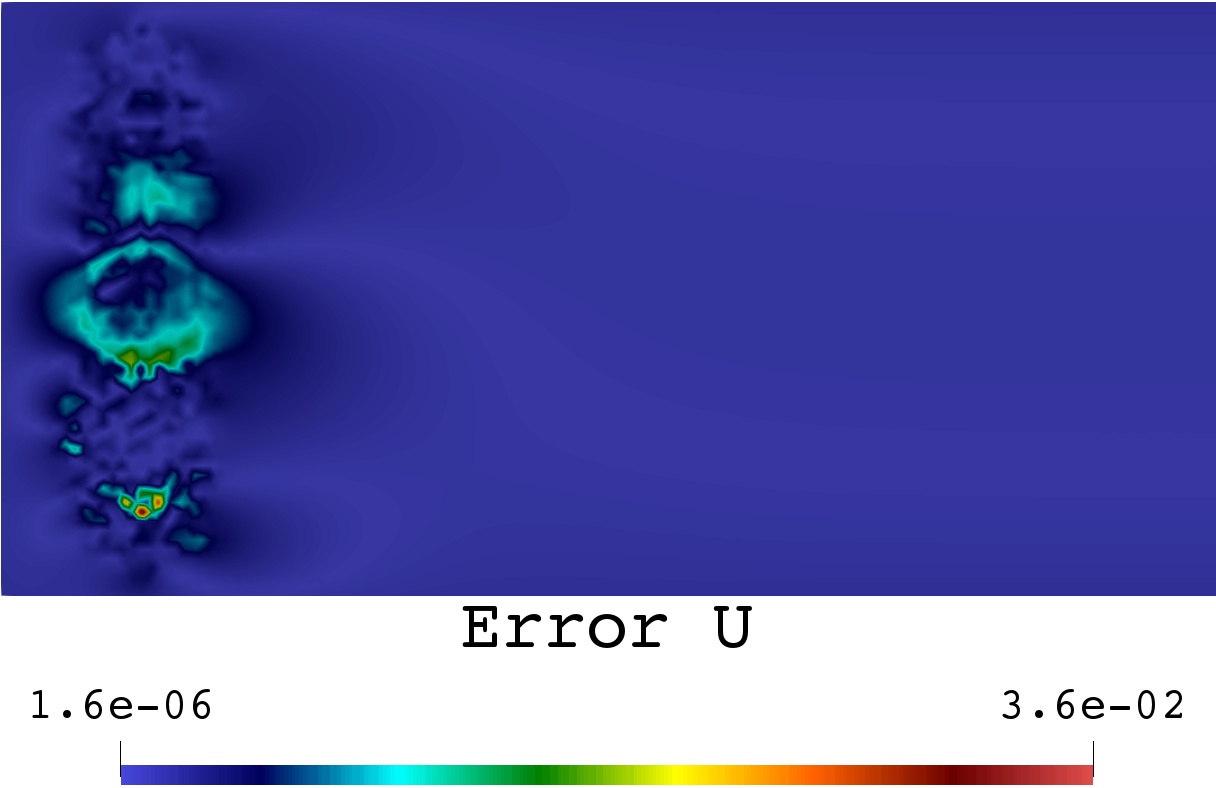}
\end{minipage}
\begin{minipage}{0.24\textwidth}
  \includegraphics[width=\textwidth]{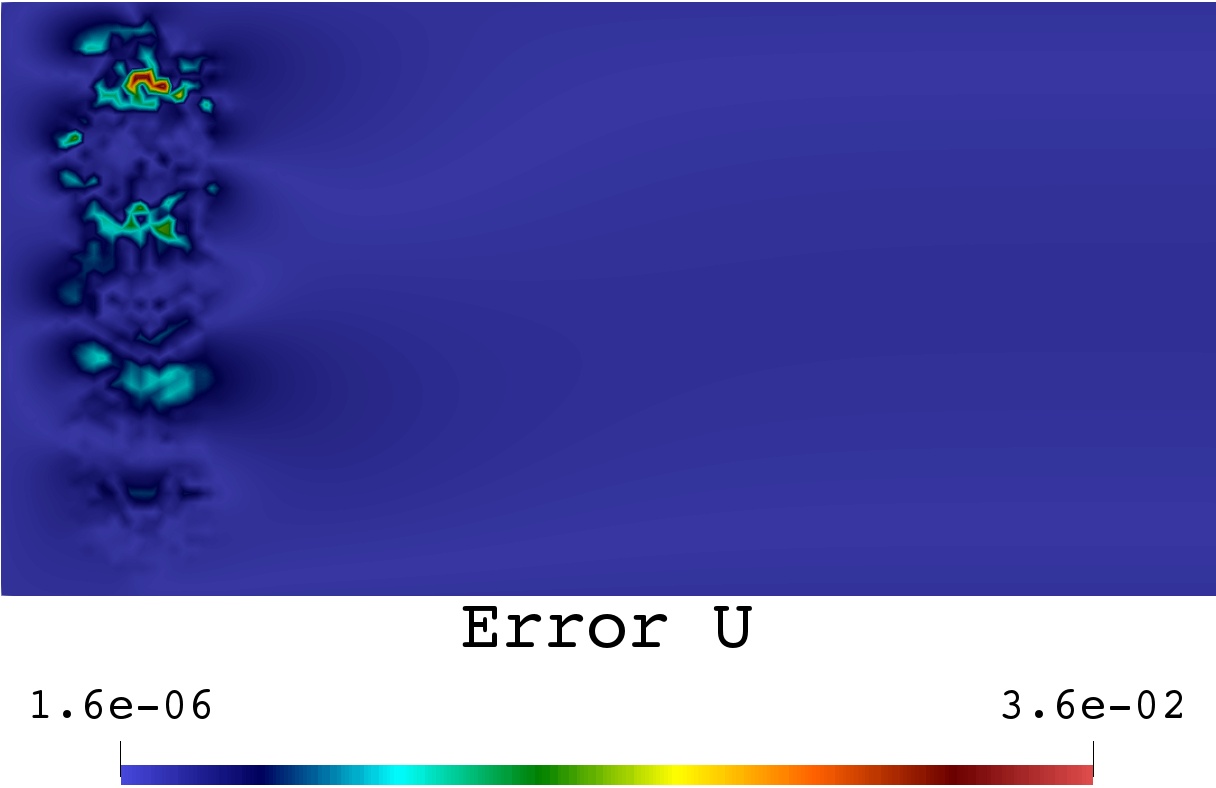}
\end{minipage}
\begin{minipage}{0.24\textwidth}
  \includegraphics[width=\textwidth]{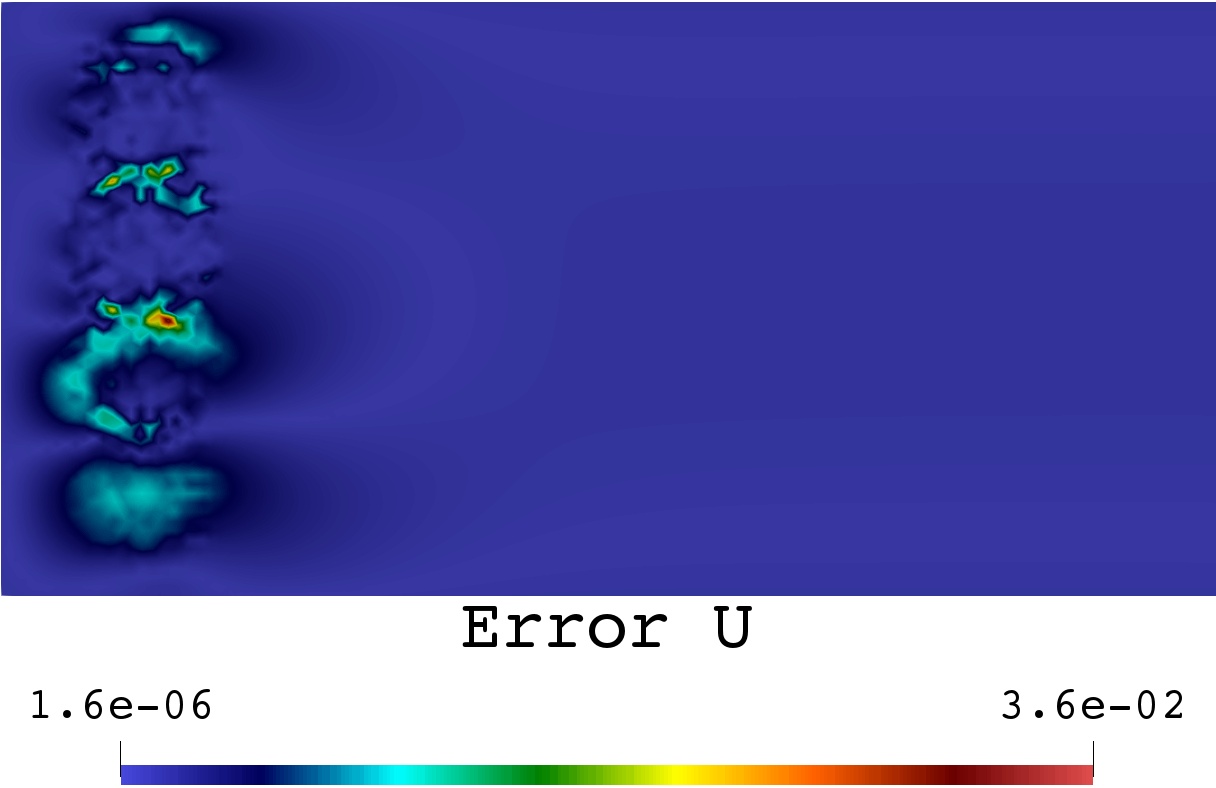}
\end{minipage}
\begin{minipage}{0.24\textwidth}
  \includegraphics[width=\textwidth]{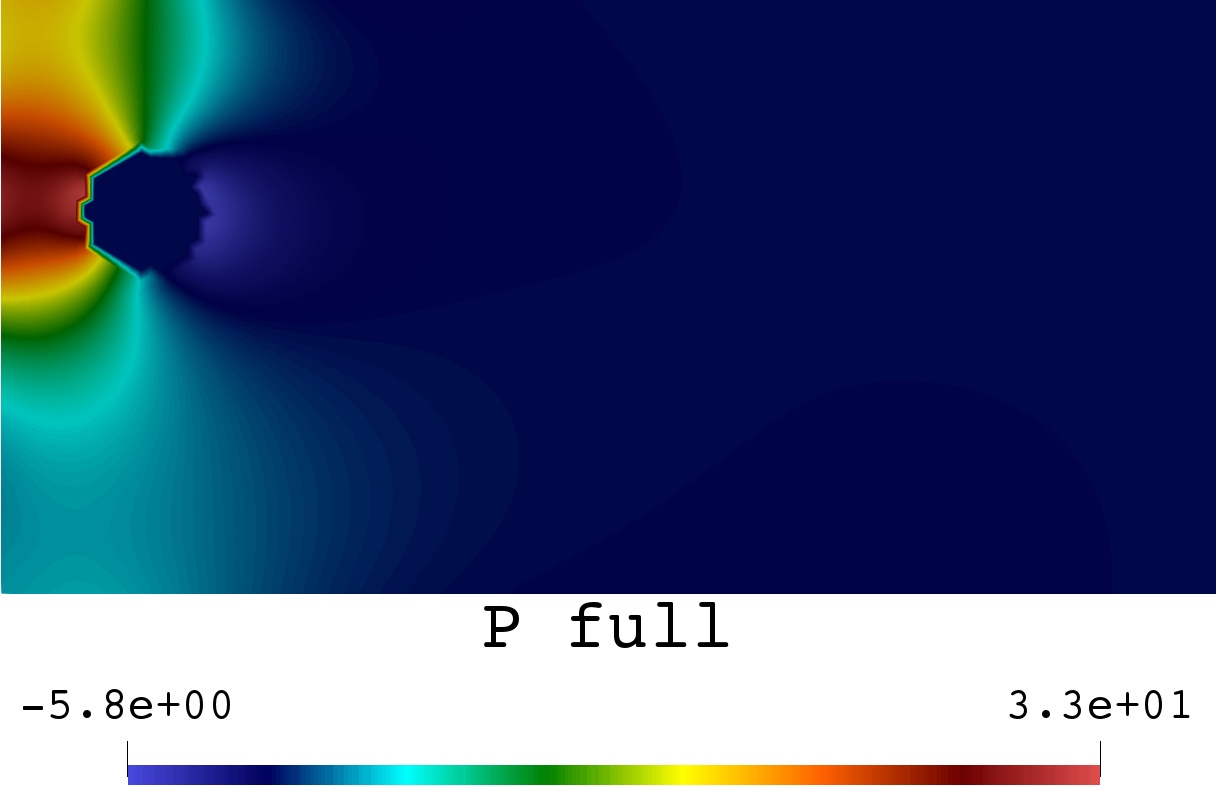}
\end{minipage}
\begin{minipage}{0.24\textwidth}
  \includegraphics[width=\textwidth]{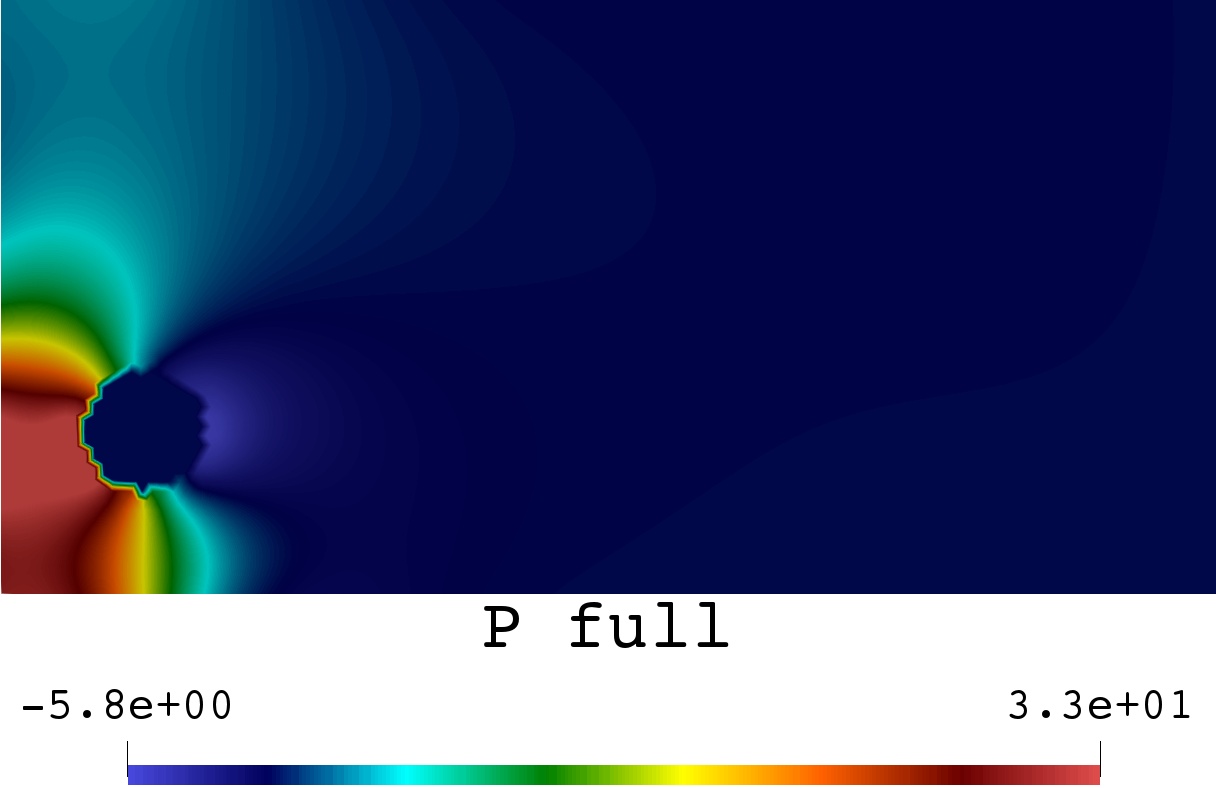}
\end{minipage}
\begin{minipage}{0.24\textwidth}
  \includegraphics[width=\textwidth]{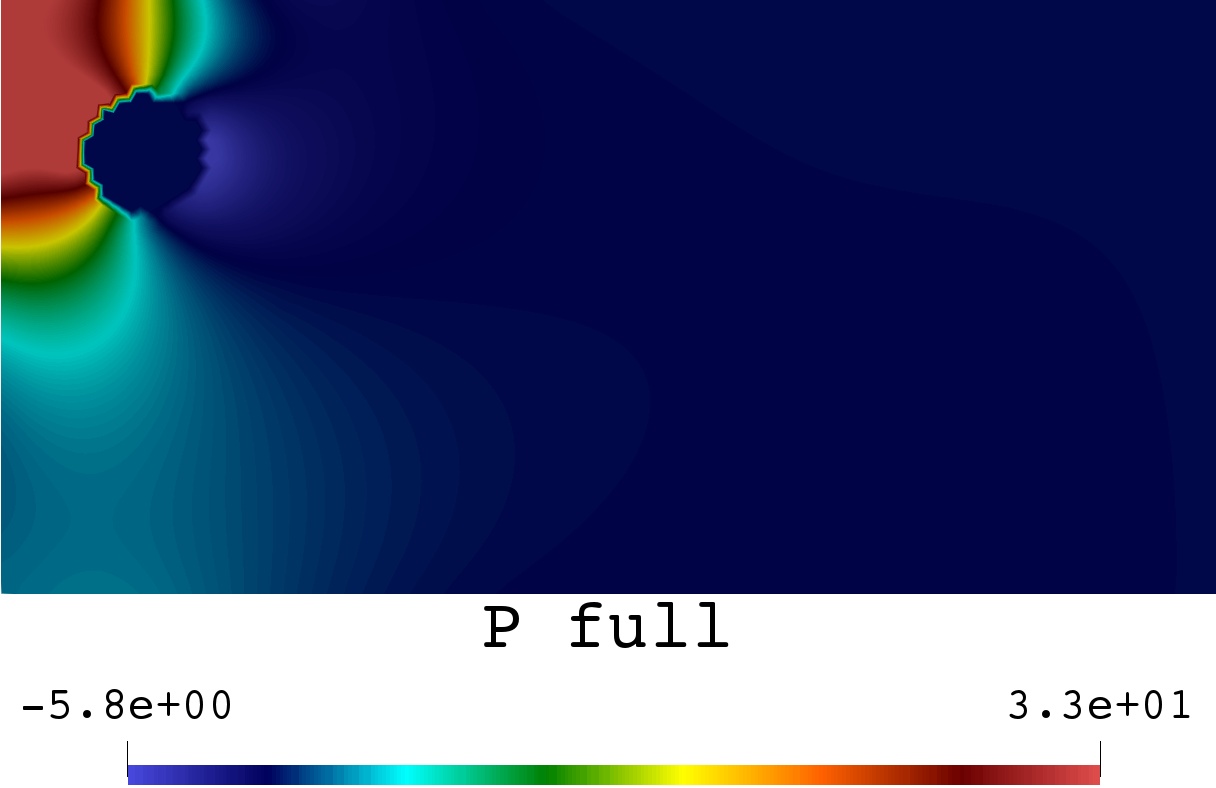}
\end{minipage}
\begin{minipage}{0.24\textwidth}
  \includegraphics[width=\textwidth]{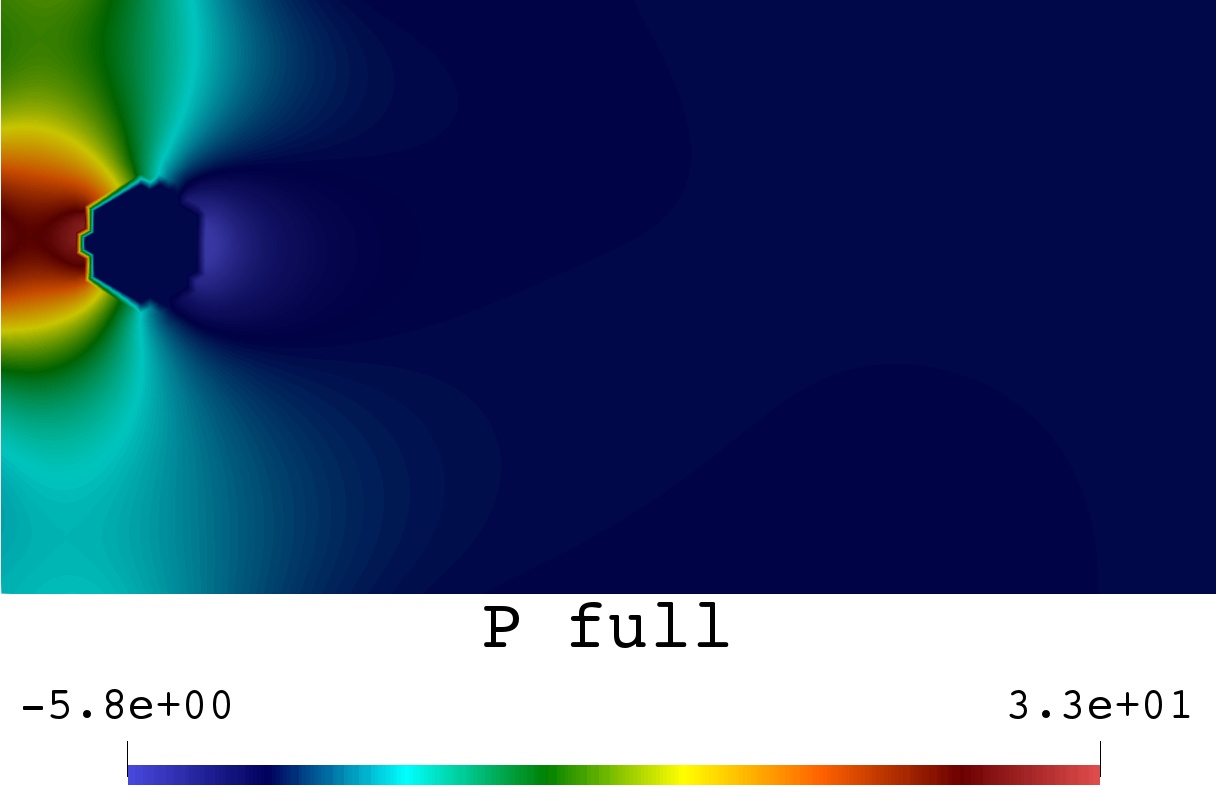}
\end{minipage}
\begin{minipage}{0.24\textwidth}
  \includegraphics[width=\textwidth]{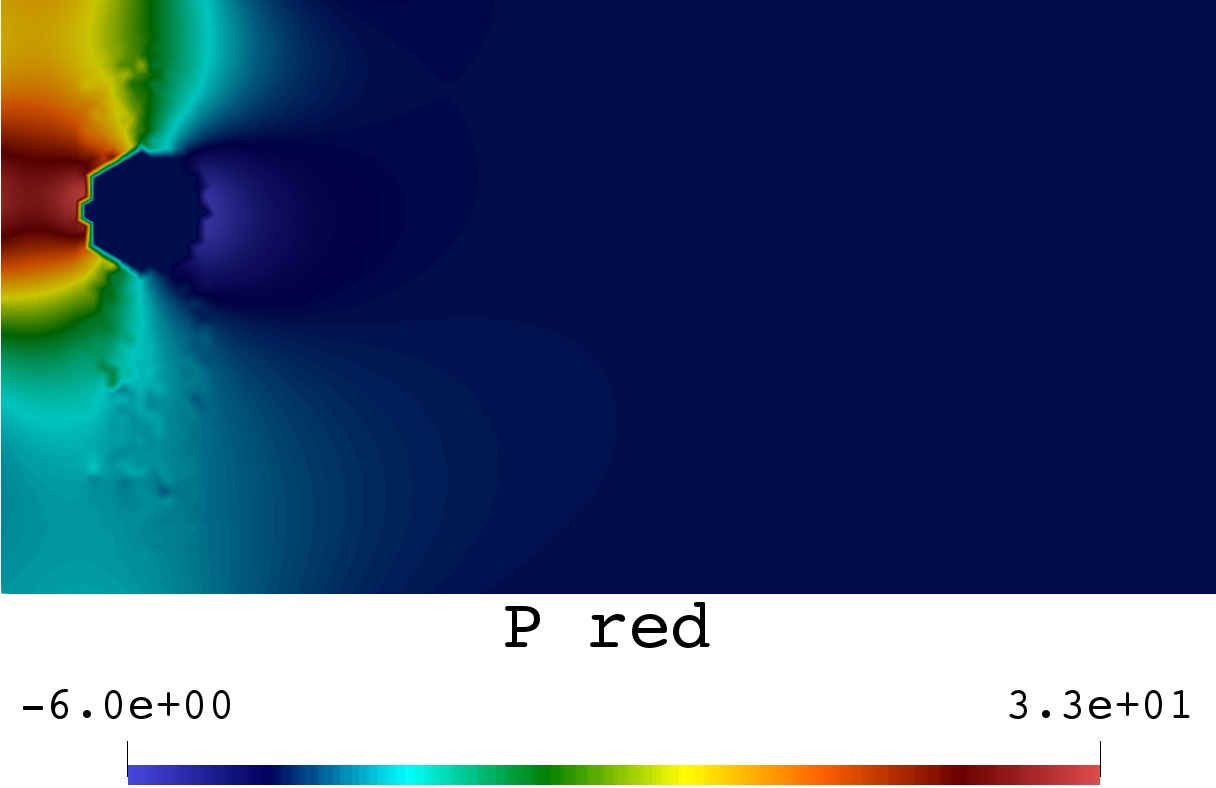}
\end{minipage}
\begin{minipage}{0.24\textwidth}
  \includegraphics[width=\textwidth]{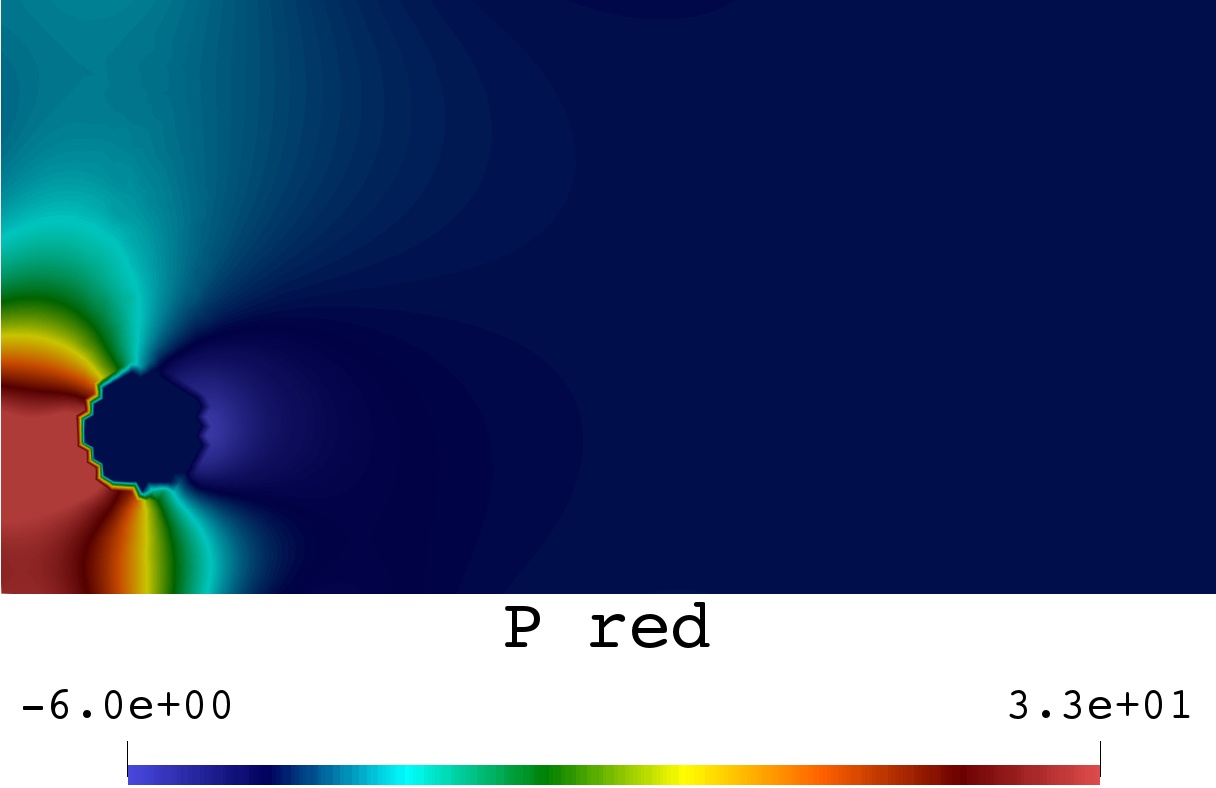}
\end{minipage}
\begin{minipage}{0.24\textwidth}
  \includegraphics[width=\textwidth]{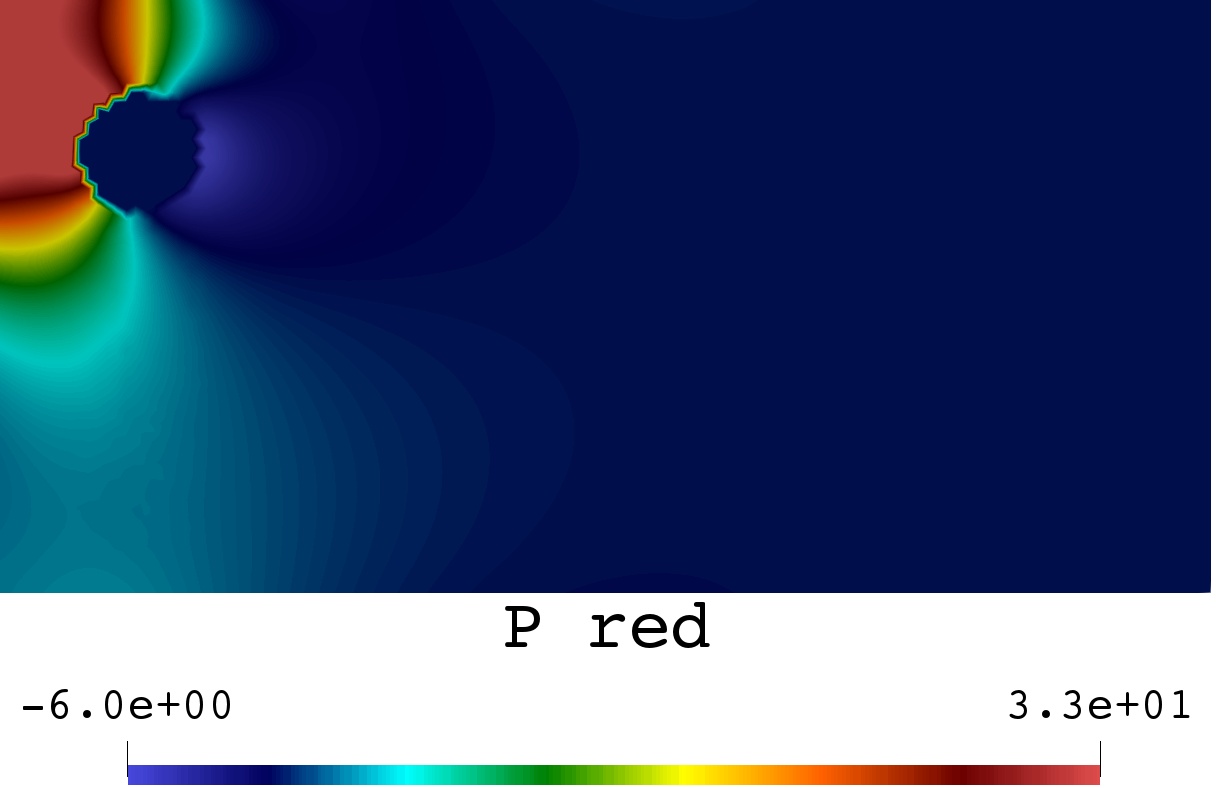}
\end{minipage}
\begin{minipage}{0.24\textwidth}
  \includegraphics[width=\textwidth]{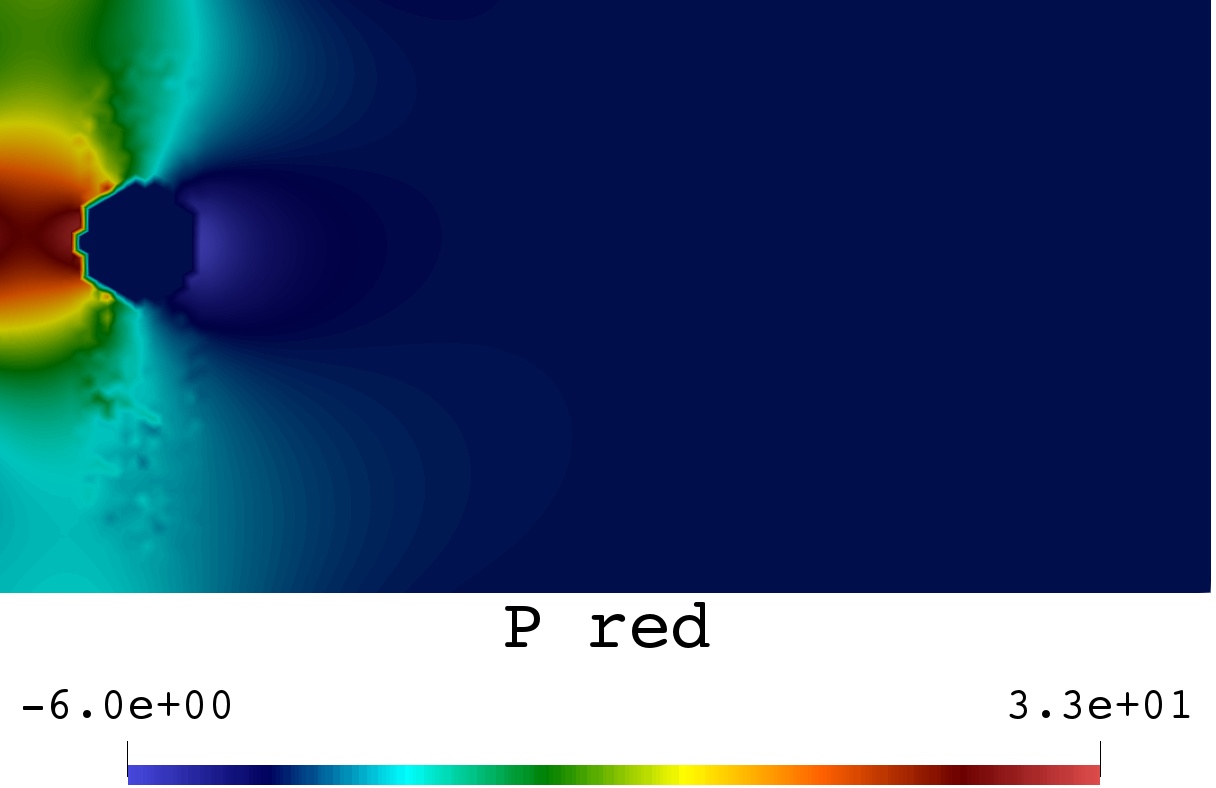}
\end{minipage}
\begin{minipage}{0.24\textwidth}
  \includegraphics[width=\textwidth]{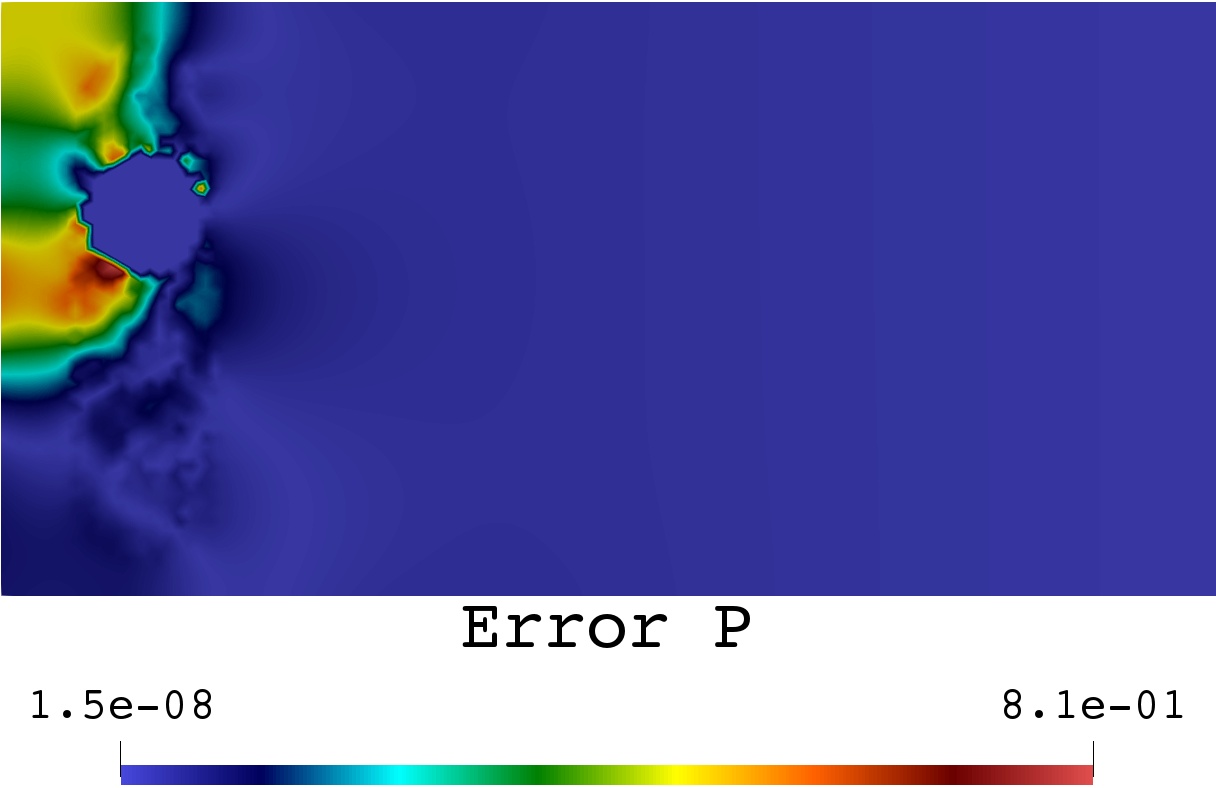}
\end{minipage}
\begin{minipage}{0.24\textwidth}
  \includegraphics[width=\textwidth]{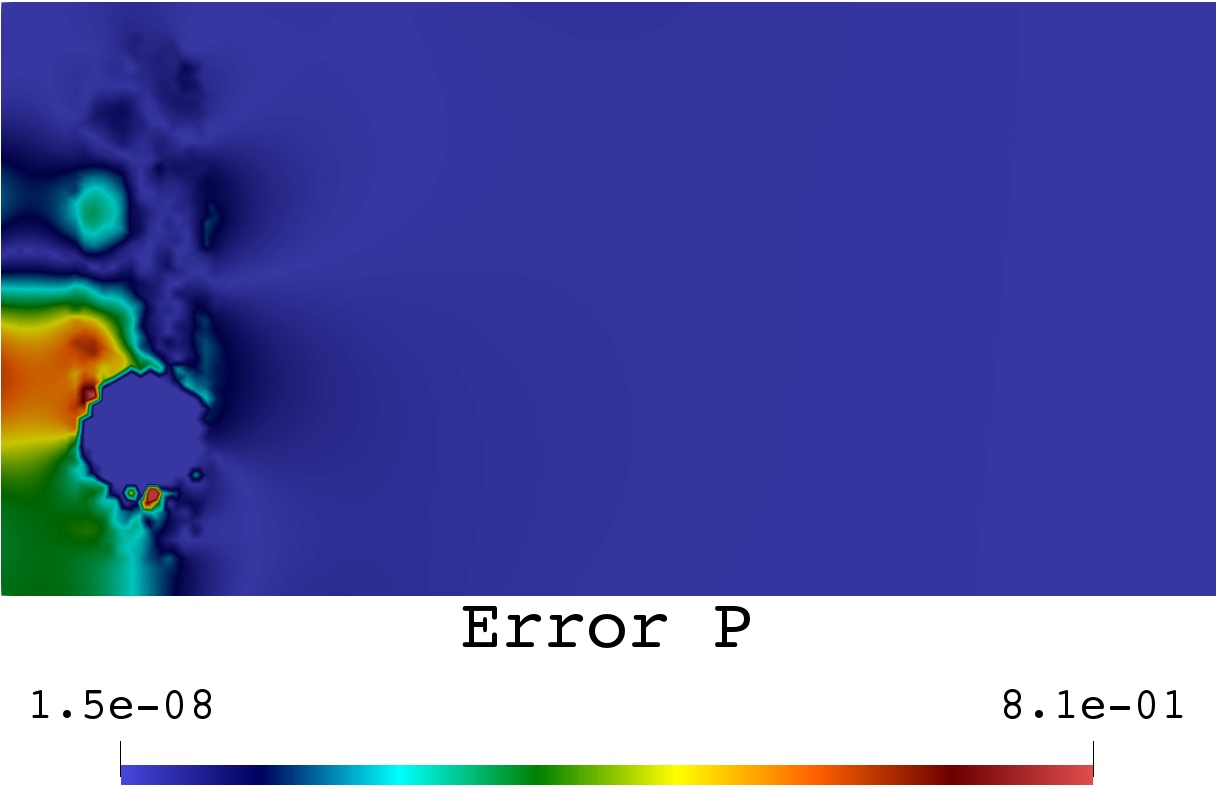}
\end{minipage}
\begin{minipage}{0.24\textwidth}
  \includegraphics[width=\textwidth]{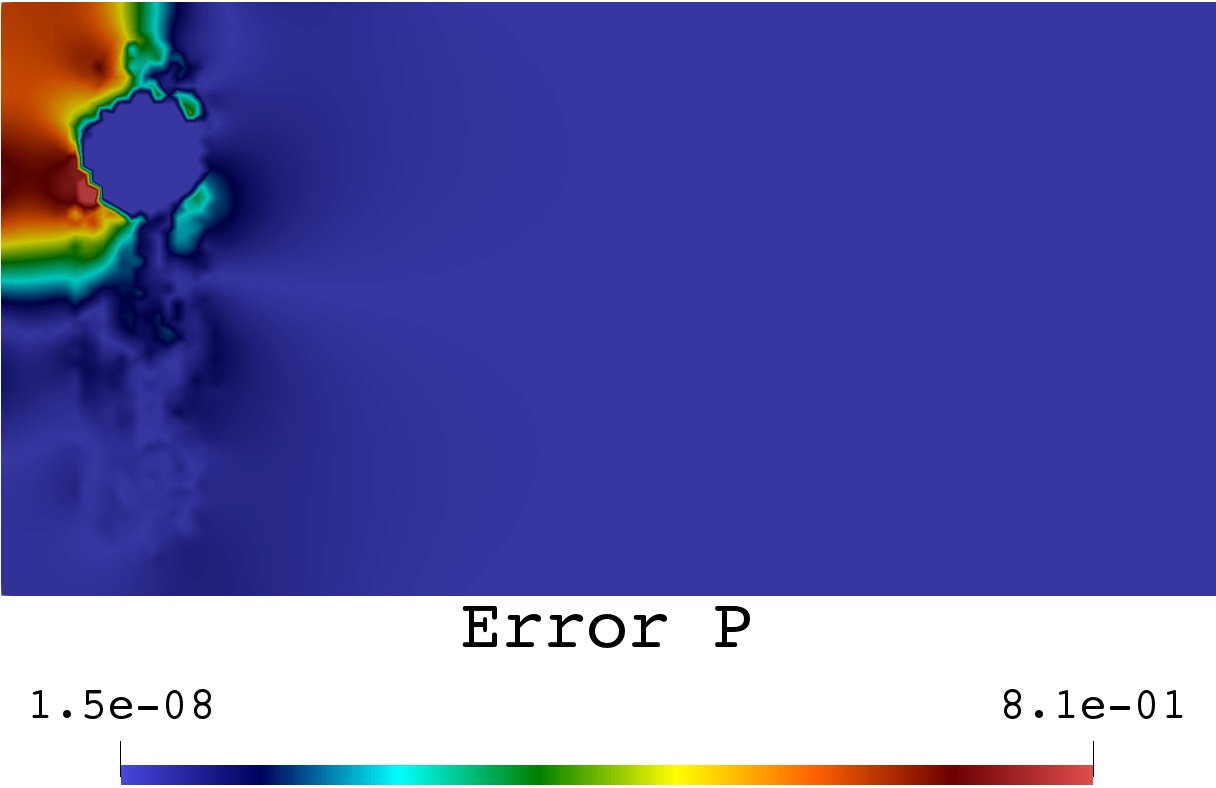}
\end{minipage}
\begin{minipage}{0.24\textwidth}
  \includegraphics[width=\textwidth]{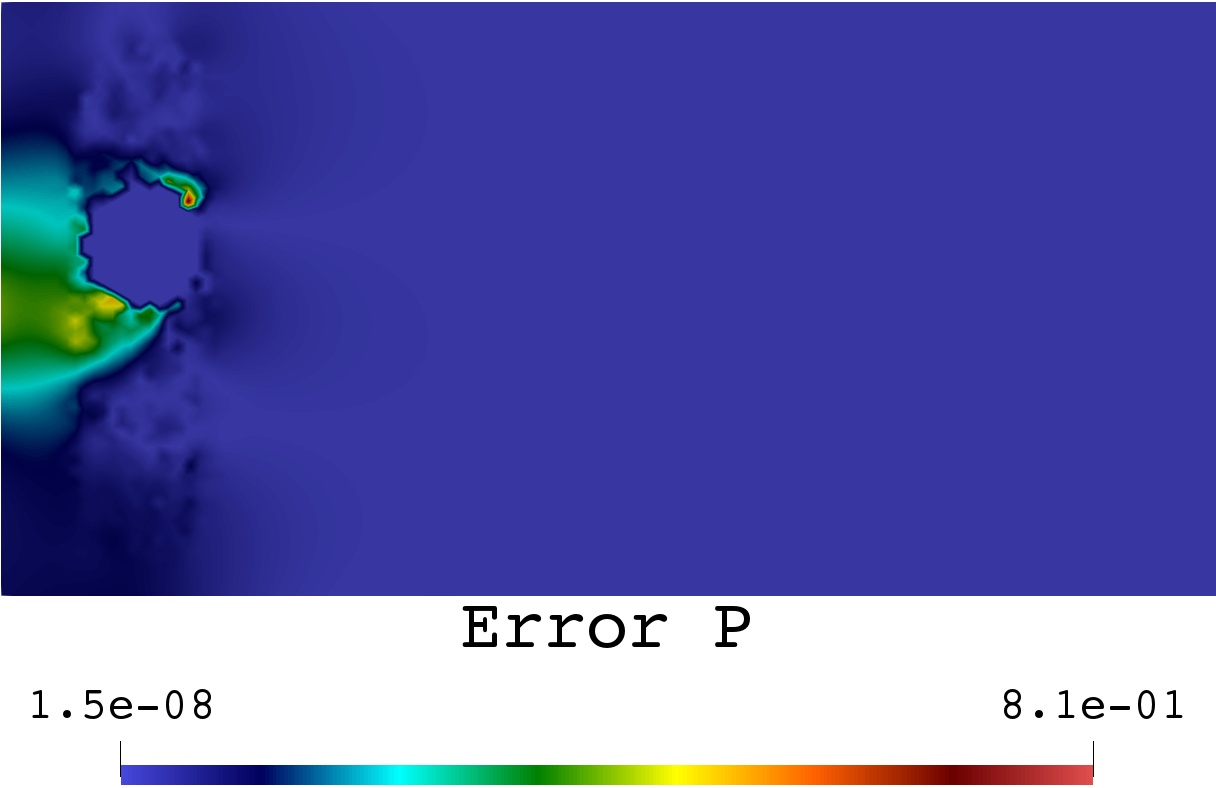}
\end{minipage}
\end{minipage}
\caption{Results for the 1D geometrical parametrization with $\mu_1 \in [-0.65, 0.65]$ {\blue and supremizer} basis enrichment. In rows $1-3$ we report the full-order solution, the reduced order solution and the absolute error plots for the velocity field while in rows $4-6$ we report the same quantities for the pressure field. The different columns are for four different values of the input parameter, $\mu_1 = [0.2876,-0.4520,0.4849,0.1780]$, while $\mu_0$ is fixed to the value $-1.5$.}
  \label{FULL_RED_ERROR_1P}
\end{figure}

In \autoref{FULL_RED_ERROR_1P} we report also the visualization of the full-order solutions, the reduced-order solutions and the absolute error for both the velocity and the pressure field in the case with a supremizer stabilization approach. The improved results for pressure are obvious and the plots clearly show at a glance that the FOM and ROM solutions cannot be easily distinguished. 
\subsubsection{Execution times}
In this subsection we analyze the execution times of the online stage and we compare them against the full-order computational times. Indicatively the computational time necessary to compute $1024$ solves of the full order problem is equal to $1$ hour $43$ minutes.
The offline stage is usually very expensive but in a reduced order modeling framework, fortunately needs to be performed only once. In \autoref{table:timers} we report the ``online'' computational time for different dimensions of the reduced order model. This consists of the time necessary for the assembly of the full-order matrices, the generation of the reduced order model and its resolution. 

\begin{figure} 
\centering
\begin{minipage}{\textwidth}
\centering
\begin{minipage}{0.24\textwidth}
  \includegraphics[width=\textwidth]{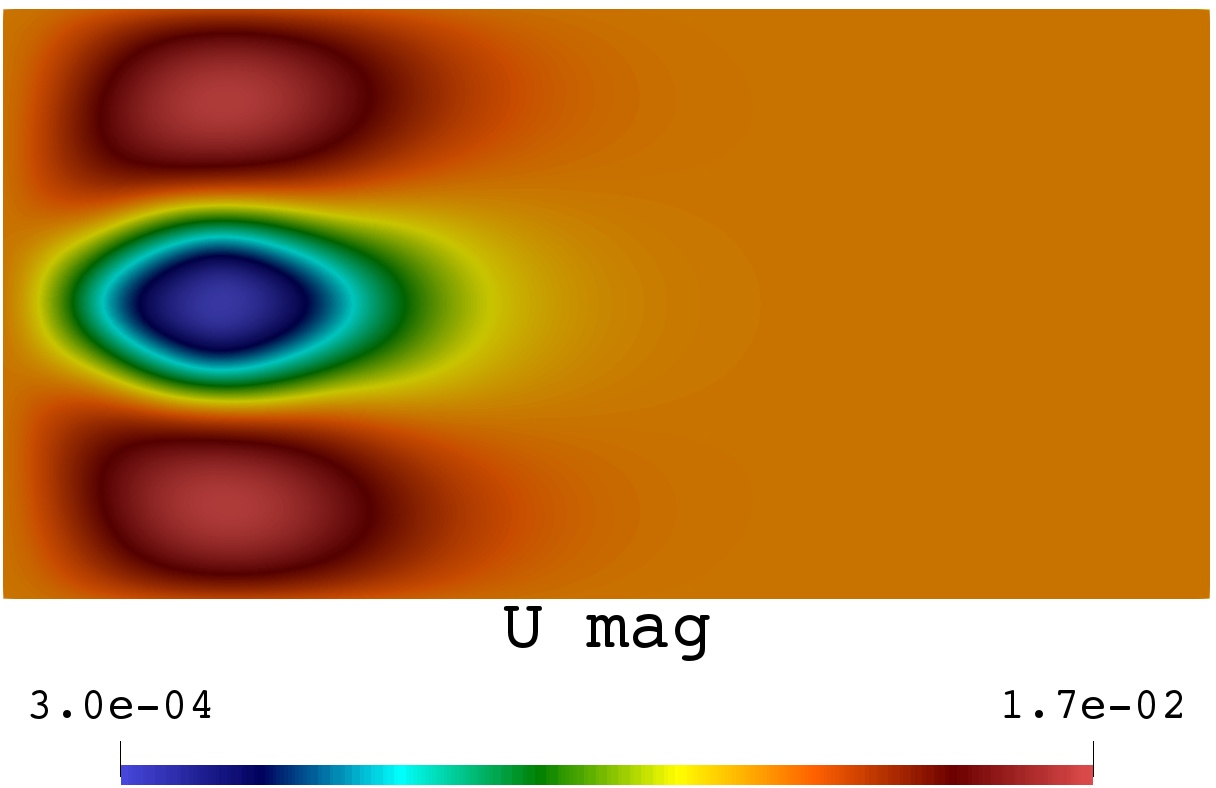}
\end{minipage}
\begin{minipage}{0.24\textwidth}
  \includegraphics[width=\textwidth]{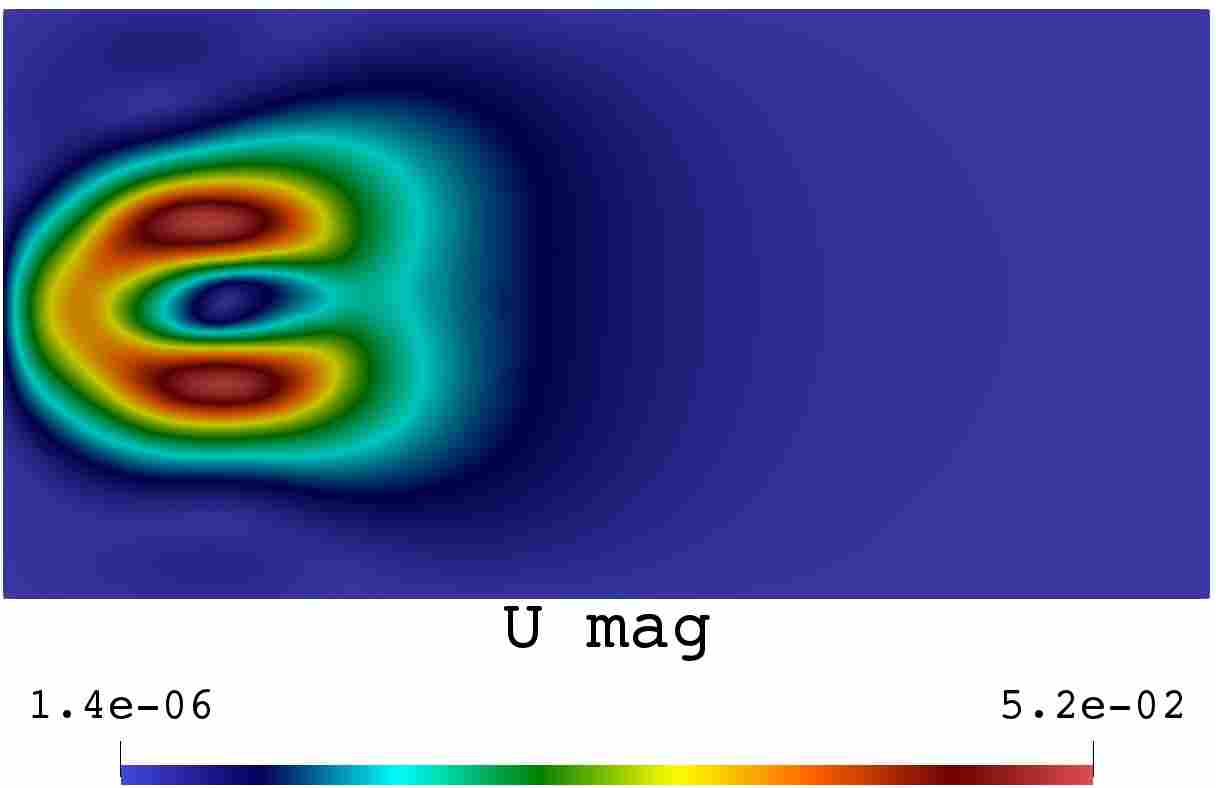}
\end{minipage}
\begin{minipage}{0.24\textwidth}
  \includegraphics[width=\textwidth]{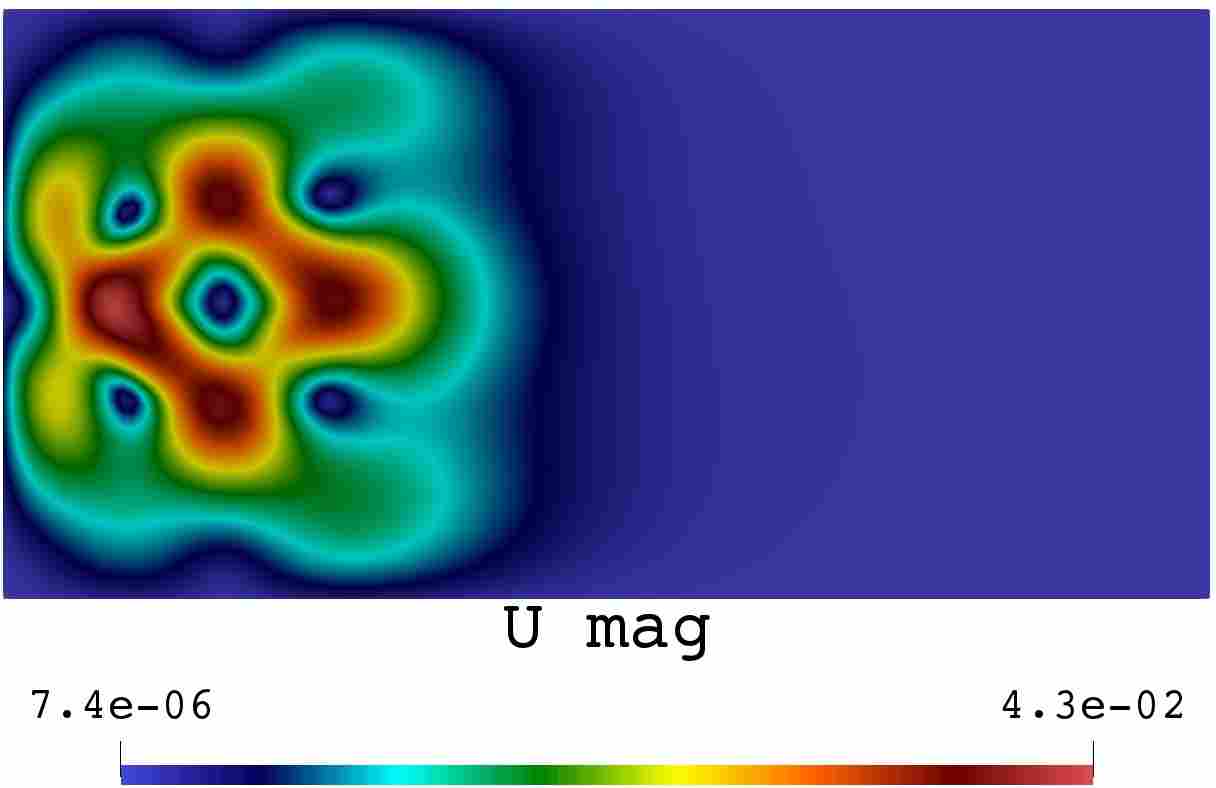}
\end{minipage}
\begin{minipage}{0.24\textwidth}
  \includegraphics[width=\textwidth]{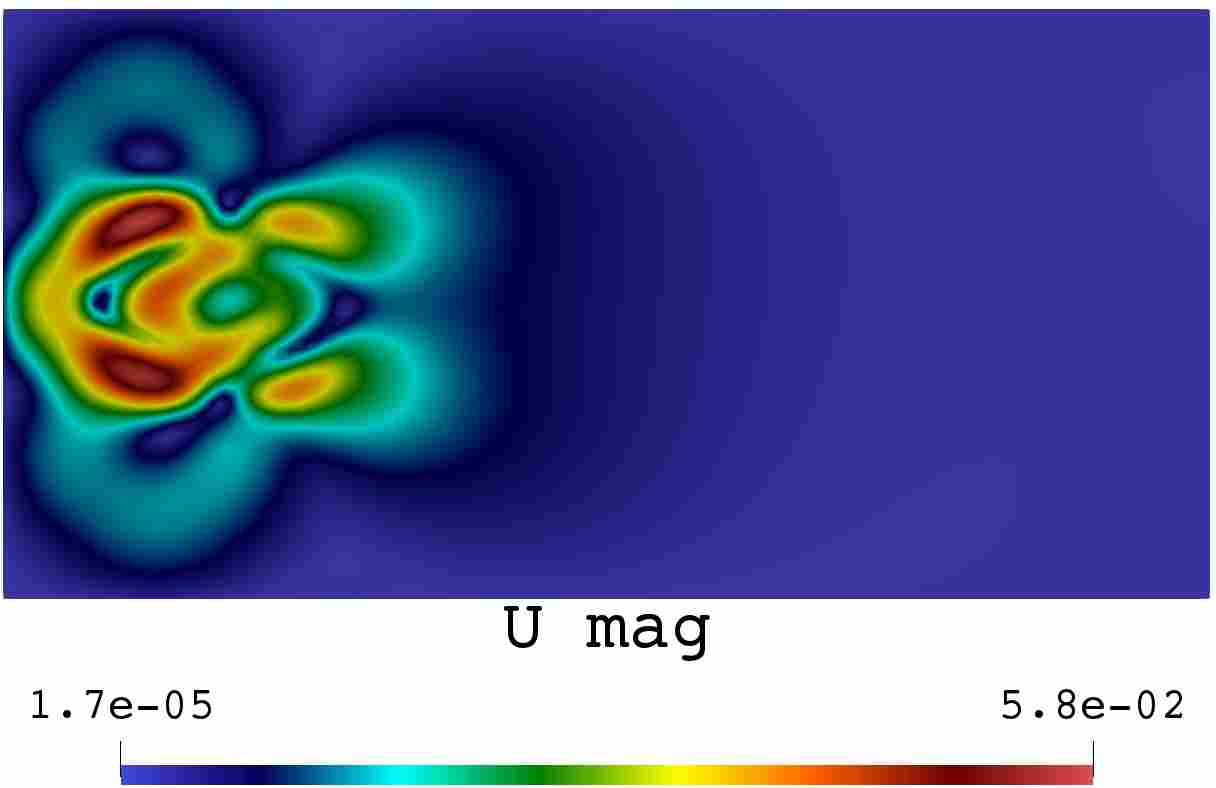}
\end{minipage}
\begin{minipage}{0.24\textwidth}
  \includegraphics[width=\textwidth]{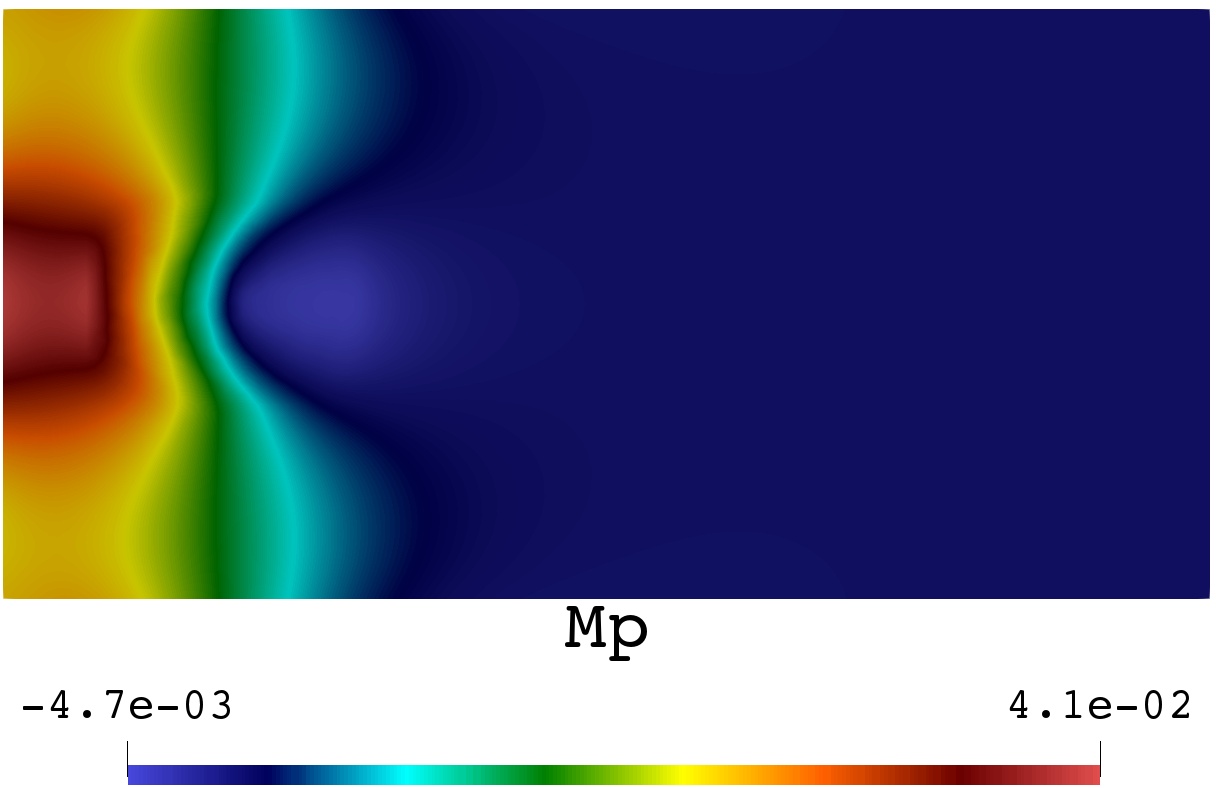}
\end{minipage}
\begin{minipage}{0.24\textwidth}
  \includegraphics[width=\textwidth]{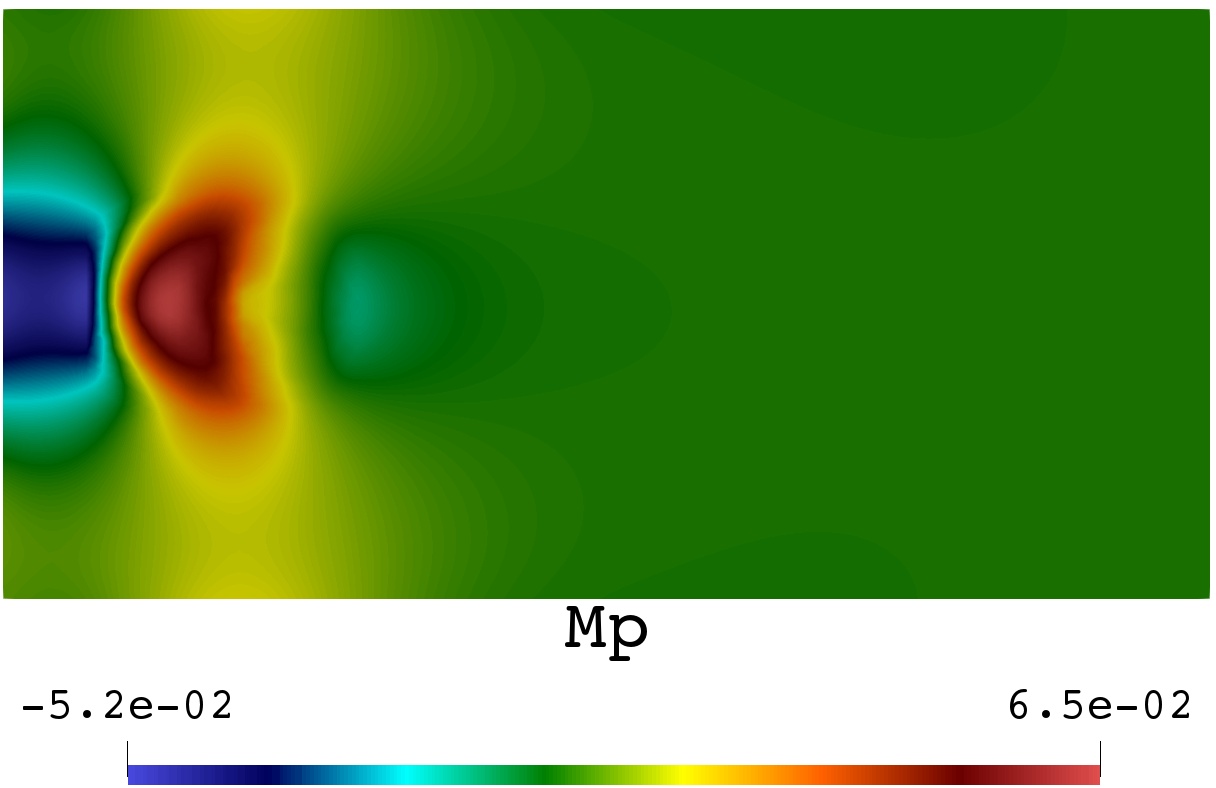}
\end{minipage}
\begin{minipage}{0.24\textwidth}
  \includegraphics[width=\textwidth]{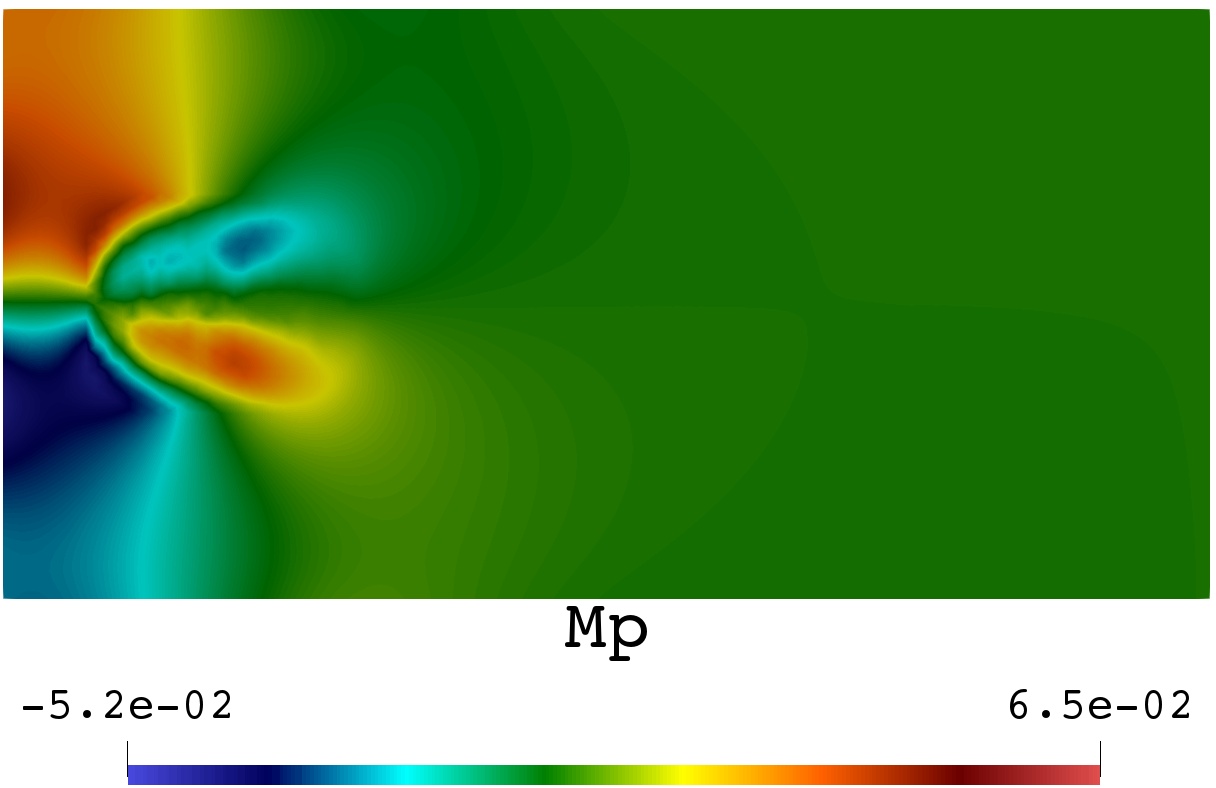}
\end{minipage}
\begin{minipage}{0.24\textwidth}
  \includegraphics[width=\textwidth]{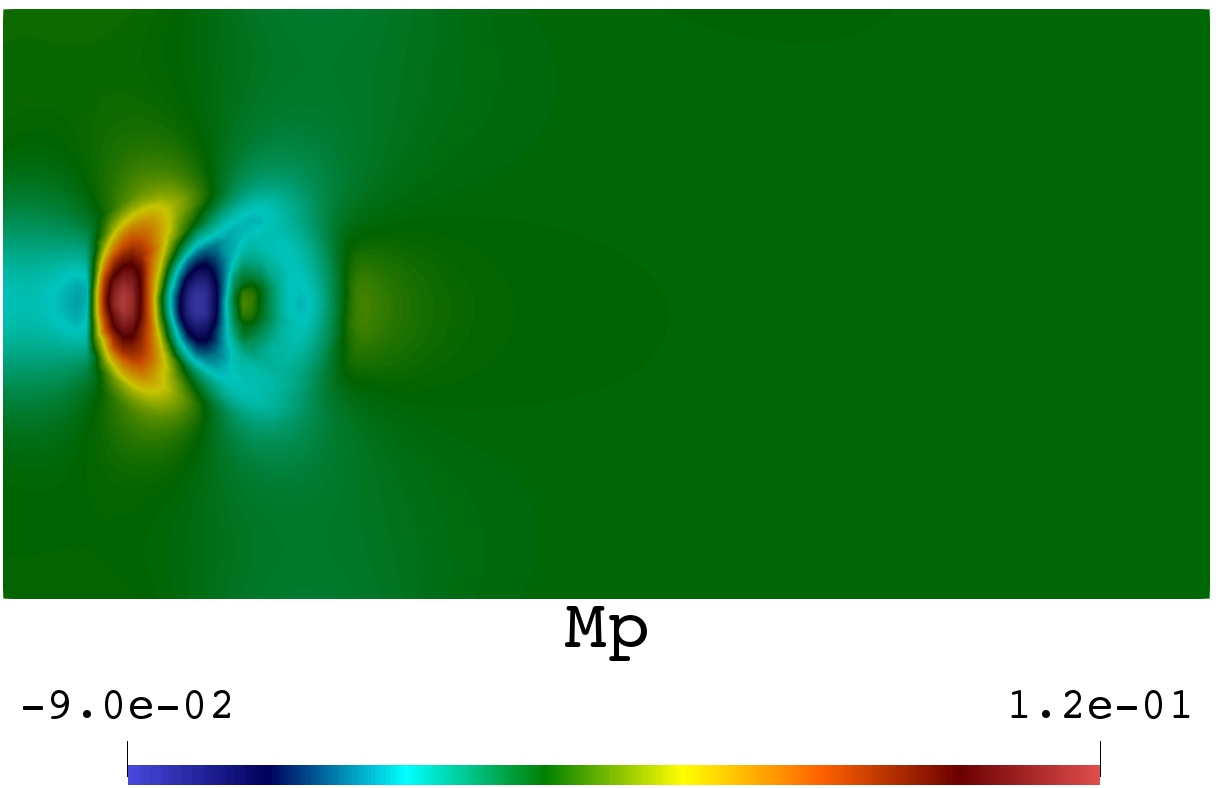}
\end{minipage}
\end{minipage}
  \caption{The first $4$ basis components for velocity (first row) and pressure (second row) in the 2D geometrical parametrization case with $\mu \in [-1.5, -1.0]\times[-0.15, 0.15]$.}
  \label{Stokes_Pressure_and_Velocity_Components_Modes_2P}
\end{figure} 
\begin{figure} 
\centering
\begin{minipage}{\textwidth}
\centering
\begin{minipage}{0.24\textwidth}
  \includegraphics[width=\textwidth]{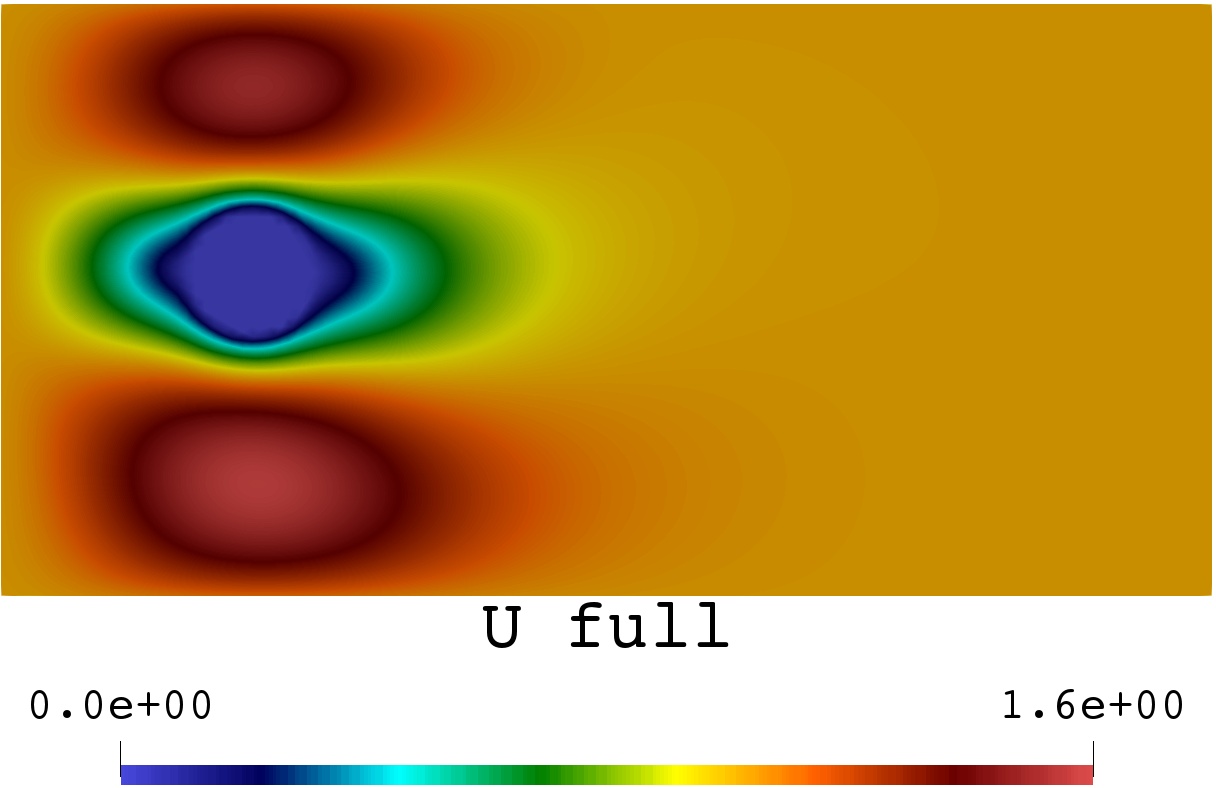}
\end{minipage}
\begin{minipage}{0.24\textwidth}
  \includegraphics[width=\textwidth]{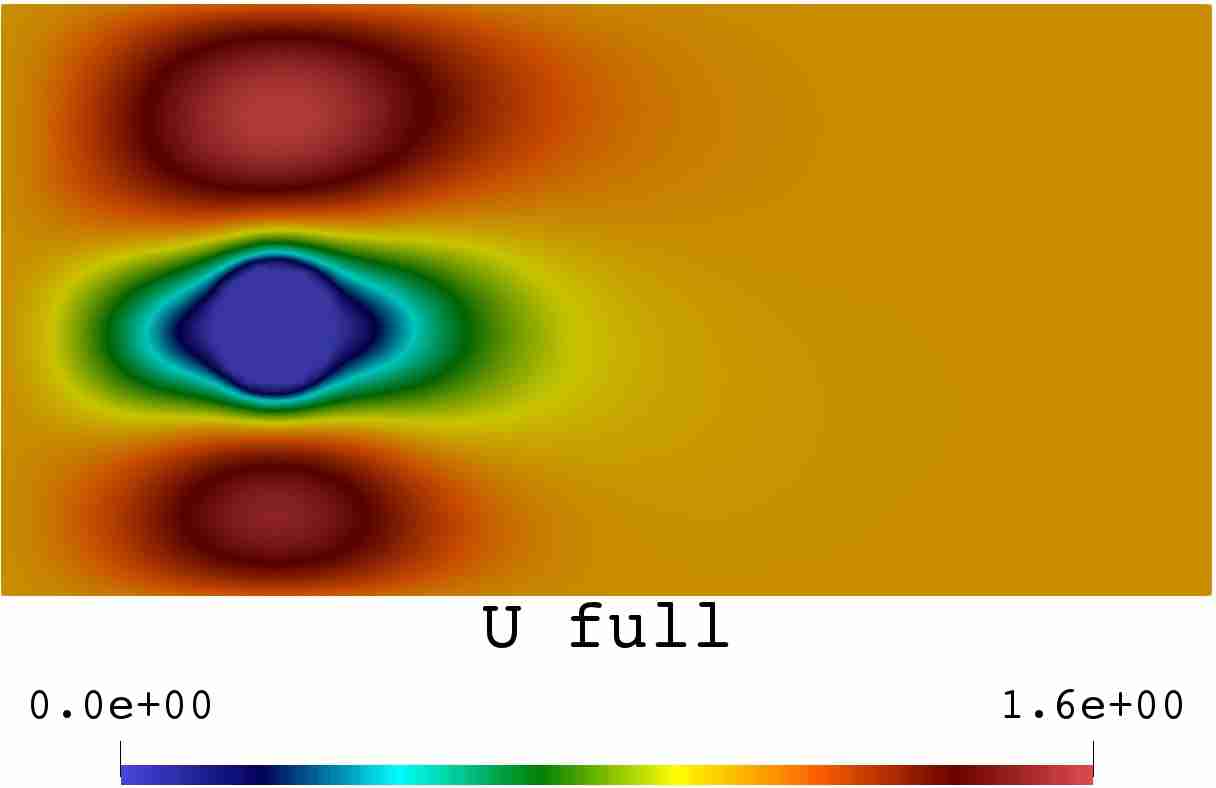}
\end{minipage}
\begin{minipage}{0.24\textwidth}
  \includegraphics[width=\textwidth]{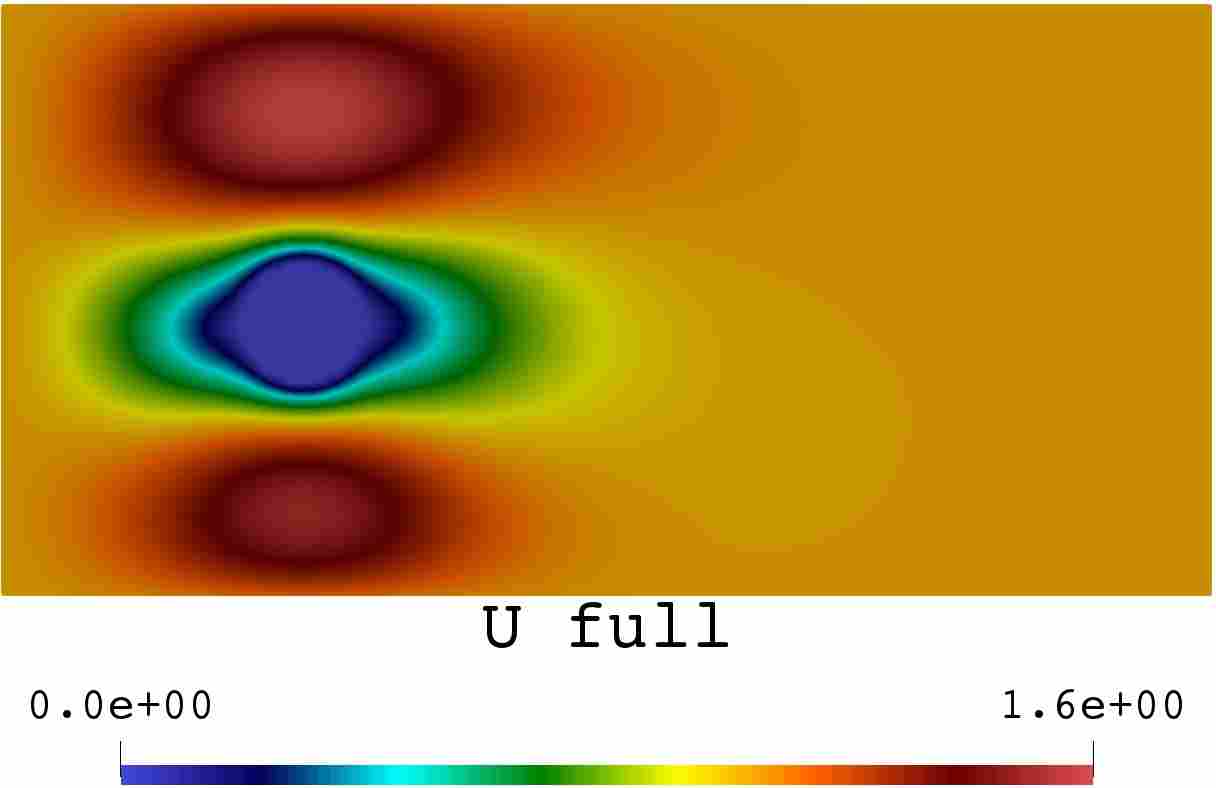}
\end{minipage}
\begin{minipage}{0.24\textwidth}
  \includegraphics[width=\textwidth]{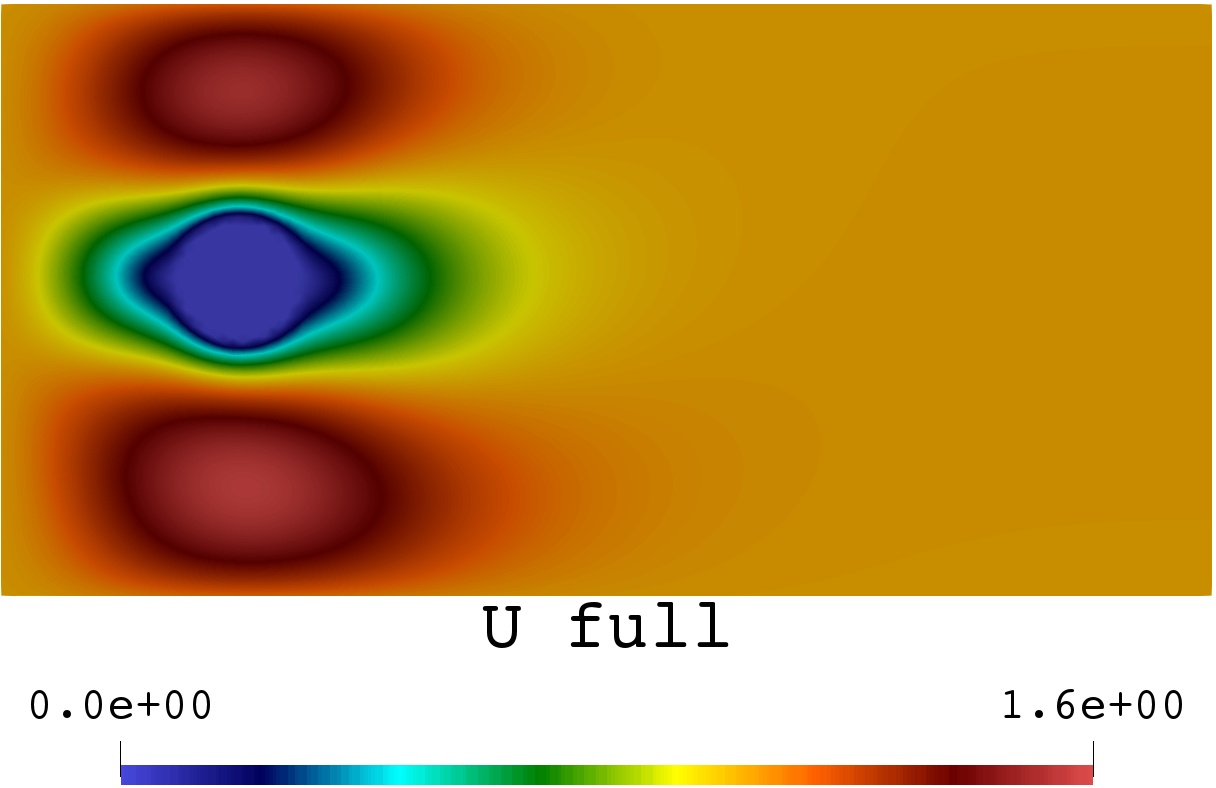}
\end{minipage}
\begin{minipage}{0.24\textwidth}
  \includegraphics[width=\textwidth]{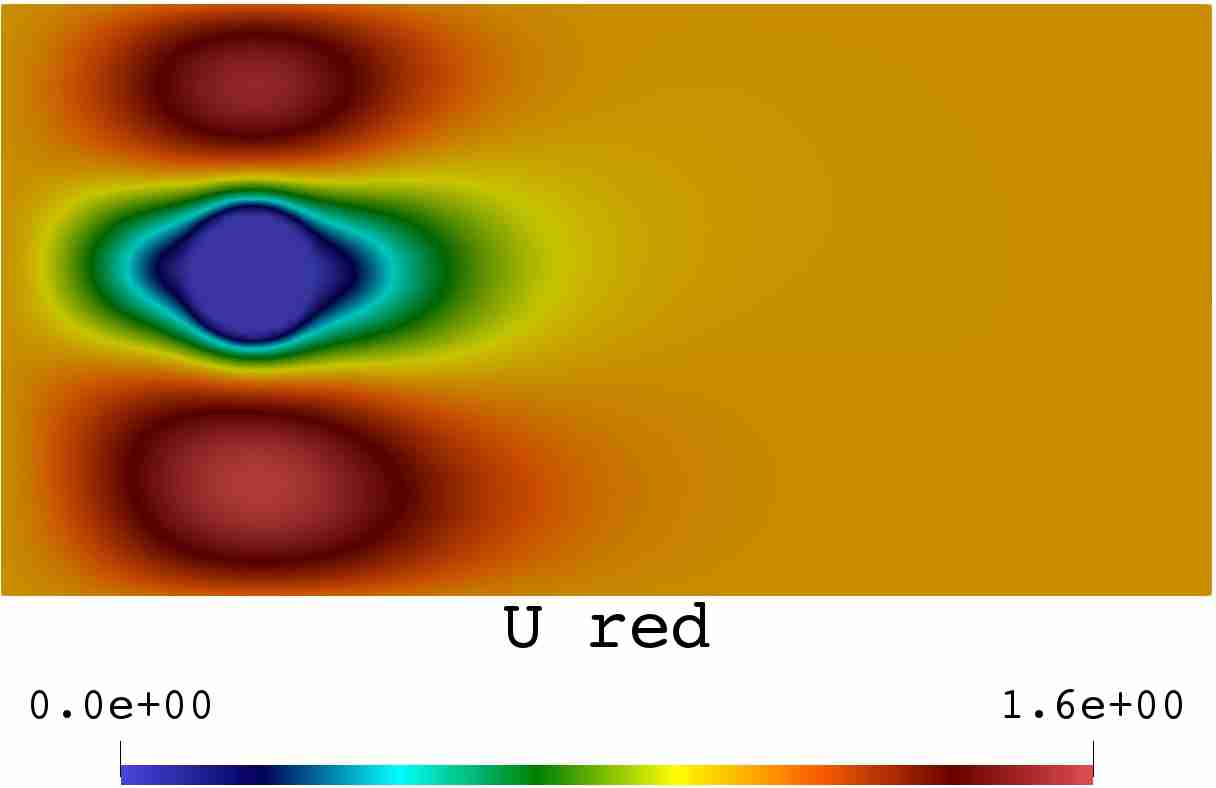}
\end{minipage}
\begin{minipage}{0.24\textwidth}
  \includegraphics[width=\textwidth]{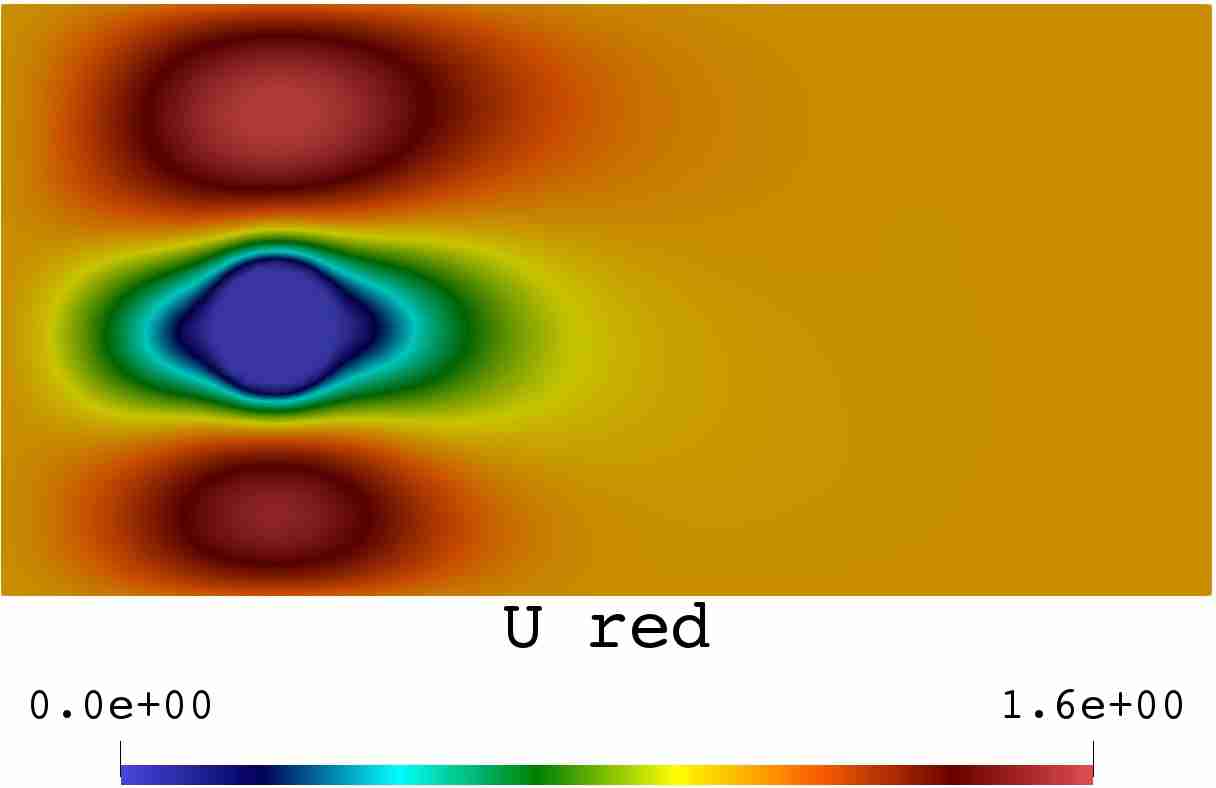}
\end{minipage}
\begin{minipage}{0.24\textwidth}
  \includegraphics[width=\textwidth]{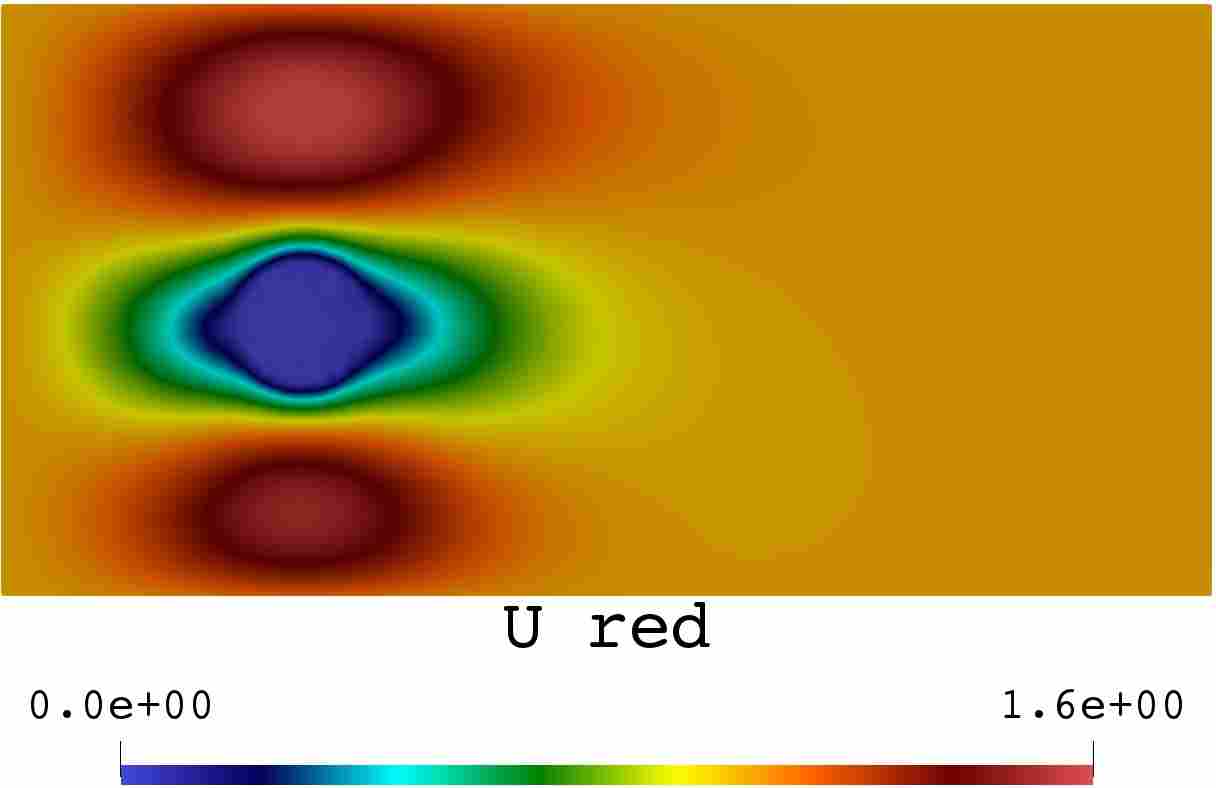}
\end{minipage}
\begin{minipage}{0.24\textwidth}
  \includegraphics[width=\textwidth]{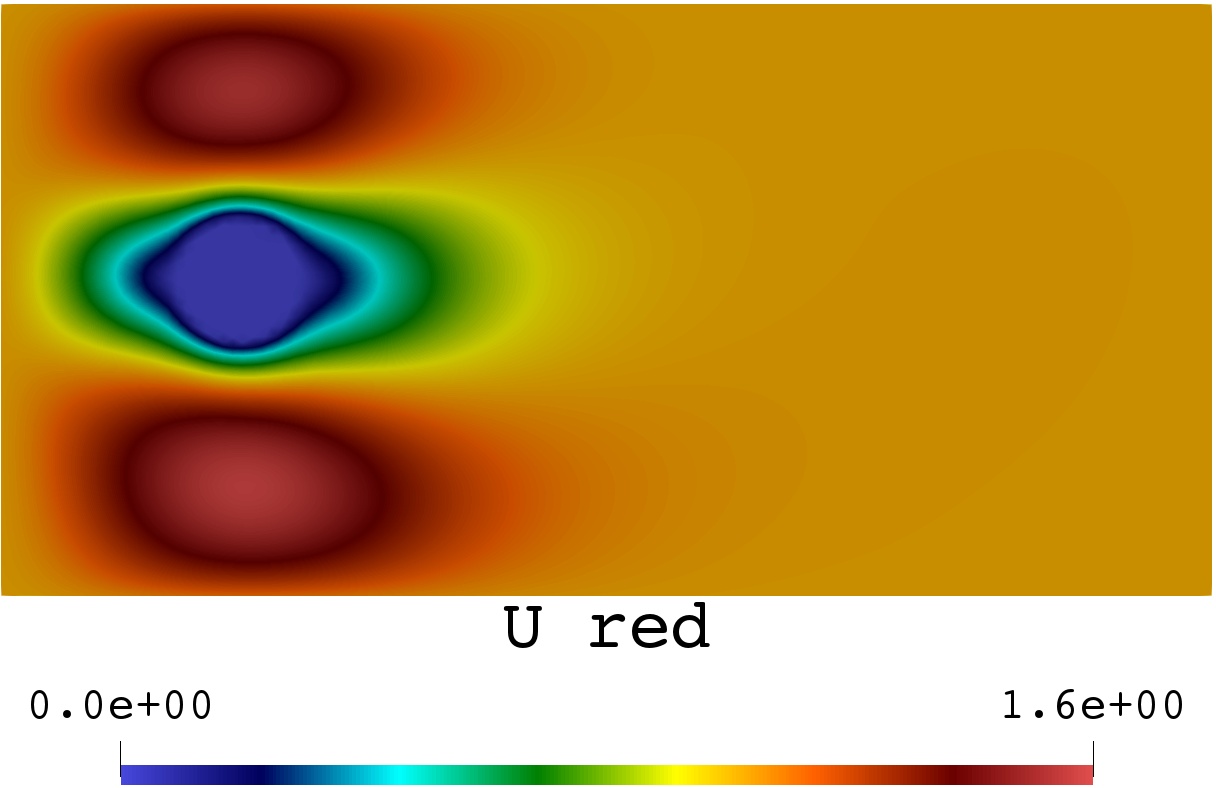}
\end{minipage}
\begin{minipage}{0.24\textwidth}
  \includegraphics[width=\textwidth]{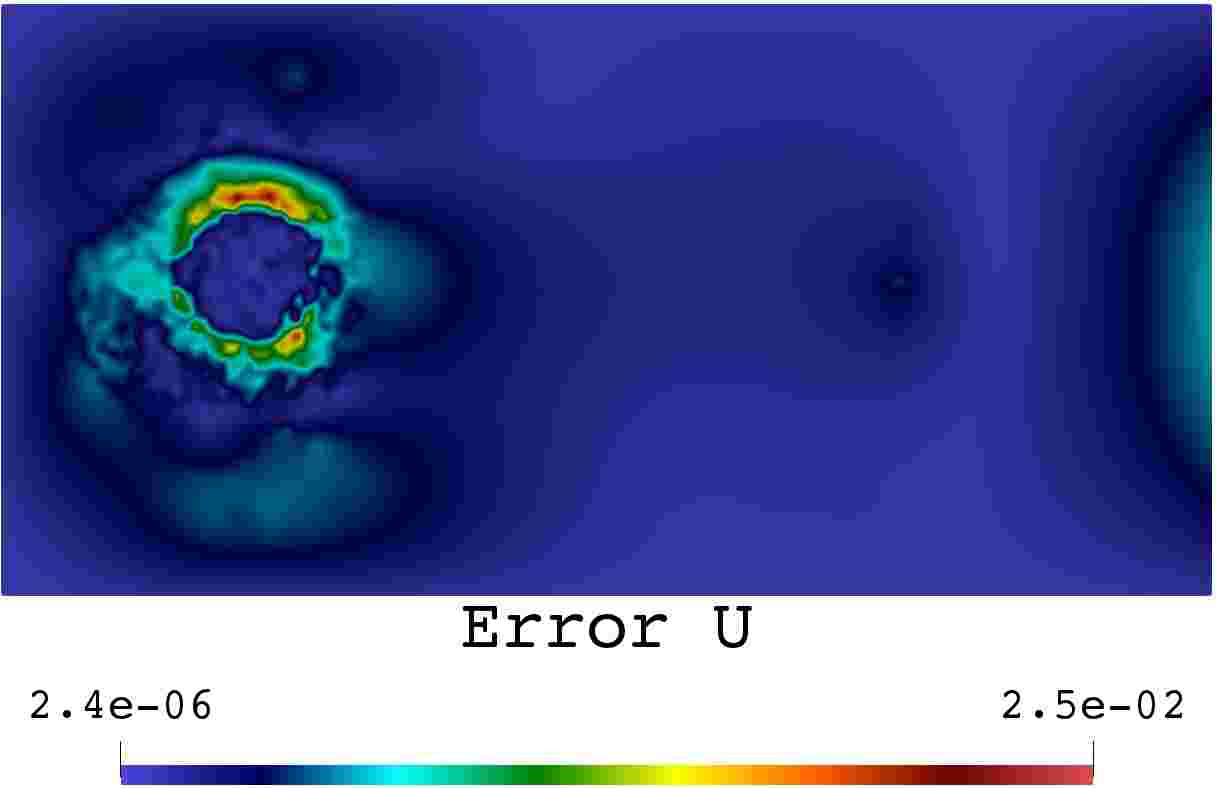}
\end{minipage}
\begin{minipage}{0.24\textwidth}
  \includegraphics[width=\textwidth]{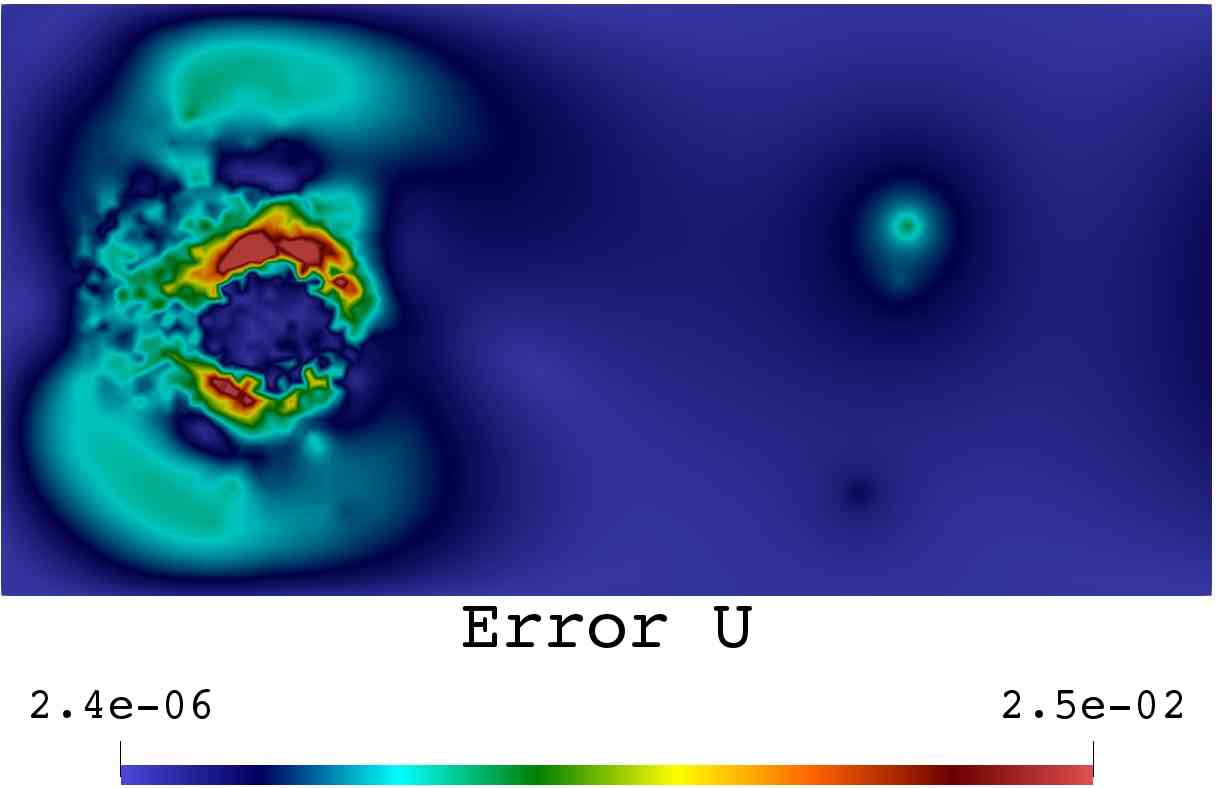}
\end{minipage}
\begin{minipage}{0.24\textwidth}
  \includegraphics[width=\textwidth]{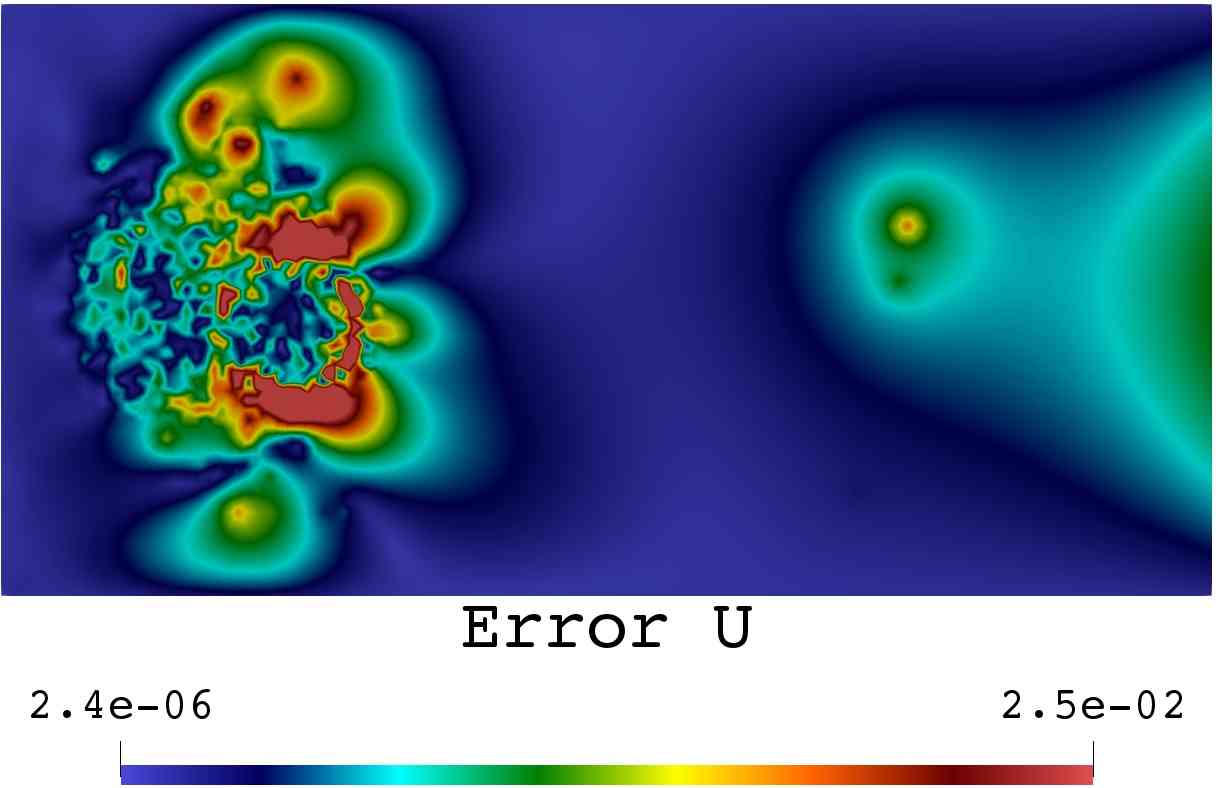}
\end{minipage}
\begin{minipage}{0.24\textwidth}
  \includegraphics[width=\textwidth]{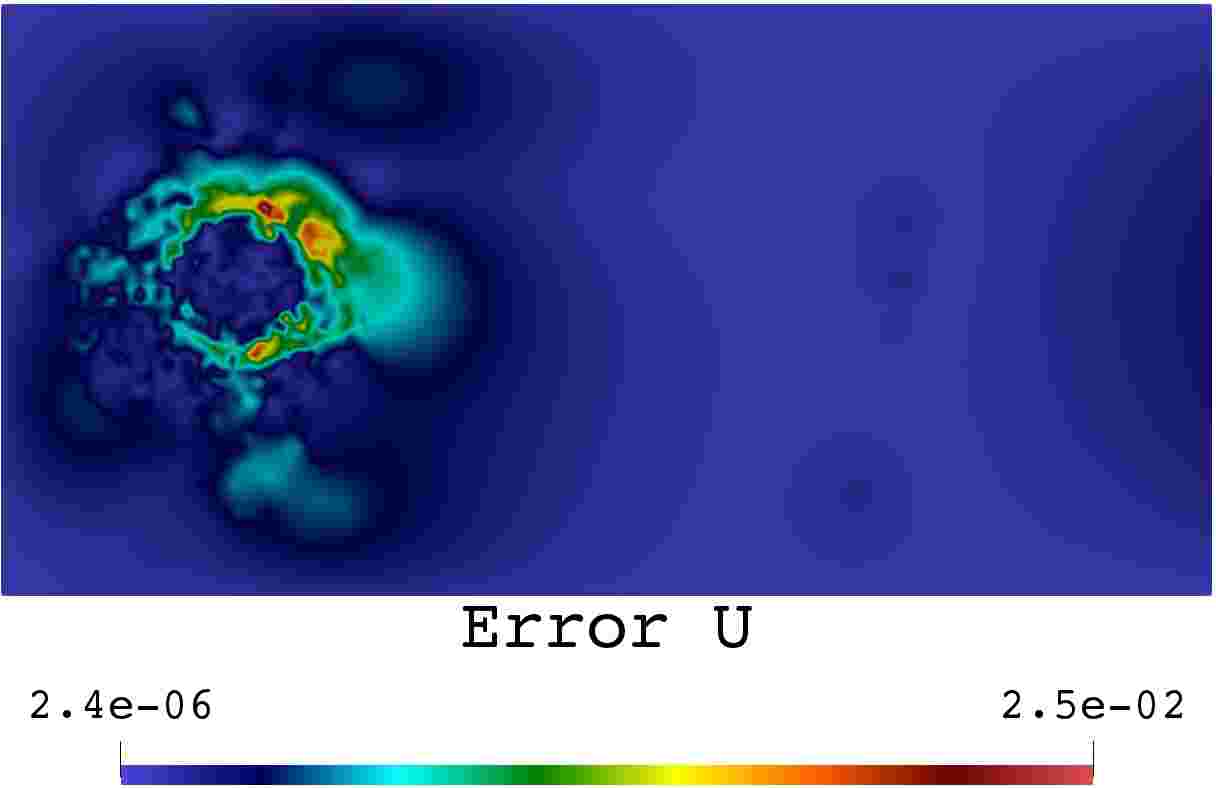}
\end{minipage}
\begin{minipage}{0.24\textwidth}
  \includegraphics[width=\textwidth]{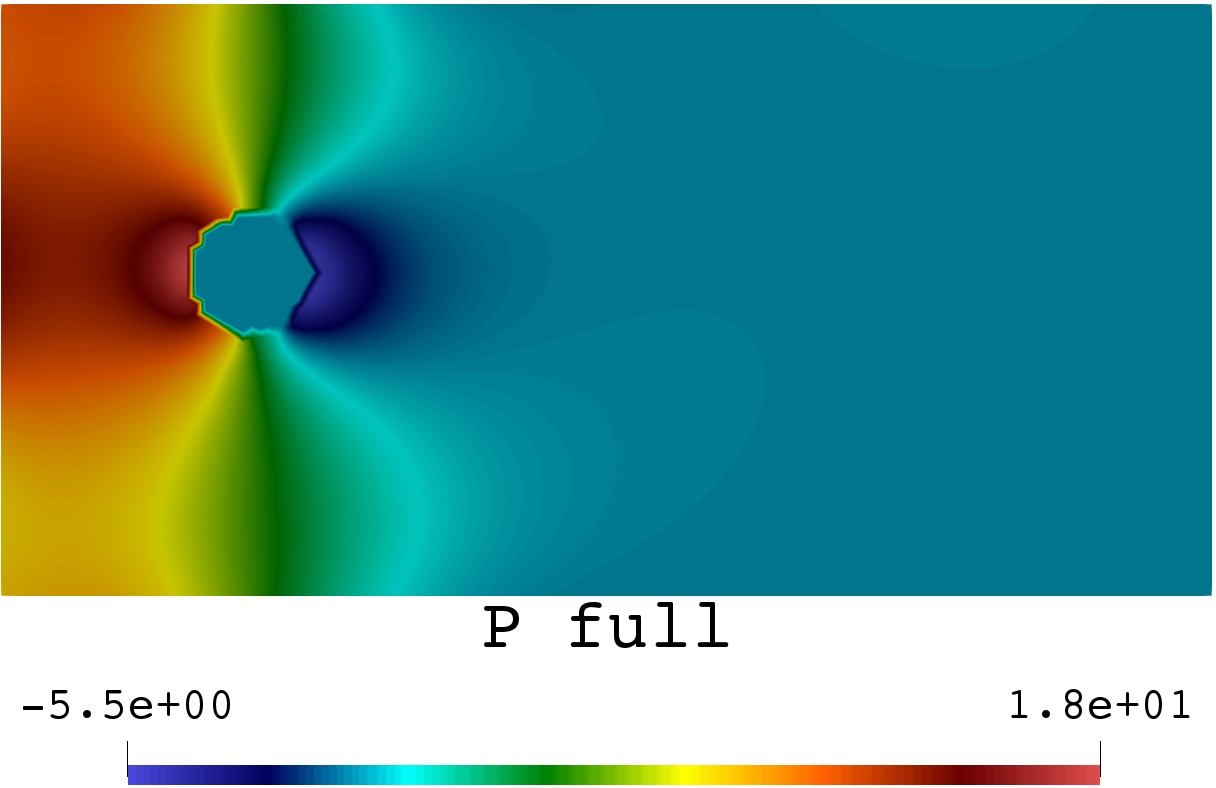}
\end{minipage}
\begin{minipage}{0.24\textwidth}
  \includegraphics[width=\textwidth]{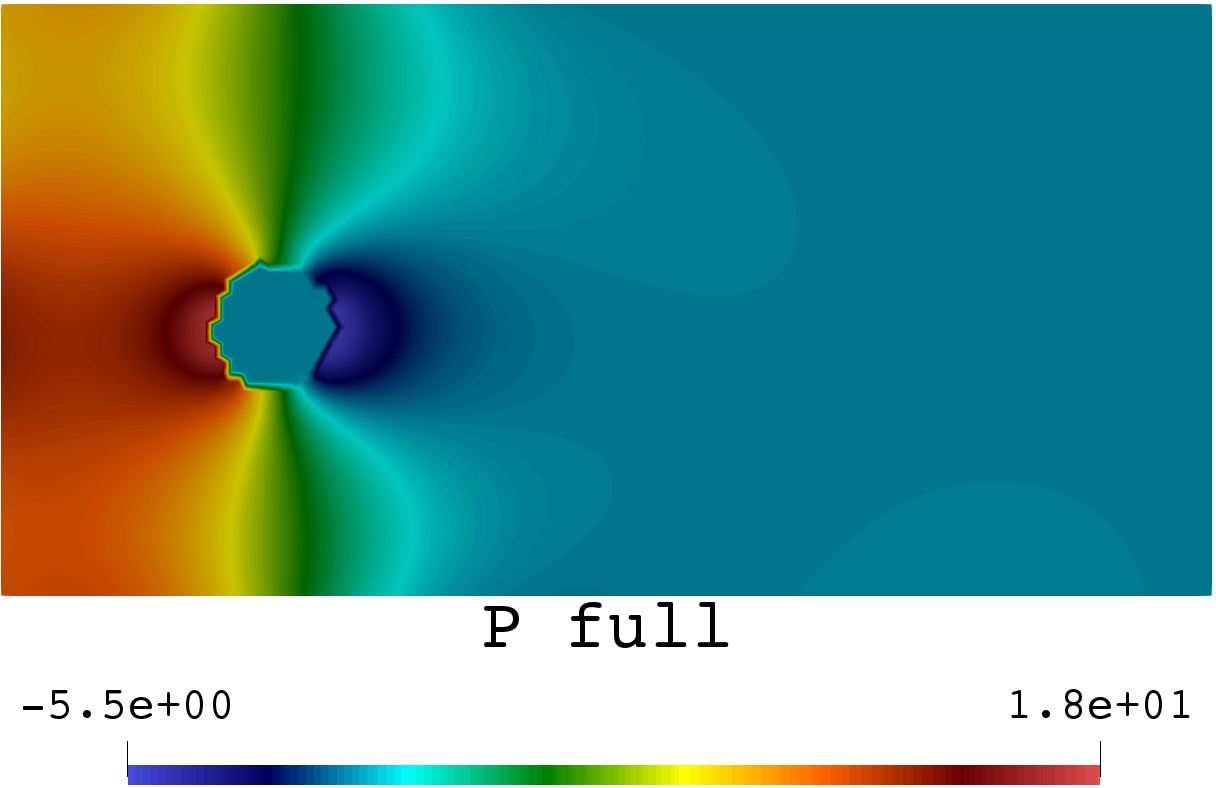}
\end{minipage}
\begin{minipage}{0.24\textwidth}
  \includegraphics[width=\textwidth]{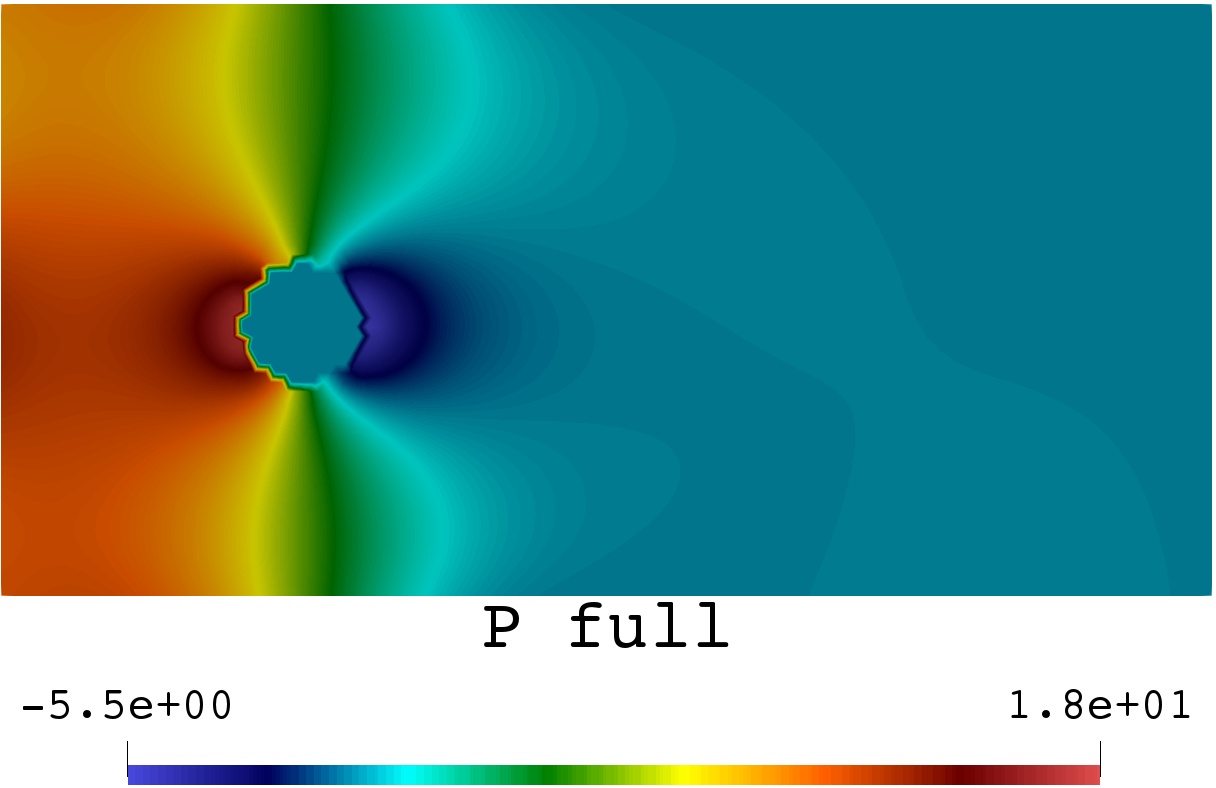}
\end{minipage}
\begin{minipage}{0.24\textwidth}
  \includegraphics[width=\textwidth]{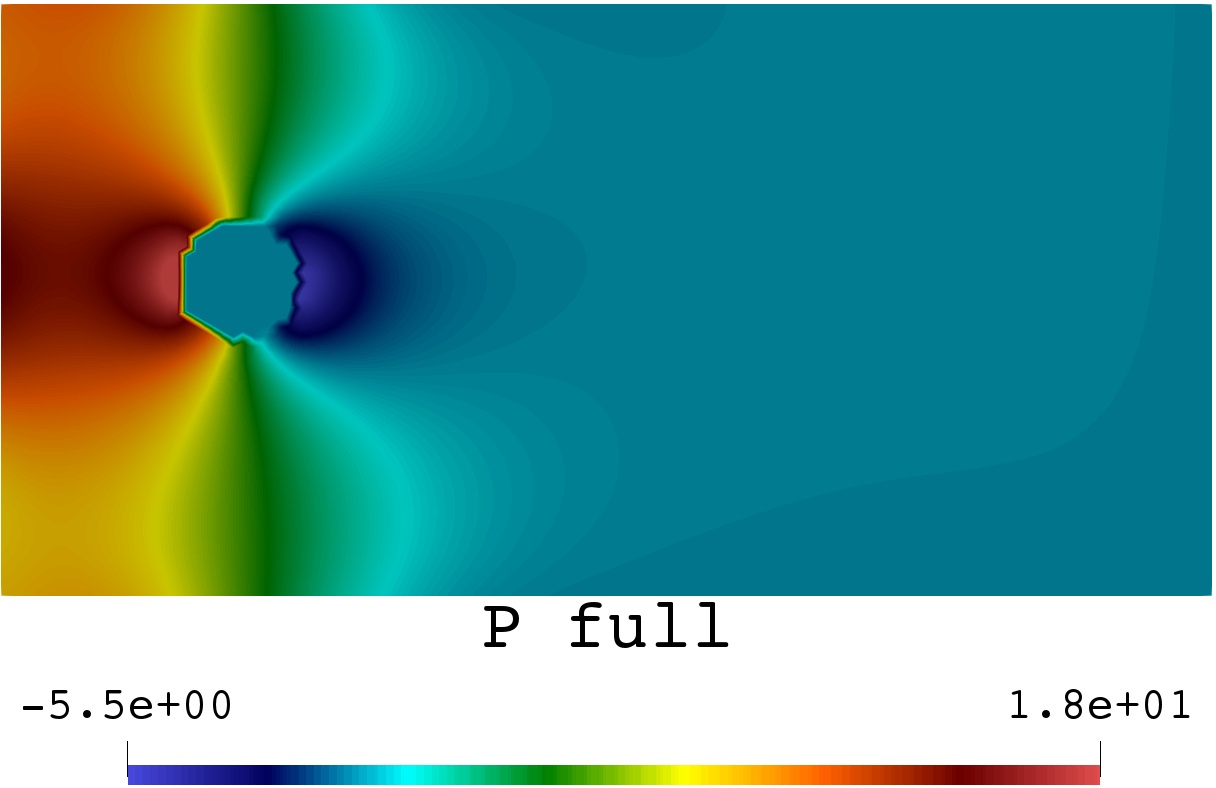}
\end{minipage}
\begin{minipage}{0.24\textwidth}
  \includegraphics[width=\textwidth]{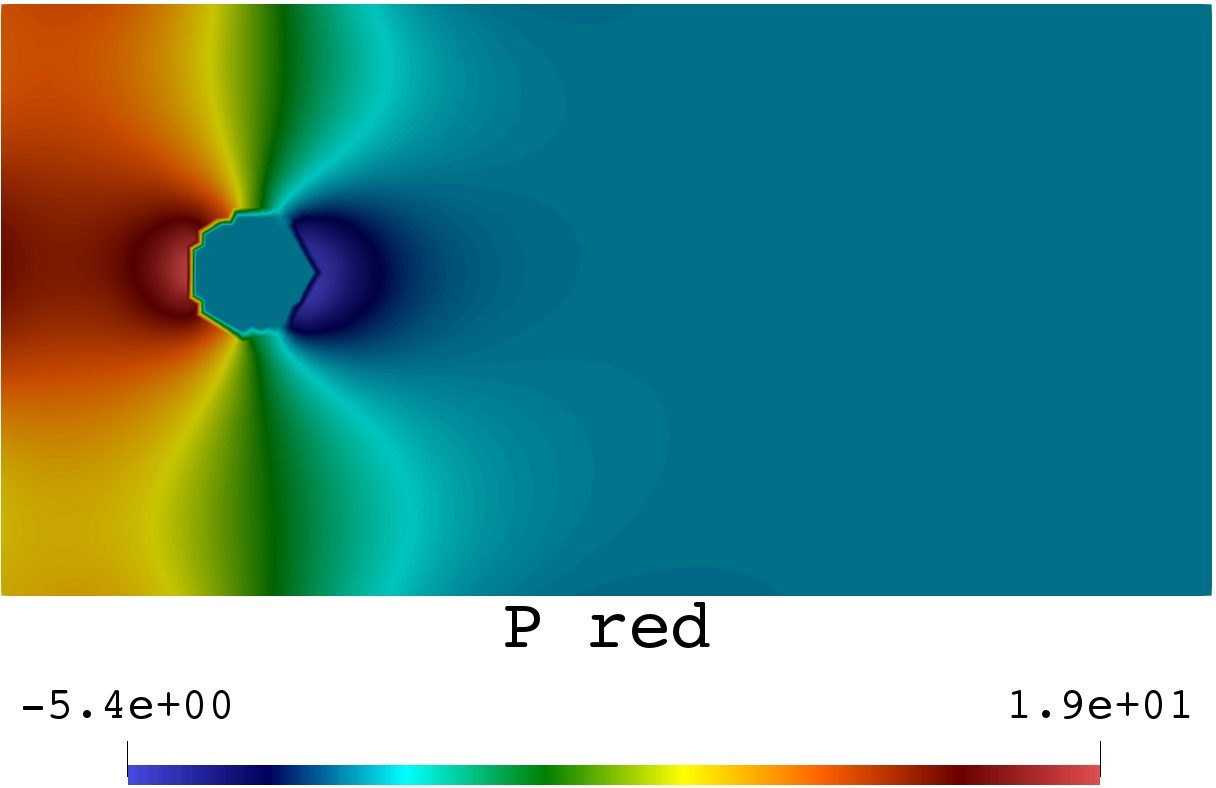}
\end{minipage}
\begin{minipage}{0.24\textwidth}
  \includegraphics[width=\textwidth]{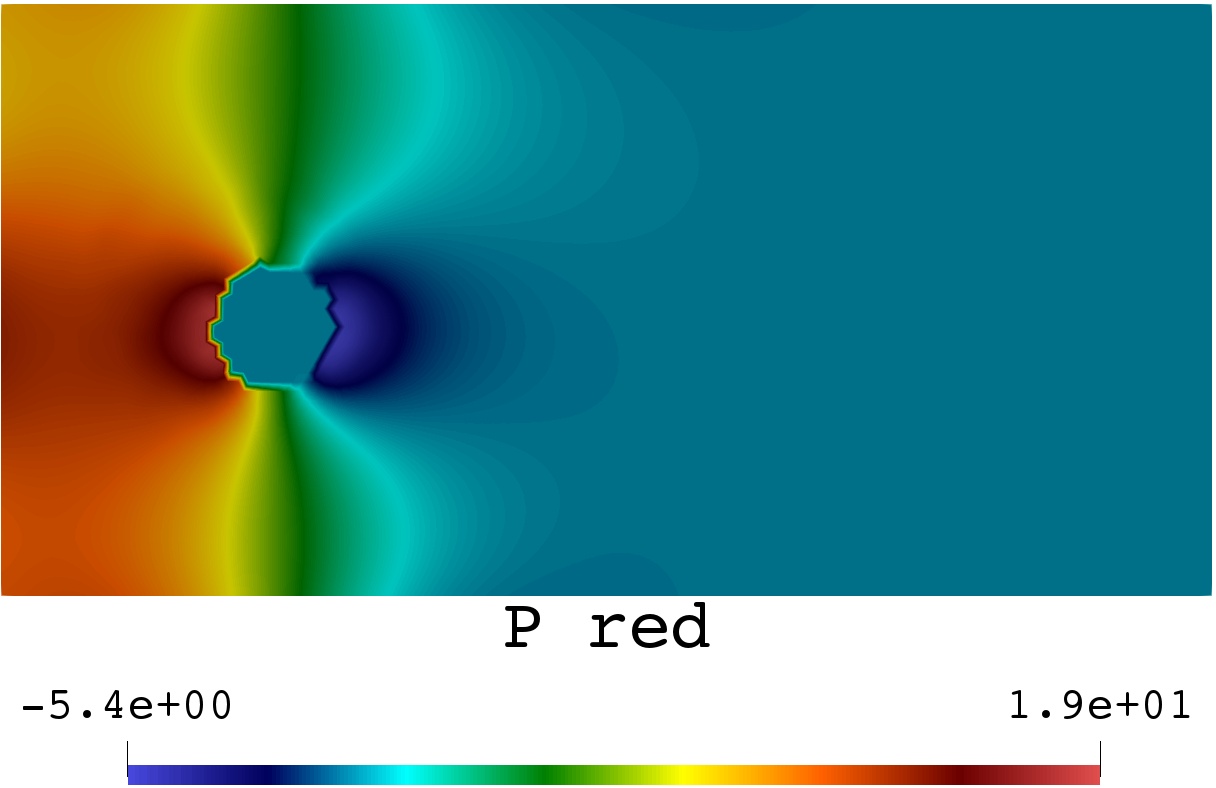}
\end{minipage}
\begin{minipage}{0.24\textwidth}
  \includegraphics[width=\textwidth]{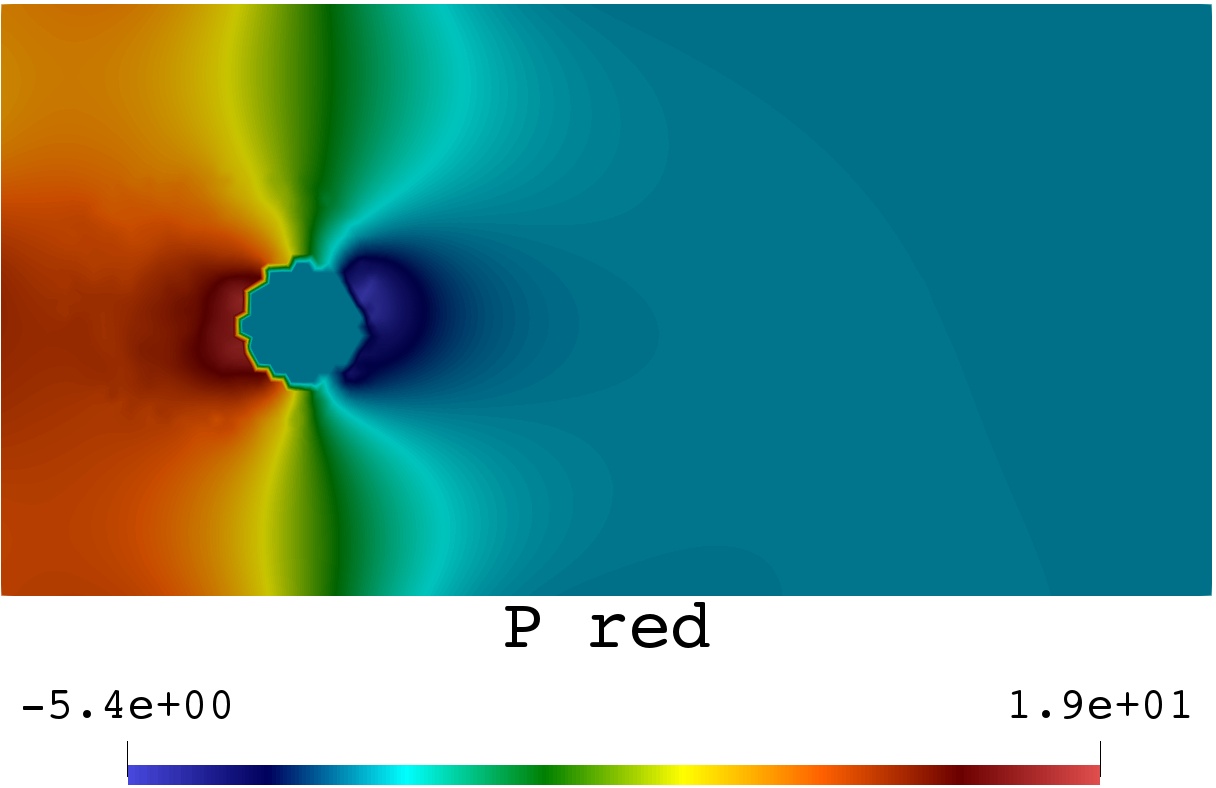}
\end{minipage}
\begin{minipage}{0.24\textwidth}
  \includegraphics[width=\textwidth]{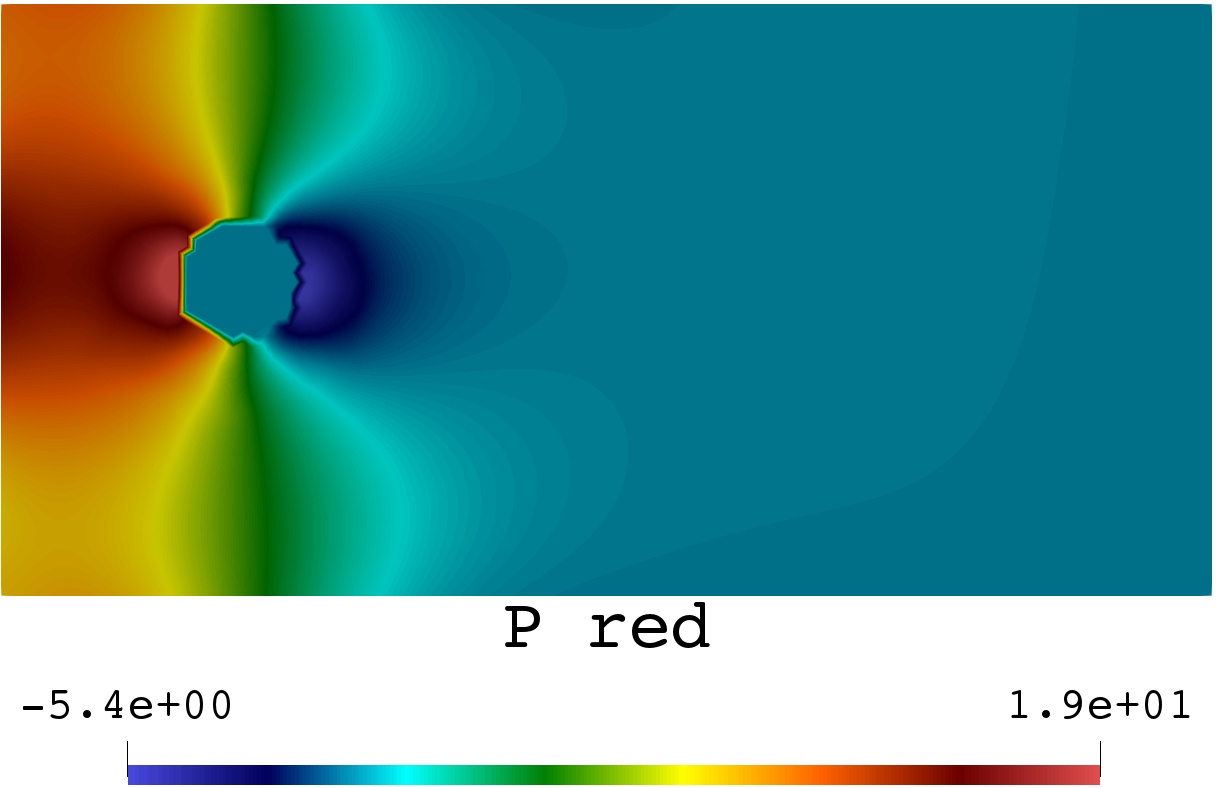}
\end{minipage}
\begin{minipage}{0.24\textwidth}
  \includegraphics[width=\textwidth]{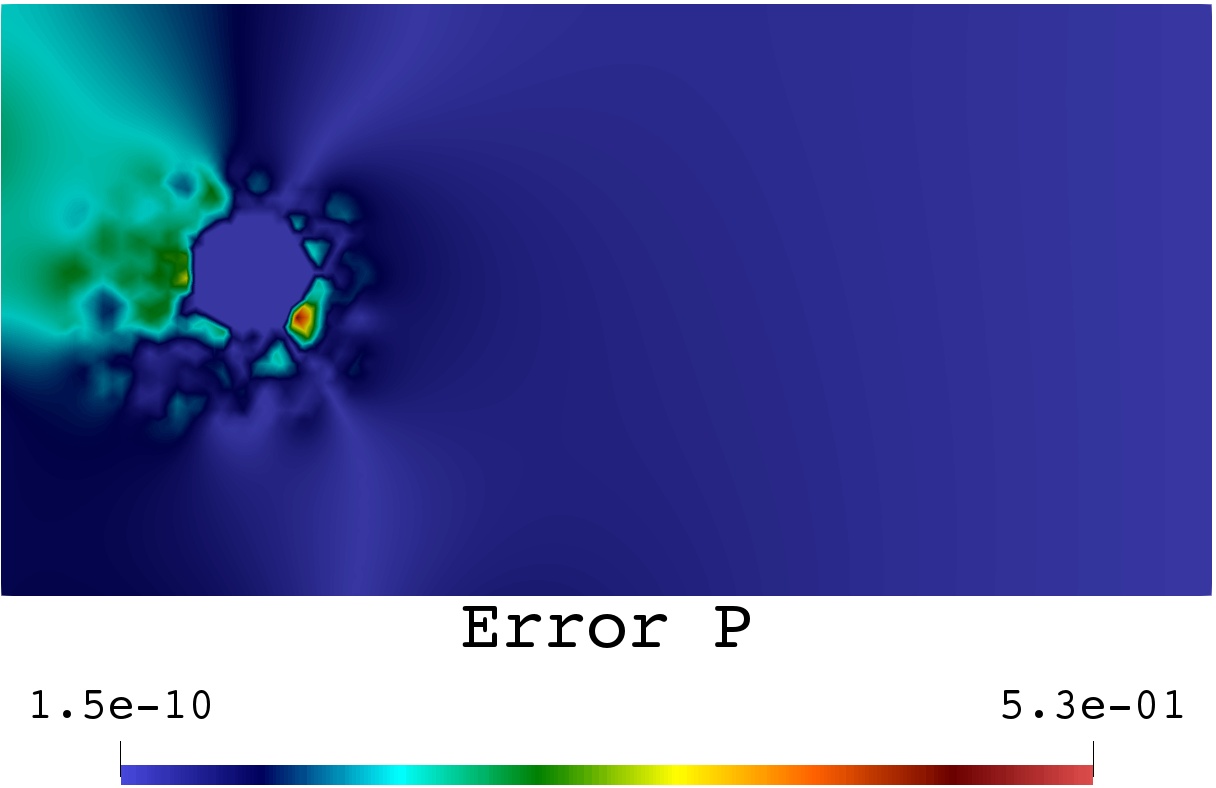}
\end{minipage}
\begin{minipage}{0.24\textwidth}
  \includegraphics[width=\textwidth]{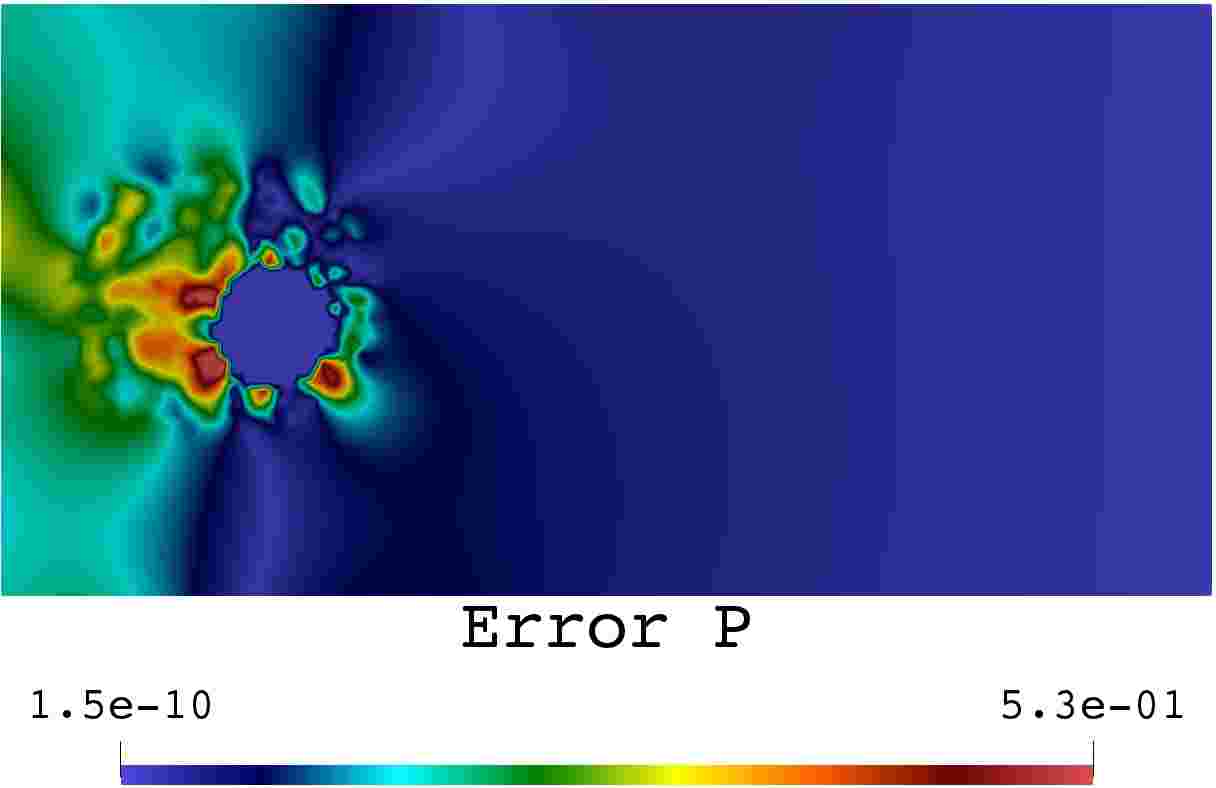}
\end{minipage}
\begin{minipage}{0.24\textwidth}
  \includegraphics[width=\textwidth]{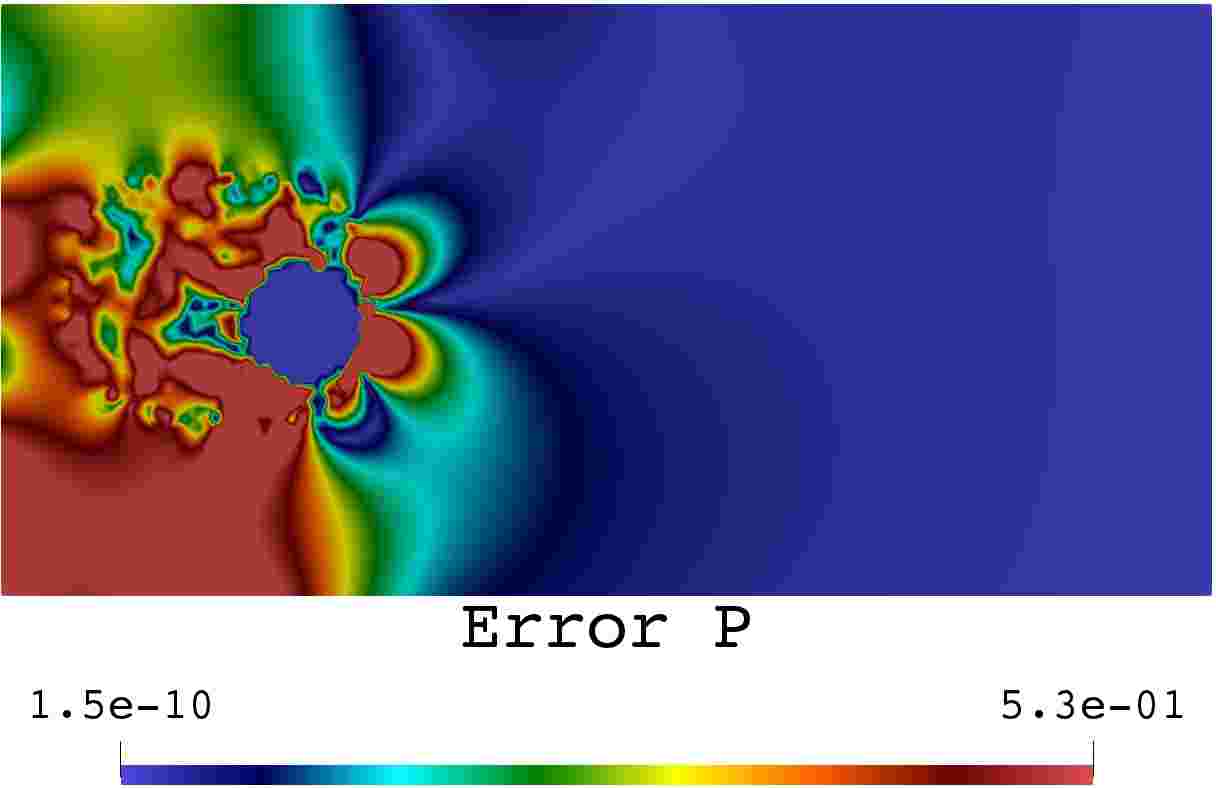}
\end{minipage}
\begin{minipage}{0.24\textwidth}
  \includegraphics[width=\textwidth]{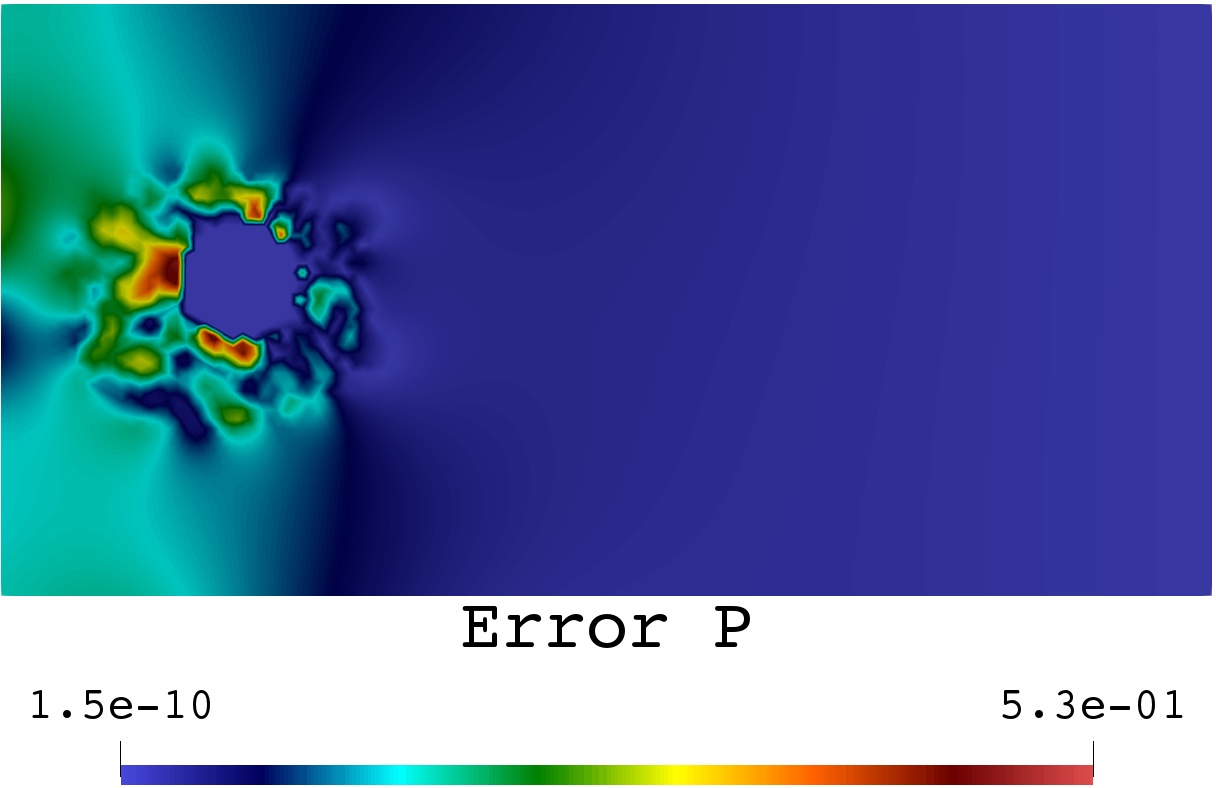}
\end{minipage}
\end{minipage}
\caption{Results for the 2D geometrical parametrization with $(\mu_0,\mu_1) \in [-1.5, -1.0]\times[-0.15, 0.15]$ {\blue and supremizer} basis enrichment. In rows $1-3$ we report the full-order solution, the reduced order solution and the absolute error plots for the velocity field while in rows $4-6$ we report the same quantities for the pressure field. The different columns are for four different values of the input parameter $\mu=[(-1.1917, 0.0958), (-1.1095, -0.0958), (-1.0273, -0.0684), (-1.2191,0.0684)]$.
}
  \label{FULL_RED_ERROR_2P}
\end{figure} 

\begin{table} \centering
  \begin{tabular}{ccccc}
   \hline
     Supremizers enrichment:& \multicolumn{2}{c}{No} & \multicolumn{2}{c}{Yes}   \\
    Number of modes & relative error u & relative error  p & relative error u & relative error  p \\
    \hline
     8  & 0.0947158 &  12.309881  &  0.2406999 &  22.319781\\
     12 & 0.0723268 &  12.133591  &  0.2078557 &  5.7159319\\
     16 & 0.0610052 &  9.6652163  &  0.1692787 &  2.6962056\\
     20 & 0.0538906 &  6.1692750  &  0.1243368 &  1.2535779\\
     25 & 0.0434925 &  3.2331644  &  0.0770726 &  0.5568314\\
     30 & 0.0396132 &  1.4693532  &  0.0437348 &  0.2504069\\
     35 & 0.0298269 &  0.7455038  &  0.0262345 &  0.1356788\\
     40 & 0.0177170 &  0.2918072  &  0.0121903 &  0.0611154\\
     45 & 0.0085905 &  0.0923509  &  0.0060355 &  0.0330206\\
     50 & 0.0053882 &  0.0473412  &  0.0046300 &  0.0279857\\
    \hline
  \end{tabular}
  \caption{Relative error between the full-order solution and the reduced basis solution for velocity and pressure in the case of the 1D geometrical parametrization. Results are reported for different dimensions of the reduced basis {\blue spaces with} and without supremizer stabilization.}
  \label{table:errors_no_supremizers}  %\label{table:errors_supremizers}
\end{table}

\begin{figure} \label{fig:comp_time} \centering
\begin{minipage}{\textwidth}
\centering
\begin{minipage}{0.48\textwidth}
\includegraphics[width=\textwidth]{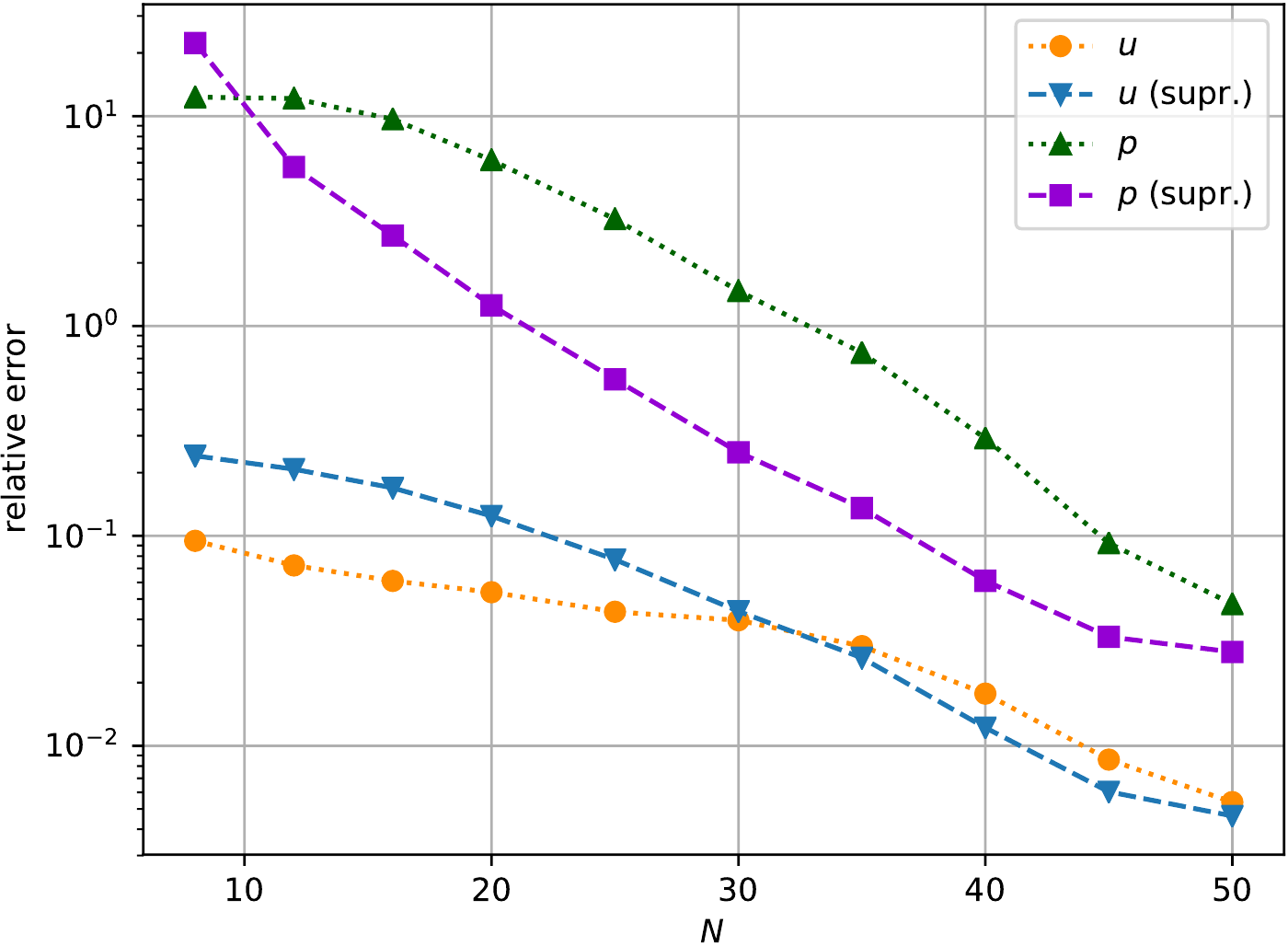}
\end{minipage}
\begin{minipage}{0.48\textwidth}
\includegraphics[width=\textwidth]{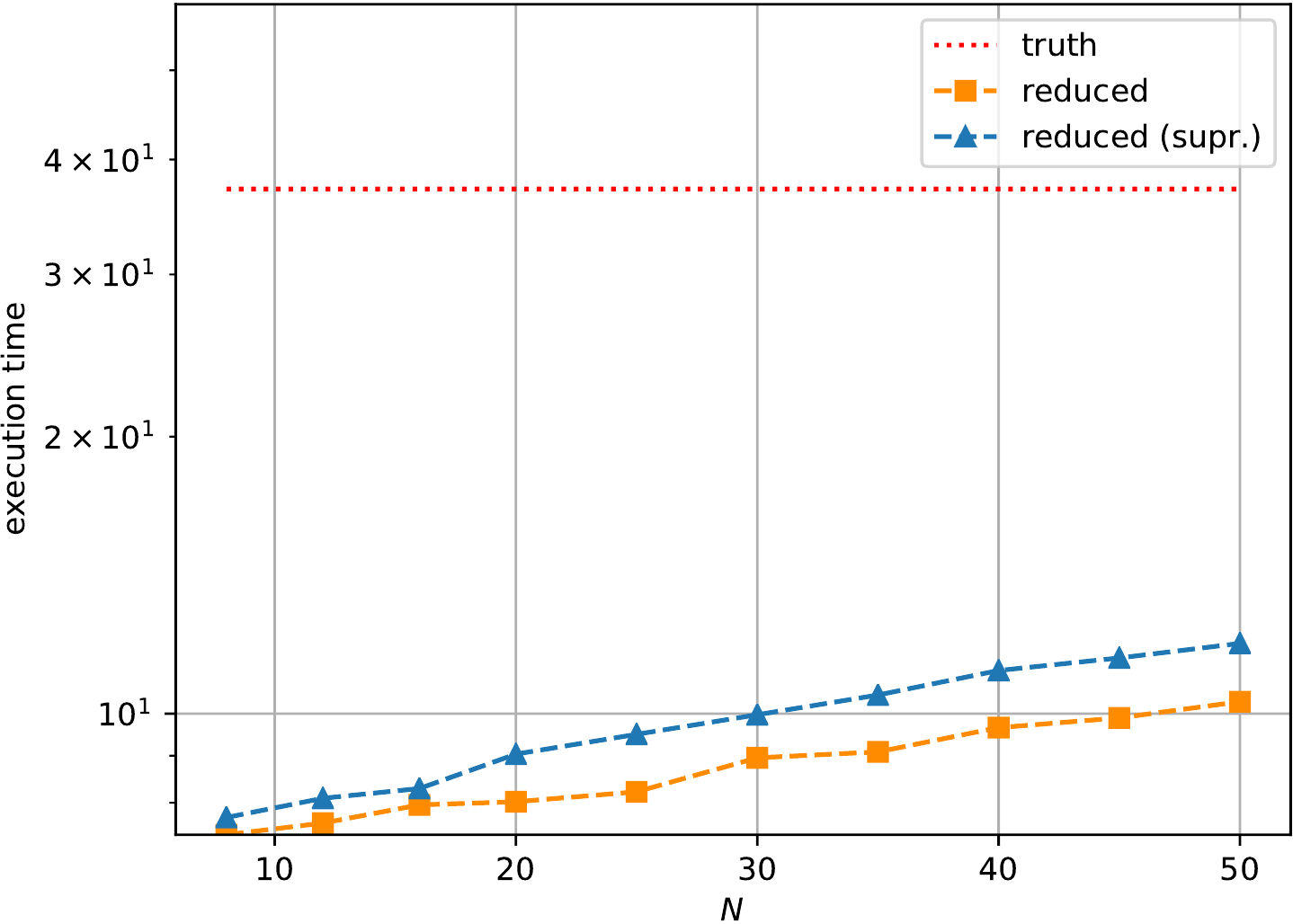}
\end{minipage}
\end{minipage}
  \caption{Visualization of the results for the case with 1D geometrical parametrization. On the left plot are depicted the relative errors for velocity and pressure with and without supremizer stabilization. On the right plot we report the execution times of the reduced order problem for 10 parameters values with and without supremizers enrichment. In both plots, results are reported for various number of modes.}
   \label{fig:timers} \label{fig:WithoughtSupremizer_solution_1000_100_1000_100_1024_j_3J_4j}\label{fig:Supremizer_solution_1000_100_1000_100_1024_j_3J_4j}
\end{figure}

\begin{table} \centering
  \begin{tabular}{ccccc}
    \hline
    Number of  modes & with no stabilization &  with stabilization   \\
                                      & execution time (sec) &  execution time (sec) \\
    \hline
     8   &  7.3858961  &   7.710907   \\
     12 &  7.6042165  &   8.091225  \\
     16 &  7.9584049  &   8.290780   \\
     20 &  8.0206915  &   9.036709    \\
     25 &  8.2229143  &   9.495323   \\
     30 &  8.9529275  &   9.972288  \\
     35 &  9.0867916  &   10.47633   \\
     40 &  9.6555775  &   11.13931    \\
     45 &  9.8934008  &   11.49422   \\
     50 &  10.302459  &   11.92024   \\
    \hline
  \end{tabular}
  \caption{Execution time, at the reduced order level, for the case with 1D geometrical parametrization. The computation time includes the assembling of the full-order matrices, their projection and the resolution of reduced problem. Times are for the resolution of 10 different values of the input parameter. The time execution at full-order level is equal to $\approx 37.16$ sec. }
  \label{table:timers}
\end{table}
In \autoref{fig:timers} (right plot) we report the total execution time for the resolution of the reduced order problem corresponding to 10 solves of the reduced order and full-order solvers respectively. One full-order solve takes $3.716$s. As one can see, the ROM leads to a considerable speed-up for { \blue all the different} analyzed configurations and for both cases with and without supremizer enrichment. 
\begin{remark}
It is important to recall that in this case, since the interest is into testing the feasibility and the accuracy of a reduced order model constructed starting from a shifted boundary method full-order solver, we did not employ any hyper reduction technique. This means that, also at the reduced order level, we assembled the full-order discretized differential operators. However, since most of the computational cost is required by the resolution of the discrete algebraic system rather than into its assembly, we can still achieve reasonable computational speedups. 
\end{remark}
\subsection{2D geometrical parametrization}
\begin{figure} \centering
\begin{minipage}{\textwidth}
\centering
\begin{minipage}{0.48\textwidth}
\includegraphics[width=\textwidth]{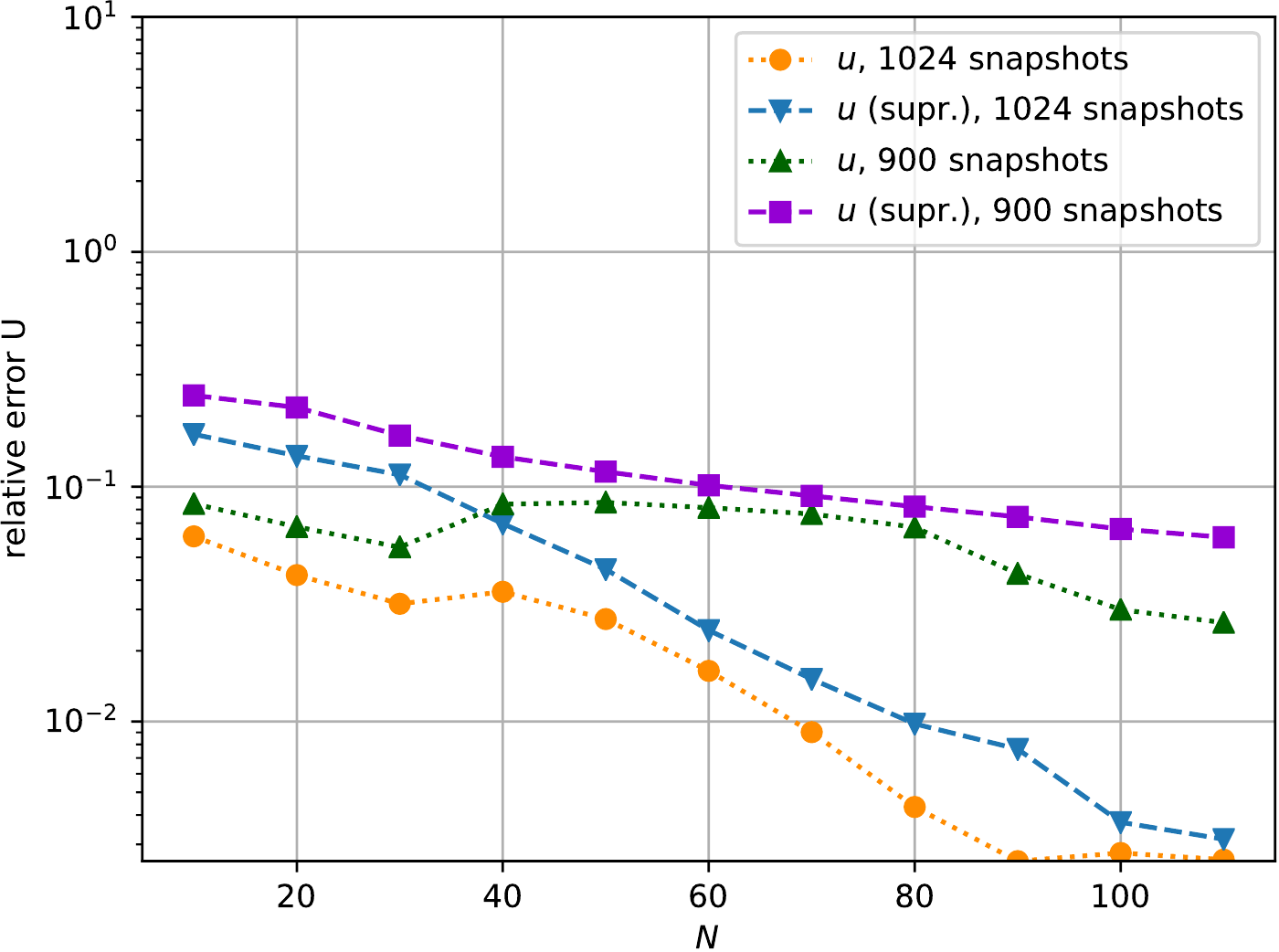}
\end{minipage}
\begin{minipage}{0.48\textwidth}
\includegraphics[width=\textwidth]{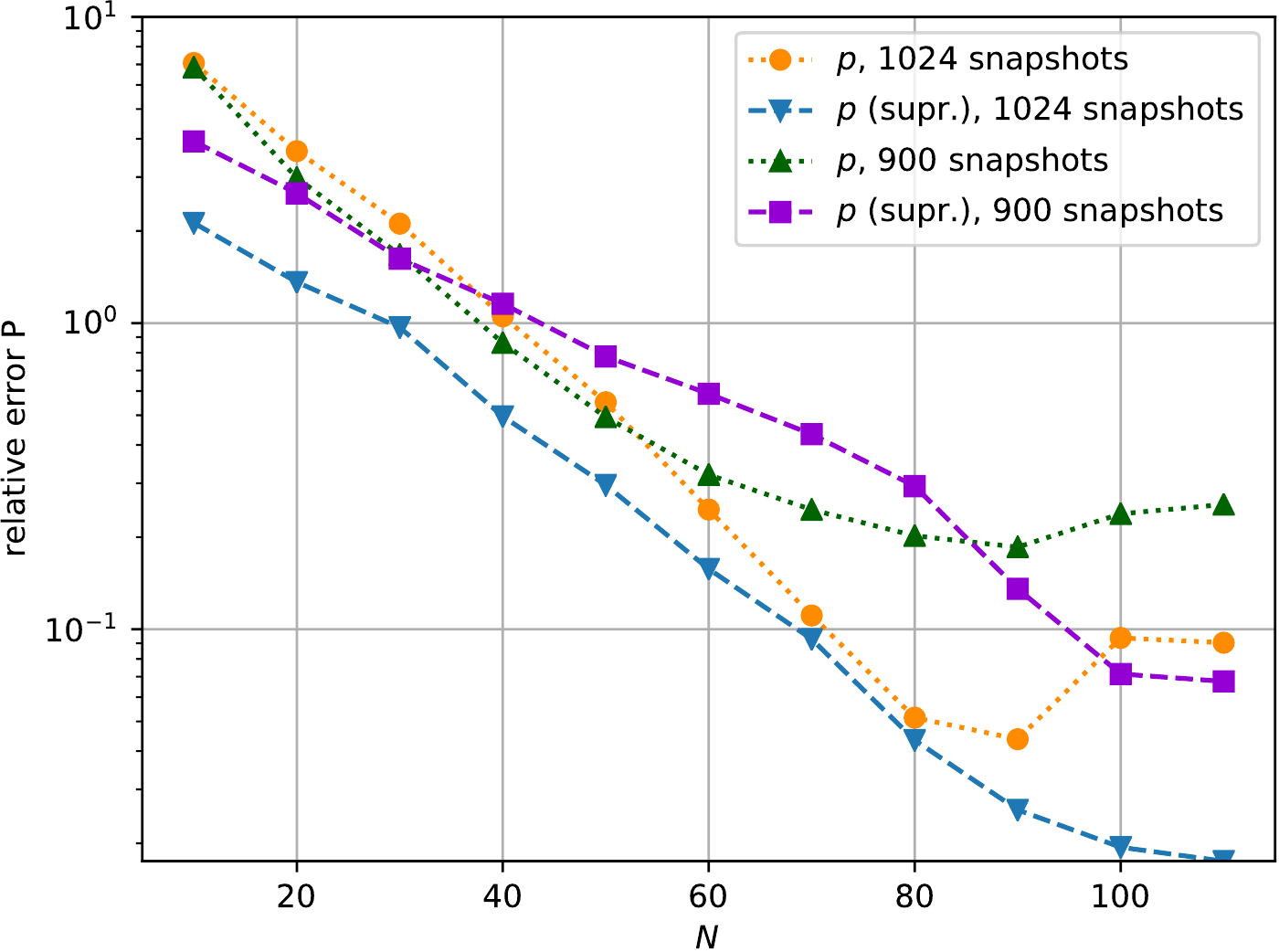}
\end{minipage}
\end{minipage}
  \caption{Visualization of the results for the case with 2D geometrical parametrization. On the left plot are depicted the relative errors for the velocity field with and without supremizer stabilization and for a different number of initial snapshots ($900$ and $1024$). On the right plot the same results are reported for the pressure field. In both plots results are reported for various number of modes.}
  \label{fig:2Dparametrization}
\end{figure}

The second case considers a 2D geometrical parametrization in the range  $\mu=(\mu_0,\mu_1)=[-1.5,-1.0]\times[-0.15,0.15]$. We perform this test to examine the performances of the methodology on a more complex and demanding scenario. It can be clearly seen in \autoref{fig:2Dparametrization} (right plot) that the supremizer enrichment is pointed out necessary for convergence stability and reliable pressure results. Finally, we notice that, without a supremizer enrichment, in this 2D geometrical parametrization case the best achieved relative error were equal to $0.0263821$ and $0.1854861$ for velocity and pressure respectively. The detailed results are omitted for shortness of space. In \autoref{table:2Dparametrization}, for { \blue the supremizer stabilization case}, we report the relative error for different dimensions of the reduced basis space and for different dimensions of the initial snapshots matrices ($900$ and $1024$ full order solutions respectively). This test serves to investigate how the number of snapshots used to compute the reduced basis space affects the accuracy of the results and to show that it is very important to use a large number of snapshots especially for reliable pressure results. The table clearly shows that a larger number of initial snapshots increase considerably the computational accuracy of the method. This is given by the obvious fact that an increase of the number of snapshots leads to a better representation of the solution manifold and therefore to a reduced basis space with better approximation properties. 
%As it is natural, we see differences and small instability in convergence for the experiments without supremizer enrichment (green and orange lines in Figure \ref{fig:2Dparametrization}), but in the cases with the pressure stabilization the differences in convergence are smaller (light blue and purple lines in Figure \ref{fig:2Dparametrization}), although not negligible, better for the more snapshots case and with faster and stabilized convergence.  
%

\begin{table} \centering
  \begin{tabular}{ccccc
}
   \hline
     Number of snapshots: & \multicolumn{2}{c}{900} & \multicolumn{2}{c}{1024}   \\
    Number of modes & relative error u & relative error  p & relative error u & relative error  p 
\\
    \hline
     10 & 0.2448511 &  3.9240637   &   0.1672753   &   2.1243228 
\\
     20 & 0.2175821 &  2.6531343   &   0.1353706   &   1.3611011 
\\
     30 & 0.1652331 &  1.6234701   &   0.1124619   &   0.9680506 
\\
     40 & 0.1340978 &  1.1560352   &   0.0696437   &   0.4958605 
\\
     50 & 0.1158443 &  0.7777786   &   0.0444991   &   0.2958338 
\\
     60 & 0.1013961 &  0.5876048   &   0.0244793   &   0.1574037 
\\
     70 & 0.0914650 &  0.4335489  &    0.0151749   &   0.0928402 
\\
     80 & 0.0822658 &  0.2933336  &    0.0097848   &   0.0434299
\\
     90 & 0.0744696 &  0.1355488   &   0.0076431   &   0.0257060
\\
     100& 0.0660493 &  0.0714350   &   0.0037280   &   0.0194051  
\\
     110& 0.0609040 &  0.0675720 &     0.0031577   &   0.0174815
\\
    \hline
  \end{tabular}
  \caption{Supremizer basis enrichment and the relative error between full-order solution and reduced basis solution for velocity and pressure for two different numbers of snapshots in the case with a 2D geometrical parametrization.}
  \label{table:2Dparametrization}  %\label{table:errors_supremizers}
\end{table}

{\blue
\subsection{Free-Form Deformation Geometrical Parametrization}

In this last numerical example we apply the developed methodology to deal with a more complex
geometrical parametrization. The Free-Form Deformation (FFD) tool is used to morph the original undeformed geometry \cite{sederbergparry1986}. The Free-Form Deformation is particularly useful in the context of parametrized reduced order model \cite{LassilaRozza2010} and has been used to deal with complex parametrized domains in both automotive and naval engineering \cite{salmoiraghi2018,tezzele2018dimension}. In this numerical example, since we wanted to investigate the applicability of the proposed method to more complex geometrical transformation we start from a simple initial shape and we deform it using free-form deformation. 

The undeformed geometry is given by the sketch reported in the left side of Figure 
\ref{fig:sketch_ffd}. The embedded boundary consists into a circular and a square boundary. The circular domain is deformed using an FFD approach following the sketch reported in the right side of Figure~\ref{fig:sketch_ffd}. The free form deformation is based on the definition of a bounding box identified by a set of properly chosen control points. For the particular case the movement of two control points has been parametrized by the parameter vector $ \mu = (\mu_1, \mu_2)$ that describes the horizontal displacement of control points located on the left of the bounding box (see Figure \ref{fig:sketch_ffd}). In Figure \ref{fig:error_ffd} we report the convergence analysis for the numerical example. The test has been conducted using an equal number of modes for velocity, pressure and supremizer space ($N^r_u = N^r_p = N^r_s = 10j$). In Figure \ref{fig:geometry_snapshots_ffd} we report the deformed geometry for different samples inside the parameter space. The reduced order model has been trained using $1024$ samples chosen randomly inside the parameter space $\mathcal{P}_{\mathrm{train}} \in [-0.6,0.6]^2$ and tested using $100$ samples inside the parameter space $\mathcal{P}_{\mathrm{test}} \in [-0.5,0.5]^2$.
 
\begin{figure} 
\centering
\begin{minipage}{\textwidth}
\centering
\begin{minipage}{0.45\textwidth}
\includegraphics[width=\textwidth]{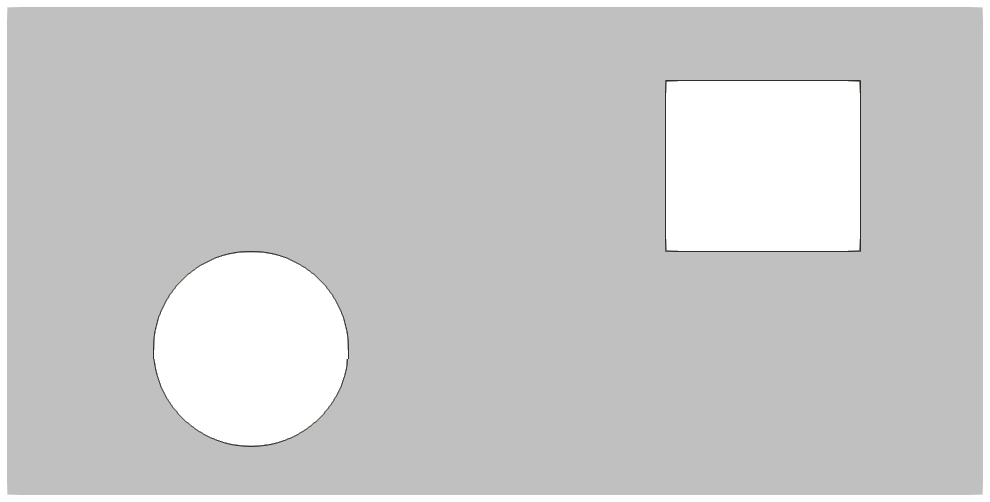}
\end{minipage}
\begin{minipage}{0.45\textwidth}
\includegraphics[width=\textwidth]{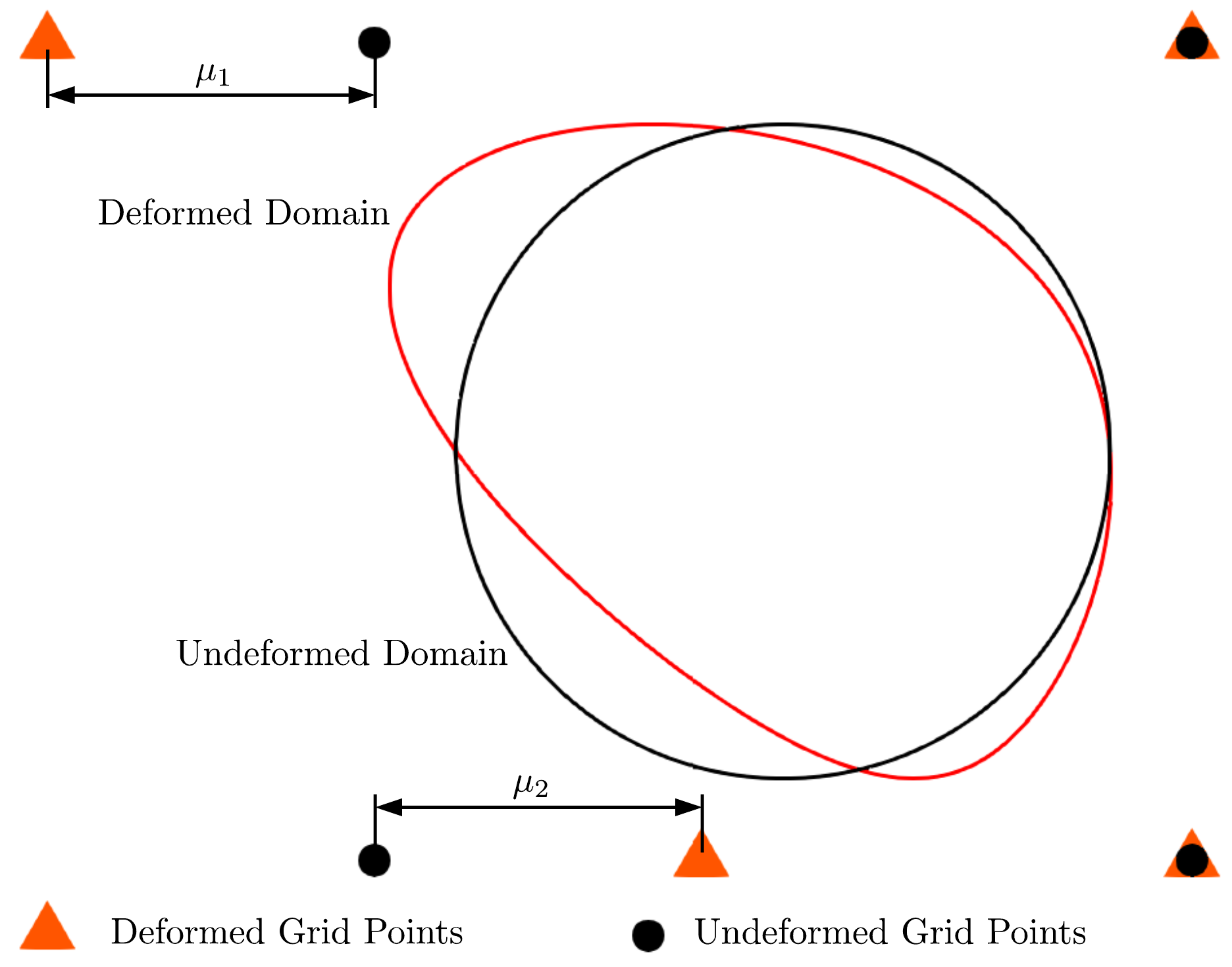}
\end{minipage}
\end{minipage}
\caption{{\blue Sketch of the undeformed embedded domain and parametrization used to morph the circular domain}} \label{fig:sketch_ffd} 
\end{figure}

\begin{figure} 
\centering
\begin{minipage}{\textwidth}
\centering
\begin{minipage}{0.245\textwidth}
\includegraphics[width=\textwidth]{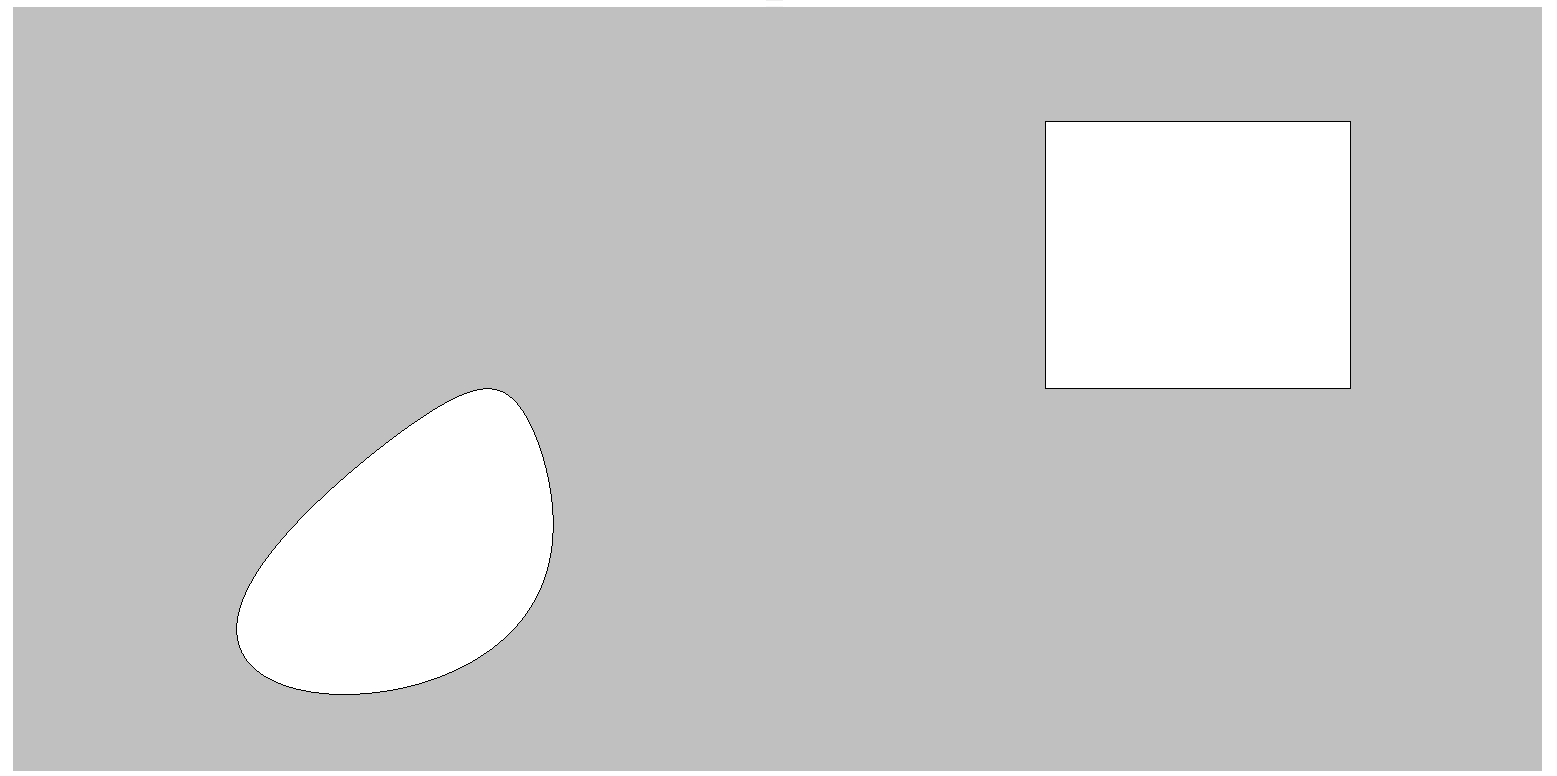}
\end{minipage}
\begin{minipage}{0.245\textwidth}
\includegraphics[width=\textwidth]{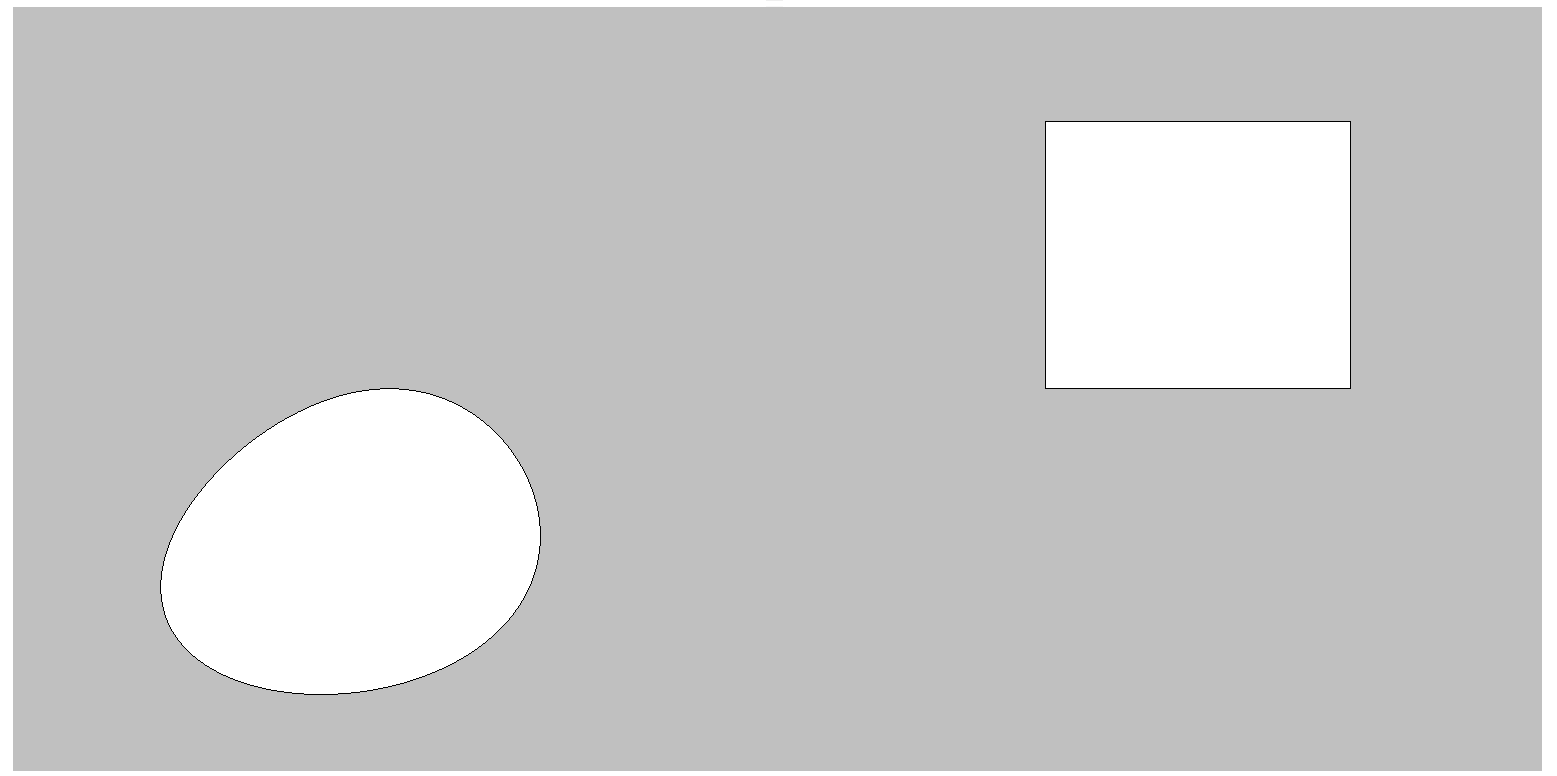}
\end{minipage}
\begin{minipage}{0.245\textwidth}
\includegraphics[width=\textwidth]{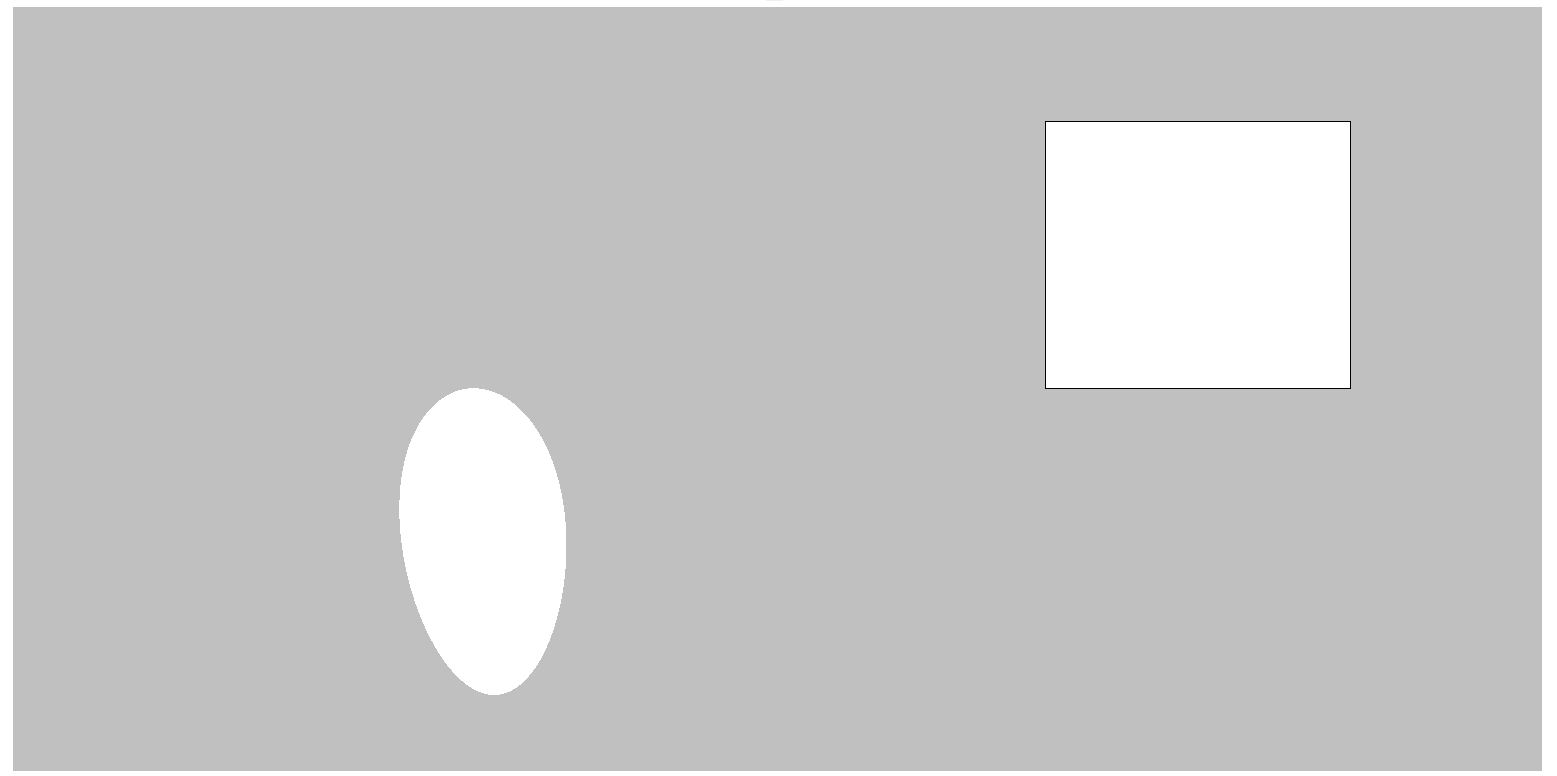}
\end{minipage}
\begin{minipage}{0.245\textwidth}
\includegraphics[width=\textwidth]{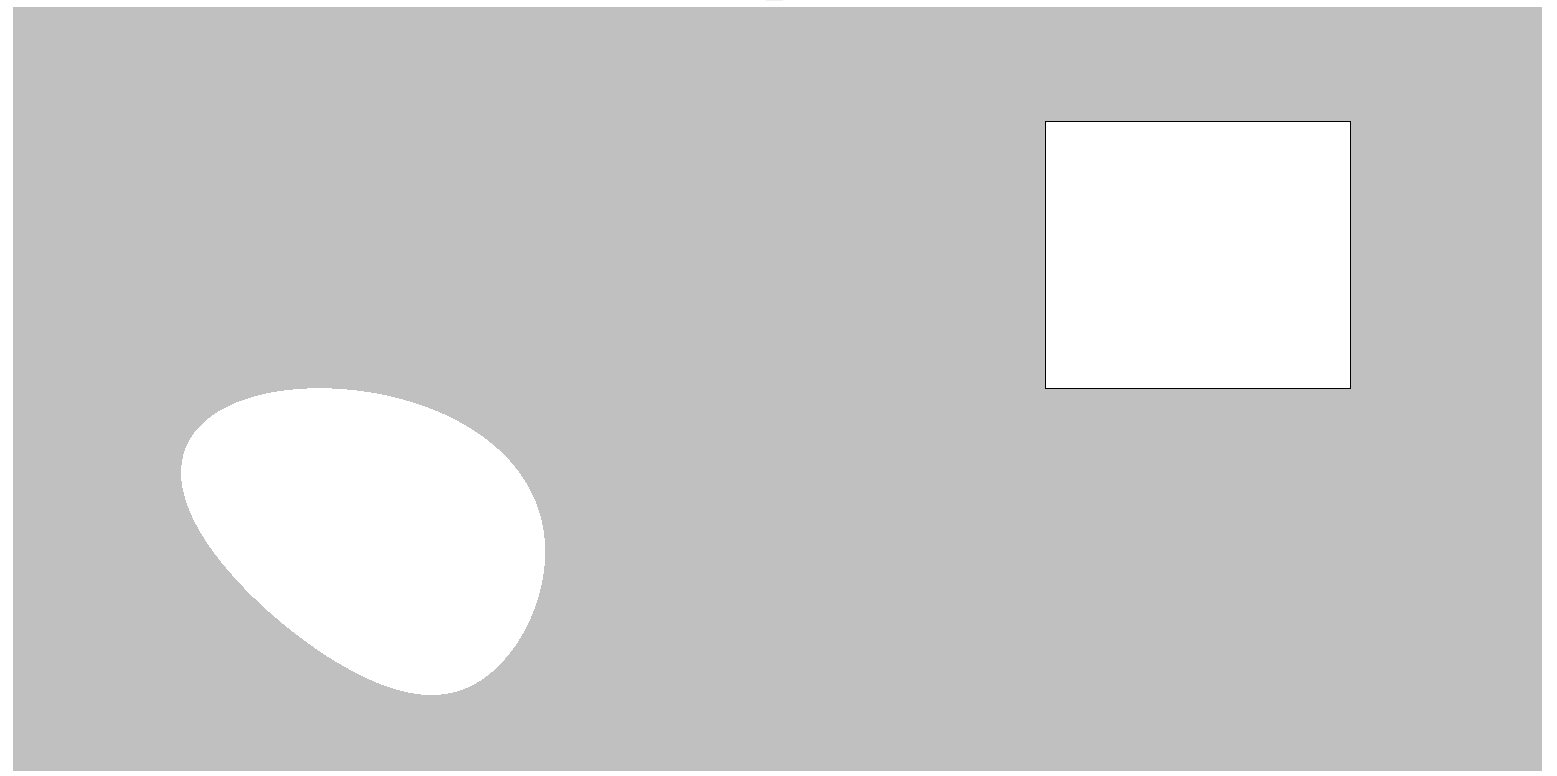}
\end{minipage}
\end{minipage}
\caption{{\blue Free-Form deformation geometrical parametrization and visualization of four samples $\mu=[(0.2613,-0.4698),(0.2613,-0.4698),(0.2613,-0.4698),(0.2613,-0.4698)]$}}
\label{fig:geometry_snapshots_ffd}
\end{figure}

\begin{figure} \centering
\begin{minipage}{\textwidth}
\centering
\begin{minipage}{0.48\textwidth}
\includegraphics[width=\textwidth]{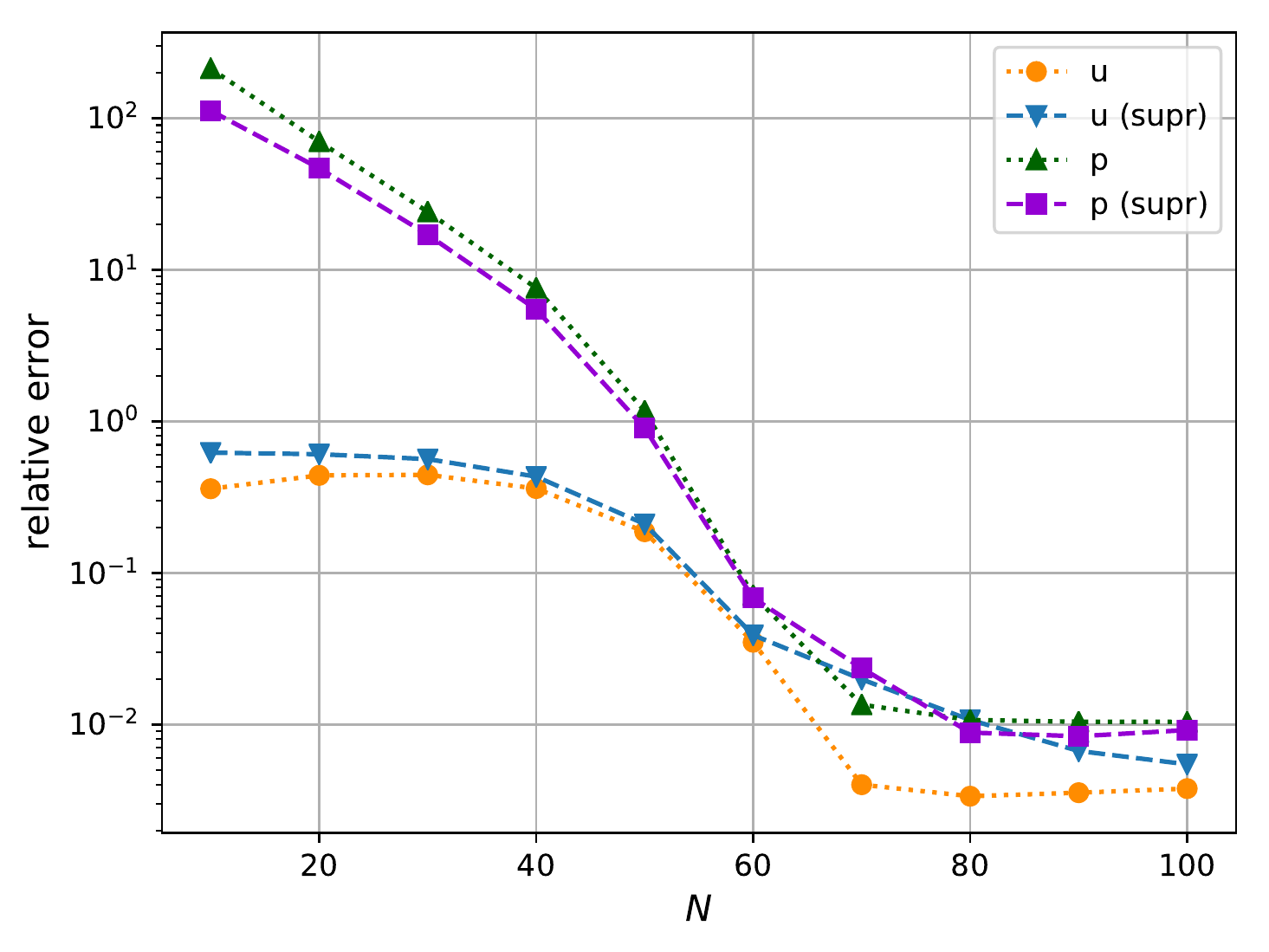}
\end{minipage}
\end{minipage}
\caption{Results for the numerical example using Free-Form deformation. The relative errors for velocity and pressure with and without supremizer stabilization are depicted for various number of modes.}\label{fig:error_ffd}
\end{figure}
}

\subsection{Some Comments}
The numerical experiments clearly show that a projection based reduced order model can be generated on top of a full-order SBM solver. The important aspect is that, especially for complex geometrical parametrization, one can avoid the somehow ``expensive'' traditional approach of the geometrical transformation to a reference domain and to only rely on a suitable smooth enough background mesh.  
The experiments have shown that it is possible to obtain accurate results for velocity and pressure even without a supremizer ``inf-sup''{\blue{-stabilization}} at the reduced level. 

After comparing the results in \autoref{fig:2Dparametrization} and \autoref{fig:Supremizer_solution_1000_100_1000_100_1024_j_3J_4j} one can observe that the supremizer stabilization permits to considerably increase the accuracy of the ROM for the pressure field without deteriorating the velocity accuracy. In the case of 1D geometrical parametrization it even permits to increase the accuracy for both the velocity and the pressure field. 

\section{Conclusions and future developments}\label{sec:conclusions}
In the present work a POD-Galerkin ROM based on SBM full-order simulations was presented. The ROM was developed to be consistent with the full-order model, and both velocity and pressure fields were considered.  The reduced order method is applied to approximate the geometrically parametrized Stokes flow around a circular embedded domain.
%{\red (additional preliminary tests on a Poisson problem are provided in \ref{sec:Preliminary_results}).} 
In particular, the focus and originality of the present work stands on the application of a velocity-pressure ROM using a background mesh in the EBM context.
A comparison of accuracy was made for simulations with and without pressure stabilization strategies for POD-Galerkin ROMs in combination with the SBM. 

The proposed ROM employs a stabilization strategy which is applied in the full order discretization formula. 
This stabilization proved to be efficient and effective also for the resulting reduced system for both the velocity and pressure fields. However, best accuracy for the pressure fields was achieved by further stabilizing the reduced system using a supremizer stabilization strategy.
The attention was in fact devoted to analyze the applicability of the proposed methods to flows. 
The ROM demonstrated to be capable of reproducing all the main features of the physical phenomenon in an accurate manner leading to a considerable computational reduction. 

Future developments will focus on different efficient methodologies for the affine decomposition of the SBM differential operator and in particular to study the applicability of well-known hyper reduction techniques, such as the empirical interpolation method, in the context of an SBM framework. 
Of future interest are also fluid-structure interaction applications as the ones presented in \cite{Stabile2018} and Navier-Stokes problems. 

\section*{Acknowledgments}
We acknowledge Dr Andrea Mola and Dr Francesco Ballarin from SISSA at International School for Advanced Studies, Mathematics Area, for fruitful discussions related to software implementation and RB methods, respectively. 
%{\red, and the PhD candidate Nabil Atallah from Department of Civil and Environmental Engineering at Duke University for the help and the information provided in the SBM Poisson case.}
This work is supported by the U.S. Department of Energy, Office of Science, Advanced Scientific Computing Research under Early Career Research Program Grant SC0012169, the U.S. Office of Naval Research under grant N00014-14-1-0311,  ExxonMobil Upstream Research Company (Houston, TX), the European Research Council Executive Agency by means of the H2020 ERC Consolidator Grant project AROMA-CFD ``Advanced Reduced Order Methods  with  Applications  in  Computational  Fluid  Dynamics'' - GA  681447, (PI: Prof. G. Rozza), INdAM-GNCS 2018 and by project FSE - European Social Fund - HEaD "Higher Education and Development" SISSA operazione 1, Regione Autonoma Friuli - Venezia Giulia.
%
%\begin{figure} 
%\centering
%\begin{minipage}{0.48\textwidth}
%\centering
%\includegraphics[width=\textwidth]{img_poisson-crop.pdf}
%\footnotesize
%(a)
%\end{minipage}
%\begin{minipage}{0.48\textwidth}
%\centering
%\includegraphics[width=\textwidth]{img/poisson_eig-crop.pdf}
%\footnotesize
%(b)
%\end{minipage}
%\caption{\red Poisson problem: The figure reports, on the left, (a) the mean relative error for the testing L2 projection of the full-order solution onto the POD basis functions, and (b) the mean relative error for the Galerkin projection ROM for various values of the number of modes. On the right, the eigenvalue decay of POD procedure.}
%\label{fig:poisson_results}
%\end{figure}
%%
%\begin{figure}
%\begin{minipage}{\textwidth}
%\centering
%\begin{minipage}{0.32\textwidth}
%\includegraphics[width=\textwidth]{img/T.jpg}
%\end{minipage} 
%\begin{minipage}{0.32\textwidth}
%\includegraphics[width=\textwidth]{img/T_red.jpg}
%\end{minipage} 
%\begin{minipage}{0.32\textwidth}
%\includegraphics[width=\textwidth]{img/errorT.jpg}
%\end{minipage} 
%\end{minipage}
%\caption{\red Poisson problem: From left to right we report the full-order results, the reduced order results and the absolute value of the error, respectively. The results are for one random value of the parameter in the range of variation.}\label{fig:poisson_field}
%\end{figure}
%}

\newpage

%\section*{References}
\bibliographystyle{amsplain_gio} % \section*{References}
\bibliography{bibfile_sissa}
\end{document}